\documentclass[12pt,oneside,british,a4wide]{amsart}
\usepackage{amsmath, amssymb,verbatim,appendix}
\usepackage[pagebackref=true]{hyperref}
\usepackage[mathscr]{eucal}
\usepackage{amscd}
\usepackage{amsthm}
\usepackage{stmaryrd}
\usepackage{enumitem}
\usepackage{comment}
\usepackage[latin1]{inputenc} 
\usepackage{tikz}
\usetikzlibrary{shapes,arrows}
\usepackage{cite}
\usepackage{url}
\usepackage{float}

\usepackage{setspace,a4wide,color,xcolor,graphicx}
\def\showauthornotes{1}

\ifnum\showauthornotes=1
\newcommand{\Authornote}[2]{{\sf\small\color{red}{[#1: #2]}}}
\else
\newcommand{\Authornote}[2]{}
\fi

\newtheorem{theorem}{Theorem}[section]
\newtheorem{lemma}[theorem]{Lemma}
\newtheorem{corollary}[theorem]{Corollary}
\newtheorem{proposition}[theorem]{Proposition}

\newtheorem{question}[theorem]{Question}

\newtheorem{fact}[theorem]{Fact}

\numberwithin{equation}{section}
\setlist[enumerate]{font={\rmfamily}}
\setlist[enumerate,1]{label={(\roman*)}}

\theoremstyle{definition}

\newtheorem{definition}[theorem]{Definition}

\newtheorem{notation}[theorem]{Notation}

\def\E{\mathbb{E}}

\def\R{\mathbb{R}}

\def\C{\mathbb{C}}

\def\F{\mathbb{F}}

\def\OP{\mathrm{OP}}

\def\NFOP{\mathrm{NFOP}}
\def\IP{\mathrm{IP}}
\def\NIP{\mathrm{NIP}}

\def\VC{\mathrm{VC}}

\newcommand{\e}{\varepsilon}

\def\ra{\rightarrow}

\newcommand{\dbar}{\bar{d}}

\newcommand{\calL}{\mathcal{L}}

\newcommand{\calC}{\mathcal{C}}

\newcommand{\calB}{\mathcal{B}}

\newcommand{\calP}{\mathcal{P}}

\newcommand{\calI}{\mathcal{I}}

\newcommand{\calQ}{\mathcal{Q}}

\newcommand{\calS}{\mathcal{S}}

\newcommand{\At}{{\mathrm{At}}}
\newcommand{\rk}{{\mathrm{rk}}}

\newenvironment{proofof}[1]{\indent{\itshape Proof of #1}.\;}{\qed}

\begin{document}

\title[The structure of subsets of $\F_p^n$ of bounded $\mathrm{VC_{2}}$-dimension]{The structure of subsets of $\F_p^n$\\ of bounded $\mathrm{VC_{2}}$-dimension}

\author{C. Terry}

\author{J. Wolf}

\address{Department of Mathematics, Statistics, and Computer Science, University of Illinois Chicago, Chicago IL 60607, USA}

\email{caterry@uic.edu}

\address{Department of Pure Mathematics and Mathematical Statistics, Centre for Mathematical Sciences, Wilberforce Road, Cambridge CB3 0WB, UK}

\email{julia.wolf@dpmms.cam.ac.uk}

\date{}

\begin{abstract}We show that a subset of $\F_{p}^{n}$ of $\mathrm{VC_{2}}$-dimension at most $k$ is well approximated by a union of atoms of a quadratic factor of complexity $(\ell,q)$ (denoting the complexities of the linear and quadratic part, respectively), where $\ell$ and $q$ are bounded by a constant depending only on $k$ and the desired level of approximation. This generalises a result of Alon, Fox and Zhao \cite{Alon.2018is} on the structure of sets of bounded $\VC$-dimension, and is analogous to contemporaneous work of the authors \cite{Terry.2021b} in the setting of 3-uniform hypergraphs.  

The main result originally appeared--albeit with a different proof--in a 2021 preprint \cite{Terry.2021av2}, which has since been split into two: the present work, which focuses on higher arity NIP and develops a theory of local uniformity semi-norms of possibly independent interest, and its companion \cite{Terry.2021a}, which strengthens these results under a generalized notion of stability.
\end{abstract}

\maketitle

\tableofcontents

\section{Introduction}

A major theme in combinatorics is the decomposition of a mathematical object (say a graph, set, or function) into a structured part and a random-looking (or \emph{pseudorandom}) part. By discarding the pseudorandom contributions and explicitly computing the contributions from the structured part, it is often possible to obtain a count of certain substructures inside the object in question. This philosophy has its roots in techniques in classical analytic number theory, and came to the fore in graph theory with Szemer\'edi's regularity lemma \cite{Szemeredi.1975} in the 1970s. The key is to define a suitable notion of pseudorandomness for the problem at hand, and then to obtain a structural counterpart. While the latter is usually regarded as the deep and difficult part of the problem, the former is evidently crucial to the success of any decomposition endeavour.

In many famous instances, such as Roth's theorem on 3-term arithmetic progressions, the $\ell^\infty$ norm of the Fourier transform provides a suitable measure of pseudorandomness. Indeed, any bounded function on a finite abelian group can be decomposed into a part that is pseudorandom in this sense, plus a structured part, which is determined by a small number of linear complex exponential phases (associated with the large Fourier coefficients of the function). 

However, it was (re)discovered by Gowers in his groundbreaking work on Szemer\'edi's theorem \cite{Gowers.1998, Gowers.2001} that the $\ell^\infty$ norm of the Fourier transform is not sufficiently sensitive for controlling, for example, arithmetic progressions of length 4. Gowers introduced a new measure of pseudorandomness, the so-called $U^3$ norm, which turns out to be more suitable for this purpose. He also proved a structure result for functions with large $U^3$ norm, showing that they correlate with suitably defined quadratic (complex exponential) phases. This opened the door to so-called \emph{quadratic regularity lemmas}, which break down a function into a quadratically pseudorandom part (as measured by the $U^3$ norm) and a part determined by a small number of quadratic phases. 

This \emph{higher-order Fourier analysis} has developed into an exceedingly active area of research over the past two decades, leading to a number of high-profile applications such as the Green-Tao theorem on long arithmetic progressions in the primes \cite{Green.2008}. It has also engendered numerous fruitful connections to other subdisciplines, including ergodic theory and theoretical computer science.

In applications, the quantitative aspects of such decompositions are often extremely important, and it is natural to ask whether there are circumstances where particularly efficient decompositions can be guaranteed. In the context of linear regularity decompositions of subsets of $\F_p^n$, the authors showed in \cite{Terry.2019} that the model-theoretic notion of \emph{stability} provides a sufficient condition for especially efficient decompositions. In subsequent work  \cite{Alon.2018is}, Alon, Fox and Zhao showed that a weaker assumption regarding the \emph{$\VC$-dimension} alone produces an efficient bound on the complexity of a linear regularity decomposition. We begin our detailed account of these and related results by recalling the definition of $\VC$-dimension for subsets of an abelian group. 

\begin{definition}[$\VC$-dimension of a subset of a group]\label{def:vcdim}
Given $k\geq 1$, an abelian group $G$ and a subset $A\subseteq G$, we say that $A$ \emph{has $\VC$-dimension at least $k$} (or \emph{has $k$-$\IP$})\footnote{$\IP$ stands for the ``independence property'' in model theory.} if there exist elements $\{a_i: i\in [k]\}\cup \{b_S:S\subseteq [k]\}$ in $G$ such that $a_i+ b_S\in A$ if and only if $i\in S$.  

The \emph{$\VC$-dimension of $A$}, denoted $\VC(A)$, is defined to be the largest $k$ such that $A$ has $\VC$-dimension at least $k$. When $\VC(A)<k$, we also say that \emph{$A$ is $k$-$\NIP$}.
\end{definition}

Throughout this paper, we shall focus on the finite abelian group $\mathbb{F}_p^n$, where the prime characteristic $p$ of the base field is fixed and the dimension $n$ of the vector space is  thought of as large. This group frequently serves as a sandbox environment in which certain analytic arguments in additive combinatorics work particularly neatly. For examples of this, we refer the reader to the surveys \cite{Green.2005s, Wolf.2015, Peluse.2024}. 

In this context, the following arithmetic regularity lemma for sets of bounded $\VC$-dimension was proved by Alon, Fox and Zhao in \cite{Alon.2018is}.

\begin{theorem}[Arithmetic regularity lemma for subsets of bounded $\VC$-dimension]\label{thm:vc}
For all primes $p$, integers $k\geq 1$, and reals $\e>0$, there exists $M=M(p,k,\e)$ such that the following holds. 

Let $n\geq 1$ be an integer. Suppose that $A\subseteq \F_p^n$ satisfies $\VC(A)<k$. Then for all $\e>0$, there exists a subgroup $H\leq \F_p^n$ of index $m\leq M$ and a set $\Xi \subseteq G/H$ such that $|\Xi|\leq \e |\F_p^n|$ and such that for all $C\notin \Xi$,  $|A\cap C|/|C|\in [0,\e)\cup (1-\e,1]$.
\end{theorem}

In fact, Alon, Fox, and Zhao showed that the constant $M$ in Theorem \ref{thm:vc} can be taken to be of the form $\e^{-k-o(1)}$, where $o(1)$ is a term that tends to $0$ as $\e$ tends to zero, at a rate possibly depending on $p$ and $k$. Their results also apply to the more general setting of abelian groups of bounded exponent.

It is not difficult to see that a subgroup of any group $G$ has $\VC$-dimension at most 1 (see e.g. \cite[Example 1]{Terry.2019}).  It is also not hard to show that this implies that a union of $k$ cosets of a subgroup has  $\VC$-dimension at most $k$.  The following argument is a straightforward adaptation of \cite[Lemma 1.5]{Sanders.2018}. 

\begin{lemma}[Unions of cosets have bounded $\VC$-dimension]
Let $G$ be an abelian group, let $H\leqslant G$ and let $Y$ be a union of $k$ cosets of $H$. Then $\VC(Y)\leq \lceil\log_2 k\rceil$.
\end{lemma}
\begin{proof} Let $s=\lceil\log_2 k\rceil$, and suppose that $a_1,\dots, a_{s+1}$, $b_S \in G$ for $S\subseteq [s+1]$ witness the $(s+1)$-independence property, that is, $a_i+b_S\in Y$ if and only if $i\in S$. By the pigeonhole principle, there are distinct sets $S,S'\subseteq [s]$ such that $a_{s+1}+b_{S\cup\{s+1\}}$ and $a_{s+1}+b_{S'\cup\{s+1\}}$ lie in the same coset of $H$, whence $b_{S\cup\{s+1\}}+H=b_{S'\cup\{s+1\}}+H$. Without loss of generality, suppose that $j\in S'\setminus S$. Then $a_{j}+b_{S\cup\{s+1\}}\in a_j+b_{S'\cup\{s+1\}}+H \subseteq Y$, a contradiction.
\end{proof}

Theorem \ref{thm:vc} implies that,  in an approximate sense, unions of subgroups are the only examples of subsets of bounded $\VC$-dimension in $\F_p^n$. Indeed, the following statement is equivalent to Theorem \ref{thm:vc}.

\begin{corollary}[Structure theorem for sets of bounded $\VC$-dimension]\label{cor:vc}
For all primes $p$, integers $k\geq 1$ and reals $\e>0$, there exists $M=M(p,k,\e)$ such that the following holds. 

Let $n\geq 1$ be an integer. Suppose that $A\subseteq \F_p^n$ satisfies $\VC(A)<k$.  Then there exists a subgroup $H\leqslant \F_p^n$ of index $m\leq M$ and a union $Y$ of cosets of $H$ such that $|A\Delta Y|\leq \e|G|$. 
\end{corollary}

Theorem \ref{thm:vc} and Corollary \ref{cor:vc} also hold in arbitrary finite groups of bounded exponent (see \cite{Alon.2018is} for the abelian case, and \cite{Conant.2018zd} for the non-abelian), although no quantitative proof exists in the non-abelian setting. For results of a similar flavour outside the setting of bounded exponent, see work of Sisask \cite{Sisask.2018} and of Conant, Pillay, and the first author \cite{Conant.2018zd}.  

We will informally refer to the cosets in $\Xi$ in Theorem \ref{thm:vc} as the ``error cosets".  Indeed, the authors have shown, using an example of Green and Sanders \cite{Green.2015}, that under the assumption of bounded $\VC$-dimension alone, it is not possible to exclude the existence of such error cosets  (see \cite[Example 4]{Terry.2019} and \cite[Proposition 3.24]{Terry.2021a}/\cite[Proposition 5.18]{Terry.2021av2}). However, the stronger assumption of \emph{stability}, as defined below, guarantees that there is a decomposition for which the set $\Xi$ is the empty set. 

\begin{definition}[Stability of a subset of a group]\label{def:stable}
Given $k\geq 1$, an abelian group $G$ and a subset $A\subseteq G$, we say that $A$ \emph{has the $k$-order property} (or \emph{has $k$-$\OP$}) if there exist $\{a_i: i\in [k]\}\cup \{b_j:j\in [k]\}$ such that $a_i+ b_j\in A$ if and only if $i\leq j$.  

When $A$ does not have the $k$-order property in $G$, then we say $A$ is \emph{$k$-stable}. 
\end{definition}

The main result of \cite{Terry.2019} was the following.

\begin{theorem}[Arithmetic regularity lemma for stable sets]\label{thm:stablegps}
For all primes $p$, integers $k\geq 1$ and reals $\e>0$, there exist $N=N(p,k,\e)$ and $M=M(p,k,\e)$ such that the following holds. 
Let $n\geq N$ be an integer. Suppose that $A\subseteq \F_p^n$ is $k$-stable. Then there exists $H\leq \F_p^n$ of codimension $m\leq M$, such that for all cosets $C$ of $H$, $|A\cap C|/|C|\in [0,\e)\cup (1-\e,1]$.
\end{theorem}

In \cite{Terry.2019}, it was shown that the bound $M$ on the codimension can be taken to be a polynomial in $\e^{-1}$. Subsequent work of Conant \cite{Conant.2017p4} established that it is in fact possible to bound the index of $H$ in this way.

Theorems \ref{thm:vc} and \ref{thm:stablegps} can be viewed as algebraic analogues of efficient graph regularity lemmas proved under the assumptions of bounded $\VC$-dimension \cite{Alon.2007, Lovasz.2010} and stability \cite{Malliaris.2014},  respectively.  Theorems \ref{thm:vc} and \ref{thm:stablegps} were subsequently strengthened and generalized by several authors in \cite{Sisask.2018, Conant.2018n5,Conant.2018zd} and \cite{Conant.2017p4,Conant.2021}, respectively.  
 
Taken together, these results provide a satisfactory answer to the question of what combinatorial restrictions ensure especially efficient linear decompositions. In light of the discussion above regarding higher order Fourier analysis, they naturally lead to the following question.
 
 \begin{question}\label{question}
 What combinatorial restrictions ensure especially efficient higher-order (specifically, quadratic) decompositions?
 \end{question}

It turns out that Question \ref{question} is extremely rich. The goal of this paper is to investigate  Question \ref{question}, focusing on higher-order analogues of $\VC$-dimension and Theorem \ref{thm:vc}. Higher-order analogues of stability and Theorem \ref{thm:stablegps} are explored in our companion paper \cite{Terry.2021a}.\footnote{The main results of this  paper and its companion originally appeared together in \cite{Terry.2021av2}.}

The main result of the present paper is a quadratic analogue of Theorem \ref{thm:vc} under the assumption that the subset in question has bounded $\VC_2$-dimension. This ternary analogue of $\VC$-dimension was first defined by Shelah in \cite{Shelah.2014}, and was subsequently shown to exhibit many pleasant properties  in \cite{Chernikov.2019}. It was related to hypergraph regularity and other combinatorial problems in \cite{Terry.2018, Chernikov.2020, Terry.2021b}, as well as bilinear forms over vector spaces in \cite{Hempel.2016, Chernikov.2019b}. However, it has received little attention to date compared to stability and NIP, the two notions of tameness discussed earlier in this introduction. Here we define $\VC_2$-dimension in the context of subsets of groups.

\begin{definition}[$\VC_2$-dimension of a subset of a group]\label{def:vc2dim}
Given $m\geq 1$, an abelian group $G$ and a subset $A\subseteq G$, we say that $A$ \emph{has $\VC_2$-dimension at least $m$} (or \emph{has $m$-$\IP_2$}) if there exist 
$$
\{a_i:i \in[m]\}\cup \{b_{j}:j\in [m]\}\cup \{c_S: i\in [m]\times [m]\}\subseteq G
$$
 such that $a_i+b_j+ c_S \in A$ if and only if $(i,j)\in S$.

The \emph{$\VC_2$-dimension of $A$}, denoted $\VC_2( A)$, is defined to be the largest $m$ such that $A$ has $\VC_2$-dimension at least $m$. When $\VC_2(A)<m$, we also say that \emph{$A$ is $m$-$\NIP_2$}.
\end{definition}

In order to formulate our main result, we will need the notion of a \emph{quadratic factor}, first introduced by Green and Tao \cite{Green.2007} in analogy to related concepts in ergodic theory (see e.g. \cite{Kra.2006}). Quadratic (and analogously defined higher-order) factors serve as descriptions of the structured part of higher-order arithmetic regularity lemmas. For the purpose of this introduction, it suffices to think of a quadratic factor $\mathcal{B}$ as a partition of $\F_p^n$ into the joint level sets of a given collection of quadratic and linear maps.\footnote{Since some additional care needs to be taken when $p=2$, in this paper we shall (sometimes tacitly) assume that the characteristic $p$ is strictly greater than 2.} As is customary, we refer to each of these level sets as an \emph{atom}, and denote the collection of all atoms associated with the quadratic factor $\mathcal{B}$ by $\At(\calB)$. A typical atom is of the form 
$$
\{x \in \F_p^n: x^T M_1 x =a_1, x^T M_2 x =a_2,\dots, x^T M_q x =a_q, x^T r_1=b_1,x^T r_2=b_2,\dots, x^T r_\ell=b_\ell \},
$$
 for some symmetric $n\times n$ matrices $M_1,\ldots, M_q$ with entries in $\F_p$, and some (usually linearly independent) vectors $r_1,\ldots, r_{\ell}$. The atom is said to be \emph{labeled} by $a_1,\dots, a_q$, $b_1,\dots, b_\ell$ in $\F_p$. The \emph{complexity} of the partition $\At(\calB)$ is given by the number of linear and quadratic maps that define it, in this case $(\ell,q)$, and its \emph{rank} is the minimum rank of any non-trivial linear combination of the matrices $M_1,\dots, M_q$. In applications, we require the rank of the factor to be sufficiently high, as this  ensures that the factor has various helpful properties (for example, that its atoms are of roughly equal size).

A partition into cosets of a given subgroup can be viewed as a degenerate case of a quadratic factor which is given by linear maps only. We will therefore henceforth refer to such a partition as arising from a \emph{linear factor}. For a more detailed introduction to linear and quadratic factors, see Section \ref{sec:qfa}.  

It is not difficult to see that a single quadratic atom has $\VC_2$-dimension at most 1 (see Section \ref{subsec:exaQA}). It follows from the fact that the class of sets of bounded $\VC_2$-dimension is closed under Boolean combinations (see \cite{Chernikov.2019}) that a union of finitely many atoms of a quadratic factor has bounded $\VC_2$-dimension (quantitatively, this is shown in \cite[Theorem 3.16]{Terry.2021a}/\cite[Corollary 5.14 ]{Terry.2021av2}).
As a corollary of our main result, Theorem \ref{thm:vc2} below, we will obtain an approximate converse to that statement, namely that sets of bounded $\VC_2$-dimension look approximately like unions of quadratic atoms (see Corollary \ref{cor:vc2}). 

We now state our arithmetic regularity lemma for sets of bounded $\VC_2$-dimension.\footnote{This result originally appeared as \cite[Theorem 1.8]{Terry.2021av2}.}
  
\begin{theorem}[Arithmetic regularity lemma for sets of bounded $\VC_2$-dimension]\label{thm:vc2}
For all primes $p>2$, all integers $m\geq 1$ and all reals $\mu>0$, there exists a polynomial growth function $\sigma_1=\sigma_1(p,m,\mu):\mathbb{R}^+\rightarrow \mathbb{R}^+$ such that for all $\sigma\geq\sigma_1$ (pointwise) there exist $n_1=n_1(p,m,\mu,\sigma)$ and $D=D(p,m,\mu, \sigma)$ such that the following holds.  

Let $n\geq n_1$ be an integer. Suppose that $A\subseteq \F_p^n$ satisfies $\VC_2(A)<m$.  Then there is a quadratic factor $\calB=(\calL,\calQ)$ on $\F_p^n$ of complexity $(\ell,q)$, and a set $\Xi \subseteq \At(\calB)$ such that 
\begin{enumerate}[label=\normalfont(\roman*)]
\item $\ell, q\leq D$; 
\item $\calB$ has rank at least $\sigma(\ell+q)$;
\item $|\Xi|\leq \mu |\At(\calB)|$; 
\item for all $B\in \At(\calB)\setminus \Xi$, $|A\cap B|/|B|\in [0,\mu)\cup (1-\mu,1]$.
\end{enumerate}
\end{theorem}

This result has explicit connections to an analogous result of the authors in the hypergraph setting \cite{Terry.2021b} (see also \cite{Chernikov.2020} for related results on hypergraphs and various higher arity versions of VC-dimension). 

It turns out that even though linear in structure, the example of Green and Sanders \cite{Green.2015} cited in the context of Theorem \ref{thm:vc} also shows that a quadratic decomposition of a set of bounded $\VC_2$-dimension must necessarily feature an error atom. This shows that the set $\Xi$ in Theorem \ref{thm:vc2} cannot be dispensed with. For a proof of this, we refer the reader to our companion paper (see \cite[Proposition 3.26]{Terry.2021a}/\cite[Proposition 5.22]{Terry.2021av2}).

The following corollary of Theorem \ref{thm:vc2} is immediate.\footnote{This result originally appeared as \cite[Corollary 1.9]{Terry.2021av2}.}

\begin{corollary}[Structure theorem for sets of bounded $\VC_2$-dimension]\label{cor:vc2}
For all primes $p>2$, all integers $m\geq 1$ and all reals $\mu>0$, there exists a polynomial growth function $\sigma_2=\sigma_2(p,m,\mu):\mathbb{R}^+\rightarrow \mathbb{R}^+$ such that for all $\sigma\geq\sigma_2$  (pointwise) there exist $n_1=n_1(p,m,\mu,\sigma)$ and $D=D(p,m,\mu,\sigma)$ such that the following holds.  

Let $n\geq n_1$ be an integer. Suppose that $A\subseteq \F_p^n$ satisfies $\VC_2(A)<m$.  Then there is a quadratic factor $\calB=(\calL,\calQ)$ on $\F_p^n$ of complexity $(\ell,q)$, and a union $Y$ of atoms of $\calB$ such that 
\begin{enumerate}[label=\normalfont(\roman*)]
\item $\ell, q\leq D$;
\item $\calB$ has rank at least $\sigma(\ell+q)$;
\item $|A\Delta Y|\leq \mu|G|$. 
\end{enumerate}
\end{corollary}

While one might expect the bounds on the complexity of the quadratic factor obtained under the additional assumption of bounded $\VC_2$-dimension to be considerably stronger than in the general case, the initial approach taken in this paper does not provide such an improvement. Indeed, our proof of Theorem \ref{thm:vc2} uses a general quadratic regularity lemma (see Section \ref{subsec:quadarl}) and thus produces Tower-type bounds (in $\mu^{-1}$) for both $\ell$ and $q$. In a forthcoming paper \cite{Terry.2024f}, we improve the bound on $q$ to $\log_p(\mu^{-O_k(1)})$ (resulting in a partition size that is polynomial in $\mu^{-1}$).

\textbf{Perspectives.} In the context of groups definable in theories of finite $\VC_k$-dimension, Shelah \cite{Shelah.2007} ($k=2$) and  Chernikov and Hempel \cite{Chernikov.2019b} ($k>2$) showed ``relative absoluteness" for the subgroup $G^{00}$.  Thus, these attempts at extending stable group theory into the realm of bounded $\VC_k$-dimension hint at the existence of bilinear structure but fall short of giving an explicit algebraic description along the lines of Theorem \ref{thm:vc2}. Moreover, the problem of producing a truly higher order analogue of stability, along with the corresponding tools of stable and NIP group theory, has remained open.  We show, via our answers to Question \ref{question} here and in \cite{Terry.2021a}/\cite{Terry.2021av2}, that a missing ingredient in this pursuit has been the machinery of higher Fourier analysis, and the correspondence it lays bare between higher \emph{arity} relations and higher \emph{degree} structure in groups.  

While the approach in this paper is combinatorial and finitary in nature, we expect a rich model-theoretic framework to emerge over the coming years that encompasses higher arity notions of $\VC$-dimension and stability, and higher degree polynomial structure. As was the case for the aforementioned work on linear decompositions, we expect the model theoretic tools, once developed, to be vastly more general than those available to arithmetic combinatorialists at this stage.

\textbf{Overview of the argument.} As indicated above, the proof of Theorem \ref{thm:vc2} begins with an application of the standard quadratic regularity lemma (see Theorem \ref{thm:gtarl}), which asserts, roughly speaking, that given any subset $A\subseteq \F_p^n$, we may decompose $\F_p^n$ into simultaneous level sets of a bounded number of linear and quadratic phases (i.e. into the atoms of a bounded-complexity quadratic factor) in such a way that $A$ is quadratically uniform (as measured by the Gowers $U^3$ norm) ``with respect to almost all atoms". Whilst the functional version of this statement is well known, to the best of our knowledge the version we give as Proposition \ref{prop:unifatom2} has no precedent in the literature. Indeed, one key innovation of this paper is the introduction of a suitably localised $U^3$ semi-norm (Definition \ref{def:localu3}) which allows us to make the aforementioned statement precise. We then prove, in Proposition \ref{prop:trivdense}, that if a set $A$ of bounded $\VC_2$-dimension is locally quadratically uniform with respect to a quadratic atom, then it must have density close to $0$ or close to $1$.  This is analogous to \cite[Proposition 1]{Terry.2019} (with the local $U^3$ semi-norm in place of the local $U^2$ semi-norm, and $\VC_2$-dimension in place of $\VC$-dimension), and implies Theorem \ref{thm:vc2}.

\textbf{Note on earlier versions of this work.} Theorem \ref{thm:vc2} and its Corollary \ref{cor:vc2} first appeared in the arXiv preprint \cite{Terry.2021av2} of the authors. Due to the length and complexity of that paper, it seemed profitable to separate the result for bounded $\VC_2$-dimension ($\NIP_2$) from the development of the generalised notion of stability ($\NFOP_2$), which constituted the main focus of \cite{Terry.2021a}. The proof of Theorem \ref{thm:vc2} has been completely reworked for the present paper, avoiding the language of hypergraphs employed in \cite{Terry.2021av2} and instead relying on local uniformity norms to control the relevant averages throughout Sections \ref{sec:localnorms}-\ref{sec:mainproof}.

\textbf{Outline of the paper.} We begin, in Section \ref{sec:qfa}, by giving the technical background on quadratic Fourier analysis over $\F_p^n$ that will be necessary for the remainder of the paper. In fact, we take some care to motivate the development of quadratic Fourier analysis with the model theoretic reader in mind. For the expert in additive combinatorics, there is little of novelty in this section. In Section \ref{sec:localnorms}, we introduce the local $U^2$ and $U^3$ semi-norms, and establish some of their basic properties. In particular, we demonstrate (in conjunction with the content of the two appendices) that the local $U^3$ semi-norm controls (at least in an approximate sense) the local $\IP_2$-operator. Using the local $U^2$ semi-norms, we also prove some basic properties of quadratic atoms. In particular, we show that they have unbounded $\VC$- but bounded $\VC_2$-dimension. Much more is true, but we leave this to our companion paper \cite{Terry.2021a}. In Section \ref{sec:mainproof}, after some preparatory work, we give the proof of Theorem \ref{thm:vc2}. Appendix \ref{app:expsums} is dedicated to proving the basic exponential sum estimates over bilinear forms that underlie much of the technical part of the paper, whilst Appendix \ref{app:counting} contains the proof of the approximate control of the local $\IP_2$-operator by the corresponding local $U^3$ semi-norm. In Appendix \ref{app:locsparseuni} we show that locally sparse sets are locally uniform, and in Appendix \ref{app:countinggen} we generalise our counting results to operators defined over possibly distinct tuples of atoms, for use in \cite{Terry.2021b} and the forthcoming work \cite{Terry.2024f}.

\textbf{Acknowledgements.}
The first author would like to thank Maryanthe Malliaris for many useful mathematical conversations. The second author is grateful to Tim Gowers for several helpful conversations during the preparation of the earlier version of this manuscript. During the development of this work, the first author was partially supported by NSF Grant DMS-1855711, NSF CAREER Award DMS-2115518, and a Sloan Research Fellowship. During the final stages of preparing this manuscript, the work of the second author was supported by an Open Fellowship from the UK Engineering and Physical Sciences Research Council (EP/Z53352X/1). The authors' collaboration has been supported by several travel grants over the years, including by the London Mathematical Society, the Simons Foundation, and the Association for Women in Mathematics.

\textbf{Notational conventions}.
For a function $f$ and positive-valued function $g$ on a domain $D$, we write $f = O(g)$ if there exists a constant $C$ such that $|f(x)| \leq Cg(x)$ for all $x\in D$. One crucial instance of this notation will be when we write $f=(1+O(g))h$ for some positive-valued function $h$ on $D$, which means that for some constant $C$, $|f(x)-h(x)|\leq Cg(x)h(x)$ holds for all $x\in D$.   When $C$ depends on another parameter, this is (usually) indicated by a subscript on $O$.  

We denote by $\mathcal{C}$ the complex conjugate operator, taking a function $f:D\rightarrow \C$ to a function $\mathcal{C}f:D\rightarrow \C$ defined by setting $\calC{f}(x)=\overline{f(x)}$.  Given an integer $k\geq 1$, $\mathcal{C}^k$ is the operator obtained by applying this operator $k$ times (so $\mathcal{C}^kf=f$ for $k$ even and $\mathcal{C}^kf=\mathcal{C}f$ for $k$ odd). Given $\e\in \{0,1\}^k$ and $1\leq i\leq k$, $\e(i)$ denotes the $i$th coordinate of $\e$, and $|\e|$  denotes the number of 1s in $\e$.

Given a set $X$, we write $\E_{x\in X}$ to mean $\frac{1}{|X|}\sum_{x\in X}$.  When the set $X=\F_p^n$ we will simply write $\E_x$ for $\E_{x\in \F_p^n}$. Given a finite index set $I$ and a set $X$, we write $\E_{x_i\in X:\; i\in I}$ to mean $\E_{x_{i_1}\in X}\ldots \E_{x_{i_k}\in X}$, where $I=\{i_1,\ldots, i_k\}$ is some enumeration.  There is some ambiguity in this notation as to which enumeration is being used, but we assure the reader that any time this notation appears, the value of the equation will be the same regardless of which enumeration is chosen.  When $X=\F_p^n$ we will abbreviate $\E_{x_i\in \F_p^n: \;i\in I}$  by writing simply $\E_{x_i:\; i\in I}$.  To save space, in some cases we will also write $\E_{\begin{subarray}{l} x_i\in X\\ i\in I\end{subarray}}$ to mean $\E_{x_i\in X:\; i\in I}$.

Given a set  indexed by $I$, say $\{y_i: i\in I\}$, we write $(y_i)_{i\in I}$ to mean  $(y_{i_1},\ldots, y_{i_k})$ where $I=\{i_1,\ldots, i_k\}$ is some enumeration. Again there is ambiguity in this notation, but this will not matter as long as we assume the enumeration chosen is consistent throughout any given argument.  Given $i^*\in I$, we will write $\E_{x_i: i\in I, i\neq i^*}$ to mean $\E_{x_i: i\in I\setminus \{i^*\}}$ and $(y_i)_{i\in I, i\neq i^*}$ to denote $(y_i)_{i\in I\setminus \{i^*\}}$.  When the ambient set $I$ is clear from context, we will sometimes write $(y_i)_{i\neq i^*}$ to mean $(y_i)_{i\in I, i\neq i^*}$. Finally, when the index set $I$ is a power set, say $I=\calP(J)$, we will often use substitute $S\subseteq J$ instead of $S\in \calP(J)$ in all of the above notation. 

Throughout the paper, $p>2$ will denote a fixed prime, and $\omega$ will denote a $p$th root of unity. We will frequently use the exponential sum identity 
\begin{align}\label{star}
\E_{y\in \F_p^n}\omega^{r^T y}=\begin{cases} 1 \text{ if }&r=0\\
0 \text{ if }& r\in \F_p^n\setminus \{0\} .\end{cases}
\end{align}

In several places, we will  use an immediate corollary of the above identity, obtained by combining several instances of (\ref{star}) in the case where $n=1$. Namely, given $r_1,\ldots, r_t\in \F_p$, we have
\begin{align}\label{starstar}
\E_{y_i\in \F_p: i\in [t]}\omega^{\sum_{i=1}^t y_i r_i}=\begin{cases} 1 &\text{ if }r_1=\ldots=r_t=0\\
0 &   \text{ otherwise. }\end{cases}
\end{align}

To aid clarity, we shall occasionally denote the group of characters of $\F_p^n$ by $\widehat{\F_p^n}$, even though the latter will always be identified with $\F_p^n$. We shall follow the convention that for $q\in [1,\infty)$, the $L^q$ norms on physical space are normalised, i.e. $\|f\|_{L^q}^q=\E_{x\in \F_p^n} |f(x)|^q$, while the $\ell^q$ norms on frequency space are not. That is, if $\widehat{f}:\F_p^n\ra \C$ is the Fourier transform of $f:\F_p^n \ra\C$ (see Definition \ref{def:ft}), then $\|\widehat{f}\|_{\ell^q}^q=\sum_{t\in \widehat{\F_p^n}}|\widehat{f}(t)|^q$. A similar convention is in force for inner products. Finally, we define $\|f\|_{\infty}=\sup_{x\in \F_p^n}|f(x)|$ and $\|\widehat{f}\|_{\infty}=\sup_{t\in \widehat{\F_p^n}}|\widehat{f}(t)|$.

\section{Background on quadratic Fourier Analysis}\label{sec:qfa}

Higher-order Fourier analysis has its roots in Gowers's work on Szemer\'edi's theorem \cite{Gowers.1998, Gowers.2001}. It has since been developed in multiple important directions and found analogues and applications in adjacent fields. For the purpose of this brief overview, we shall only mention the groundbreaking work of Host and Kra in ergodic theory \cite{Host.2005} and efforts of Szegedy and co-authors (culminating in \cite{Candela.2019}, see also references therein) to generalize the central results beyond a small class of finite abelian groups.

The recent survey \cite{Gowers.2016s0e} by Gowers gives an excellent overview of these various developments, with a strong focus on motivating problems in arithmetic combinatorics and number theory. A more technical but nevertheless outstanding introduction to the setting of vector spaces over a fixed finite field in particular was given by Green \cite{Green.2007}. For a more infinitary viewpoint, we recommend the book by Tao \cite{Tao.2012cv}.

We return to the underlying theme set out in the introduction, namely our goal to decompose a set (or a function) into a structured and a pseudorandom part. For certain applications, such as counting 3-term arithmetic progressions (3-APs), measuring pseudorandomness by the maximum size of the Fourier cofficients of the function is perfectly adequate, and the associated structured part will be pleasantly linear in nature (in a sense to be made precise momentarily). For other applications, such as counting 4-term and longer progressions, a more sophisticated measure of pseudorandomness is required, and the object of primary interest in this paper is the resulting higher-order polynomial (specifically, quadratic) structure that forms its counterpart.

In the moderately more detailed discussion that follows, which was written with the model-theoretic reader in mind, we shall restrict our attention to the setting of bounded functions defined on $\F_p^n$ with $p>2$. Almost all of the concepts introduced can be generalized to other abelian groups, albeit often only with significant effort. The reader familiar with the literature may skip Sections \ref{subsec:uknorms}-\ref{subsec:quadarl} but should be aware that we will be using some results in the form in which they are stated here (which may not always be entirely standard).

\subsection{The Gowers $U^2$ norm and the IP operator}\label{subsec:uknorms}

As far back as the 1950s, Roth \cite{Roth.1953} observed that the $\ell^\infty$ norm of the Fourier transform of the characteristic function of a set $A$ can be used to control the count of 3-APs contained in $A$. While we will make only minimal use of the Fourier transform in this paper, we define it for completeness and to facilitate comparison with \cite{Terry.2019}. Recall that we write $\omega$ for a $p$th root of unity.

\begin{definition}[Fourier transform]\label{def:ft}
Given a function $f:\F_p^n\ra \C$ and $t\in \F_p^n$, define the \emph{Fourier coefficient of $f$ at $t$} to be
\[\widehat{f}(t)=\E_{x\in \F_p^n}f(x)\omega^{-x^T t}.\]
\end{definition}

Here and elsewhere, the notation $\E_{x\in B}$ will be used to denote the normalise average $\frac{1}{|B|}\sum_{x\in B}$. Where the set $B$ is omitted, the average is understood to be taken over the entire group $\F_p^n$. It is easy to prove the \emph{inversion formula}, i.e. for all $x\in \F_p^n$
\[f(x)=\sum_{t\in \widehat{\F_p^n}} \widehat{f}(t)\omega^{x^T t},\]
as well as \emph{Parseval's identity} $\|f\|_{L^2}=\|\widehat{f}\|_{\ell^2}$. (We refer the reader to the normalisation conventions at the end of the introduction.)

In today's language, and adapted to our setting, Roth's starting point can then be stated as follows.

\begin{fact}[Fourier transform controls 3-APs]\label{fact:3apcontrol}
Let $f:\F_p^n\ra \C$ be such that $\|f\|_\infty \leq 1$. Then 
\[|\E_{x,d}f(x)f(x+d)f(x+2d)|\leq \|\widehat{f}\|_\infty.\]
\end{fact}

Since a function taking values in $\{-1,1\}$ with equal probability has vanishingly small Fourier coefficients at all frequencies with high probability, it is reasonable to regard functions $f$ for which $\|\widehat{f}\|_\infty$ is small as \emph{pseudo}- or \emph{quasi}-random (to add further terminology, they are often referred to as \emph{Fourier uniform}).

It turns out to be more useful for the study of general linear patterns to replace the $\ell^\infty$ norm of the Fourier transform by a more combinatorial measure of uniformity\footnote{The definition of the $U^2$ norm is in close correspondence with a well known measure of quasirandomness in graphs, namely the count of 4-cycles. The insight that it makes a fruitful measure of quasirandomness in the setting of subsets of groups is due to Gowers.}. We give the definition in two parts, as this will be useful to us later. (Again, we refer the reader to the end of the introduction for the notation used.)

\begin{definition}[$U^2$ inner product]\label{def:ipu2}
Given functions $(f_{\e})_{\e\in \{0,1\}^2}:\F_p^n\ra \C$, we define their \emph{$U^2$ inner product} by
\[\langle (f_{\e})_{\e\in \{0,1\}^2}\rangle_{U^2}=\E_{\begin{subarray}{l} x_0,x_1\end{subarray}}\E_{\begin{subarray}{l}y_0,y_1\end{subarray}}\prod_{\e\in \{0,1\}^2}\calC^{|\e|}f_{\e}(x_{\e(1)}+y_{\e(2)}).\] 
\end{definition}

Note that by definition, $\langle (f_{\e})_{\e\in \{0,1\}^2}\rangle_{U^2}$ is linear in $f_{00}$ and $f_{11}$ and conjugate linear in $f_{01}$ and $f_{11}$. In the same way as the standard inner product allows us to define the $L^2$ norm, the $U^2$ inner product allows us to define an associated norm.

\begin{definition}[$U^2$ norm]\label{def:u2}
Given $f: \F_p^n\ra \C$, define its \emph{$U^2$ norm} by
\[\|f\|_{U^2}= \langle (f)_{\e\in \{0,1\}^2}\rangle_{U^2}^{1/4}.\]
\end{definition}

It is not a priori obvious that this is well-defined, or that it defines a norm. This follows immediately from the following fact, which is an easy exercise (see e.g. \cite[Proposition 1.9]{Green.2007}).

\begin{fact}\label{fact:fourierequiv}
For $f:\F_p^n\ra\C$, we have $\|f\|_{U^2}^4=\|\widehat{f}\|_4^4$. Furthermore, if $\|f\|_\infty\leq 1$, we have $\|\widehat{f}\|_\infty^4 \leq \|\widehat{f}\|_4^4=\|f\|_{U^2}^4\leq \|\widehat{f}\|_\infty^2$.
\end{fact}

This means that, at least loosely speaking, a function $f$ has small $U^2$ norm if and only if its Fourier transform has small $\ell^\infty$ norm\footnote{Indeed, the letter ``U" in $U^2$ stands for ``uniformity".}. As a corollary of Facts \ref{fact:3apcontrol} and  \ref{fact:fourierequiv}, we see that the 3-AP count, weighted by a function $f$, is controlled by the $U^2$ norm of an associated auxiliary function. 

In order to interpret this statement in the context of a subset $A\subseteq \F_p^n$ of density $\alpha=|A|/|\F_p^n|$, write $1_A$ for its characteristic function and $f_A=1_A-\alpha$ for its so-called \emph{balanced function}. It is not hard to infer from the above discussion that if $\|f_A\|_{U^2}$ is small, then $A$ contains roughly the number of 3-APs expected in a random\footnote{For the model theorist, by a \emph{random set of the same density} we mean a subset $A$ of $\F_p^n$ such that each element of $\F_p^n$ has been chosen to lie in $A$ independently with probability $\alpha$.} subset of $\F_p^n$ of the same density (see e.g. \cite[Proposition 1.8]{Green.2007}).

Something similar is true of a configuration that is very much of interest to us in this paper, namely copies of the independence property (see Definition \ref{def:vcdim}), which we will count with the help of the following operator.

\begin{definition}[$\IP$-operator]\label{def:ipop}
Let $m\geq 2$. Given functions $(f_{i,S})_{i\in [m], S\subseteq [m]}:\F_p^n\ra \C$, we define 
\[T_{m-\IP}((f_{i,S})_{i\in [m],S\subseteq [m]})=\E_{x_i:i\in[m]}\E_{y_S:S\subseteq [m]}\prod_{i\in [m],S\subseteq [m]}f_{i,S}(x_i+y_S).\]
\end{definition}

\begin{notation}\label{not:ipop}
For all  operators related to the independence property and its generalisations, we adopt the following convention. When all the arguments $(f_{i,S})_{i,j\in [m],S\subseteq [m]}$ are identical, we only write one function as input for the operator, that is,
\[T_{m-\IP}(f)=T_{m-\IP}((f)_{i\in [m],S\subseteq [m]}).\]
We also write $T_{m-\IP}((f_{i,S})_{i\in [m],S\subseteq [m]})$ as
\[T_{m-\IP}((f_{i,S})_{i\in [m],i\in S\subseteq [m]}|(f_{i,S})_{i\in [m],i\notin S\subseteq [m]}).\]
When all $(f_{i,S})_{i\in [m],i\in S\subseteq [m]}$ are identical and equal to $f$ say, and all $(f_{i,S})_{i \in [m],i \notin S\subseteq [m]}$ are identical and equal to $g$ say, we will write 
\[T_{m-\IP}(f|g)=T_{m-\IP}((f)_{i\in [m],i\in S\subseteq [m]}|(g)_{i\in [m],i\notin S\subseteq [m]}).\]
\end{notation}

In particular, we will later on apply this operator with $f=1_A$ and $g=1_{A^C}$, where $A$ is a subset of $\F_p^n$ and $A^C=\F_p^n\setminus A$. In this case, $T_{m-\IP}(1_A|1_{A^C})$ counts instances of induced copies of the $\IP$-graph in the sum-graph associated with $A$.

For the expert, the following result follows straight from the (now standard) theory of Cauchy-Schwarz complexity developed by Green and Tao in \cite{Green.2010ta4}. Because it is instructive, and because we shall need it as a blueprint for a somewhat less standard argument later, we include a direct proof here. 

\begin{lemma}[$\IP$ is controlled by $U^2$]\label{lem:ipcontrol}
Let $m\geq 2$. Suppose that for each $i\in [m]$ and $S\subseteq [m]$,  $f_{i,S}:\F_p^n\ra \C$ satisfies $\|f_{i,S}\|_{\infty}\leq 1$. Then 
\[|T_{m-\IP}((f_{i,S})_{i\in [m],S\subseteq [m]})|\leq \min_{i\in [m],S\subseteq [m]}\|f_{i,S}\|_{U^2}.\]
\end{lemma}

\begin{proof} Fix $I\in [m]$ and $\calS\subseteq [m]$ achieving $\min_{i\in [m],S\subseteq [m]}\|f_{i,S}\|_{U^2}$. Then 
\begin{align}\label{al:firstcs}
|T_{m-\IP}&((f_{i,S})_{i\in [m],S\subseteq [m]})|^4\\\nonumber
=&\big|\E_{x_i:i\in[m], i\neq I}\E_{y_S:S\subseteq [m]}g((x_i)_{i\in[m],i\neq I},(y_S)_{S\subseteq [m]})\E_{x_{I}}f_{I,\calS}(x_{I}+y_{\calS})\prod_{S\subseteq [m],S\neq \calS} f_{I,S}(x_{I}+y_{S})\big|^4,
\end{align}
where $g((x_i)_{i\in[m],i\neq I},(y_S)_{S\subseteq [m]})=\prod_{i\in [m],i\neq I}\prod_{S\subseteq [m]}f_{i,S}(x_i+y_S)$ is a 1-bounded function that does not depend on $x_{I}$.
To (\ref{al:firstcs}), we may apply Cauchy-Schwarz to obtain an upper bound of 
\[\big(\E_{x_i:i\in[m],i\neq I}\E_{y_S:S\subseteq [m]}\big|\E_{x_{I}}f_{I,\calS}(x_{I}+y_{\calS})\prod_{S\subseteq [m],S\neq \calS} f_{I,S}(x_{I}+y_{S})\big|^2\big)^2.\]
Observing that we may drop the first expectation, upon expanding the inner square and rearranging we obtain
\[\big(\E_{y_S:S\subseteq [m],S\neq \calS}\E_{x_{I},x_{I}'}h(x_{I},x_{I}',(y_S)_{S\neq \calS})\E_{y_{\calS}}f_{I,\calS}(x_{I}+y_{\calS})\overline{f_{I,\calS}(x_{I}'+y_{\calS})}\big)^2,\]
where $h(x_{I},x_{I}',(y_S)_{S\neq \calS})=\prod_{S\subseteq [m],S\neq \calS} f_{I,S}(x_{I}+y_{S})\prod_{S\subseteq [m],S\neq \calS} \overline{f_{I,S}(x_{I}'+y_{S})}$ is a 1-bounded function that does not depend on $y_{\calS}$. Therefore, applying Cauchy-Schwarz a second time, we obtain an upper bound of
\[\E_{x_{I},x_{I}'}\big|\E_{y_{\calS}}f_{I,\calS}(x_{I}+y_{\calS})\overline{f_{I,\calS}(x_{I}'+y_{\calS})}\big|^2=\|f_{I,\calS}\|_{U^2}^4,\] 
as claimed.
\end{proof}

As was the case for 3-APs, the above lemma allows one to deduce that if the balanced function of a subset $A$ of $\F_p^n$ has small $U^2$ norm, then $A$ contains approximately the number of induced $\IP$-configurations one would expect in a random set.

Having shown that functions with small $U^2$ norm have well controlled $\IP$-count, it remains to extract a useful conclusion from the complementary event, namely that $\|f\|_{U^2}$ is non-negligible. This follows straight from Fact \ref{fact:fourierequiv} and the definition of the Fourier transform: if $\|f\|_{U^2}\geq \eta$, then there exists a linear complex exponential phase $\gamma(x)=\omega^{-x^T t}$ for some $t\in \F_p^n$ such that $\E_x f(x)\gamma(x)  \geq \eta^2$.

\subsection{An arithmetic regularity lemma for the $U^2$ norm}\label{subsec:arl}

The conclusions arrived at in the preceding section represent the two sides of the structure-randomness coin: either the set $A$ in question is uniform (in the sense that its balanced function has small $U^2$ norm), or we are able to identify some structure (in the form of a non-negligible correlation with a linear phase). These two statements can be refined and combined, in more or less sophisticated ways, to give a proof of Roth's theorem in $\F_p^n$ (see e.g. \cite{Green.2007, Wolf.2015} and references therein). 

A particularly relevant way of formulating this dichotomy, reminiscent of Szemer\'edi's regularity lemma \cite{Szemeredi.1975} in graph theory, is given by Theorem \ref{thm:linearreg} below.  To bring it in line with later sections, we will state Theorem \ref{thm:linearreg} in terms of the following notion.

\begin{definition}[Linear factor]\label{def:linfac}
Let $\ell\geq 1$ be an integer. A \emph{linear factor $\calL$ in $\F_p^n$ of complexity $\ell$} is a set of linearly independent\footnote{It is more standard to not insist on linear independence, and to say instead that the factor has complexity at most $\ell$. We find the present approach more convenient for our purposes.} vectors $r_1,\ldots, r_{\ell}\in \F_p^n$. 
 The \emph{atoms of $\calL$} are sets of the form
$$
\{x\in \F_p^n: x^T r_1=a_1,\ldots, x^T r_{\ell}=a_{\ell}\},
$$
for $a_1,\ldots, a_{\ell}\in \F_p$.  We let $\At(\calL)$ denote the set of atoms of $\calL$.  
\end{definition}

Note that given a linear factor $\calL$ in $\F_p^n$ of complexity $\ell$, $\At(\calL)$ is a partition\footnote{The term ``factor" usually refers to the partition generated by the atoms of $\calL$ (see e.g. \cite{Green.2007}).  For convenience, we identify the factor with its set of constraints.} of $\F_p^n$ with $p^{\ell}$ parts, each of which has size $p^{n-\ell}$.   It is immediate that $\At(\calL)$ is in fact the set of cosets of a subgroup of codimension $\ell$.

As is standard in the field (see e.g. \cite[Definition 3.3]{Green.2007}), we write $\E(f|\mathcal{\calL})$ for the projection of $f$ onto $\calL$, that is, $\E(f|\mathcal{\calL})$ is the function whose value at $x$ is the average of $f$ over the atom of $\calL$ that contains $x$.  More specifically, given $x\in \F_p^n$, $\E(f|\mathcal{\calL})(x)=\E_{y\in L}f(y)$, where $L$ is the element of $\At(\calL)$ containing $x$. 

We also define a \emph{growth function} to be any increasing function $\sigma:\mathbb{R}^+\rightarrow \mathbb{R}^+$.  Given two growth functions $\sigma, \sigma'$, we write $\sigma \geq \sigma'$ to denote that $\sigma$ is bounded below pointwise by $\sigma'$, i.e. $\sigma(x)\geq \sigma'(x)$ for all $x\in \mathbb{R}^+$. 

\begin{theorem}[$U^2$ arithmetic regularity]\label{thm:linearreg}
For every $\e>0$ and for any growth function $\sigma:\R^+\rightarrow\R^+$ (which may depend on $\e$), there exist $n_0=n_0(\e,\sigma)$ and $C=C(\e,\sigma)$ such that the following holds.
Let $n>n_0$ and let $f:\F_p^n \ra [-1,1]$ be a function. 
Then there is a linear factor $\calL$ of complexity $\ell$ such that
\begin{enumerate}[label=\normalfont(\roman*)]
\item the complexity of $\calL$ satisfies $\ell \leq C$;
\item we can write 
\[f=f_1+f_2+f_3,\]
with $f_1=\E(f|\mathcal{\calL})$, $\|f_2\|_{U^2} \leq 1/\sigma(\ell) $, and $\|f_3\|_2\leq \e$.
\end{enumerate}
\end{theorem}

The first such regularity lemma was proved by Green \cite{Green.2005}, and provided the starting point for the investigations in our first joint paper \cite{Terry.2019}. At first glance, the main result of \cite{Green.2005} looks rather different, being stated in term of a \emph{local} notion of (Fourier) uniformity: given $A\subseteq \F_p^n$ and $H\leqslant \F_p^n$, a coset $g+H$ was defined to be \emph{$\e$-uniform with respect to $A$} if $\sup_{t\neq 0}|\widehat{1_{A-g}|_H}(t)|\leq \e$, where we have written $g|_X$ for the restriction of a function $g:\F_p^n\ra \C$ to a set $X\subseteq \F_p^n$, and the Fourier transform is taken relative to the group $H$.\footnote{In \cite[Definition 4]{Terry.2019} we observed that, equivalently, a coset $g+H$ is $\e$-uniform with respect to $A$ if $\|\widehat{f_A^g}\|_\infty\leq \e$, where $f_A^g(x)=(1_{(A-g)\cap H}(x)-\alpha_{g+H})\mu_H(x)$ is the balanced function of $A$ localised to $g+H$, $\alpha_{g+H}=|(g+H)\cap A|/|H|$ is the density of $A$ on $g+H$, $\mu_H(x)=(|\F_p^n|/|H|)1_H(x)$ is the so-called \emph{characteristic measure} of $H$, and the Fourier transform is taken relative to $\F_p^n$.}

In this terminology, the arithmetic regularity lemma of Green \cite{Green.2005} asserts that given  $A\subseteq \F_p^n$ and $\e>0$, there is a subspace $H\leqslant \F_p^n$ of codimension bounded only in terms of $\e$, such that all but an $\e$-proportion of the cosets of $H$ are $\e$-uniform with respect to $A$. This follows from Theorem \ref{thm:linearreg} above with only a little work.

\subsection{The Gowers $U^3$ norm and the $\IP_2$-operator}\label{subsec:uknorms2}

It is a well known fact that the Fourier transform (and, by Fact \ref{fact:fourierequiv}, the $U^2$ norm) is not sensitive enough to control longer arithmetic progressions or other configurations of similar \emph{complexity}; see, for instance, the introduction to \cite{Gowers.2010tc}. There, the authors showed that the Fourier transform does not control patterns for which the squares of the linear forms defining the pattern are linearly dependent\footnote{Such patterns are said to be of true complexity (at least) 2.}, and went on to show that linear \emph{in}dependence of the squares is not only a necessary but a sufficient condition for a linear pattern to be controlled by the Fourier transform (or, equivalently, the $U^2$ norm). This deep phenomenon also lies (indirectly) at the heart of the present paper.

However, unlike the Fourier transform it is possible to generalize the $U^2$ norm in such a way that its generalization does control 4-APs and other patterns of higher complexity. Indeed, our presentation of the definition of the $U^2$ norm in Section \ref{subsec:uknorms} was chosen to suggest a fruitful generalisation, which is again due to Gowers.

\begin{definition}[$U^3$ inner product]\label{def:ipu3}
Given functions $(f_{\e})_{\e\in \{0,1\}^3}:\F_p^n\ra \C$, we define their \emph{$U^3$ inner product} by
\[\langle (f_{\e})_{\e\in \{0,1\}^3}\rangle_{U^3}=\E_{\begin{subarray}{l} x_0,x_1\end{subarray}}\E_{\begin{subarray}{l}y_0,y_1\end{subarray}}\E_{\begin{subarray}{l}z_0,z_1\end{subarray}}\prod_{\e\in \{0,1\}^3}\calC^{|\e|}f_{\e}(x_{\e(1)}+y_{\e(2)}+z_{\e(3)}).\]
\end{definition}

Observe that by definition, $\langle (f_{\e})_{\e\in \{0,1\}^3}\rangle_{U^3}$ is linear in $f_\e$ with $|\e|$ even and conjugate linear in $f_\e$ with $|\e|$ odd.

\begin{definition}[$U^3$ norm]\label{def:u3}
Given $f: \F_p^n\ra \C$, define its \emph{$U^3$ norm} by
\[\|f\|_{U^3}= \langle (f)_{\e\in \{0,1\}^2}\rangle_{U^3}^{1/8}.\]
\end{definition}

More generally, we may define the $U^d$ norm for $d>2$ in a similar fashion (or, alternatively, by induction via $\|f\|_{U^d}^{2^d}=\E_h\|\Delta_h f\|_{U^{d-1}}^{2^{d-1}}$, where $\Delta_h f(x)=f(x)\overline{f(x+h)}$ \footnote{By defining the semi-norm $\|f\|_{U^1}=|\E_x f(x)|$, Definition \ref{def:u2} becomes consistent with the inductive definition.}). From the inductive definition, it is not difficult to see that the $U^d$ norms are nested, in the sense that for every $d\geq 1$, $\|f\|_{U^{d}}\leq \|f\|_{U^{d+1}}$. 

On the other hand, it is not a priori obvious that the formula in Definition \ref{def:u3} defines a norm. However, equipped with the Gowers-Cauchy-Schwarz inequality below, this is not too hard to prove (see, for example, \cite[Section 11]{Tao.2010}).

\begin{lemma}[Gowers-Cauchy-Schwarz inequality]\label{lem:gcs}
Let $(f_{\e})_{\e\in \{0,1\}^3}:\F_p^n\ra \C$ be functions. Then
\[\left|\langle (f_{\e})_{\e\in \{0,1\}^3}\rangle_{U^3}\right|\leq \prod_{\e \in \{0,1\}^3} \|f_\e\|_{U^3}.\]
\end{lemma}

The $U^3$ norm allows us to generalise the statement in Fact \ref{fact:3apcontrol} to 4-APs. Alongside the definition of the $U^3$ norm and the Gowers-Cauchy-Schwarz inequality above, the following statement first appeared in \cite{Gowers.1998}. Its proof consists of repeated applications of the Cauchy-Schwarz inequality.

\begin{proposition}[$U^3$ controls 4-APs]\label{prop:4apcontrol}
For $f:\F_p^n\ra \C$ with $\|f\|_2\leq 1$, we have
\[|\E_{x,d} f(x)f(x+d)f(x+2d)f(x+3d)|\leq \|f\|_{U^3}.\]
\end{proposition}

We say informally that a set is \emph{quadratically uniform} if its balanced function has small $U^3$ norm. It is not difficult to see (but not entirely obvious either) that Proposition \ref{prop:4apcontrol} implies that sets whose balanced functions are quadratically uniform contain roughly the number of 4-APs expected in a random subset of $\F_p^n$ of the same density.

Again, something similar is true of the following operator.

\begin{definition}[$\IP_2$-operator]\label{def:ip2op}
Let $m\geq 2$. Given functions $(f_{i,j,S})_{i,j\in [m], S\subseteq [m]^2}:\F_p^n\ra \C$, we define
\[T_{m-\IP_2}((f_{i,j,S})_{i,j\in [m],S\subseteq [m]^2})=\E_{x_i:i\in[m]}\E_{y_j:j\in[m]}\E_{z_S:S\subseteq [m]^2}\prod_{i,j\in [m],S\subseteq [m]^2}f_{i,j,S}(x_i+y_j+z_S).\]
\end{definition}

\begin{notation}\label{not:ip2op}
In analogy with Notation \ref{not:ipop}, we write $T_{m-\IP_2}((f_{i,j,S})_{i,j\in [m],S\subseteq [m]^2})$ as 
\[T_{m-\IP_2}((f_{i,j,S})_{i,j\in [m],(i,j)\in S\subseteq [m]^2}|(f_{i,j,S})_{i\in [m],(i,j)\notin S\subseteq [m]^2})\]
to emphasize the product over $(i,j)\in [m]^2, S\subseteq [m]^2$ is taken over functions from the first tuple when $(i,j)\in S$, and functions from the second tuple when $(i,j)\notin S$. 

In particular, if all $(f_{i,j,S})_{i,j\in [m],(i,j)\in S\subseteq [m]^2}$ are identical and equal to $f$ say, and all $(f_{i,j,S})_{i,j \in [m],(i,j) \notin S\subseteq [m]^2}$ are identical and equal to $g$ say, then we will write 
\[T_{m-\IP}(f|g)=T_{m-\IP}((f)_{i,j\in [m],(i,j)\in S\subseteq [m]^2}|(g)_{i,j\in [m],(i,j)i\notin S\subseteq [m]^2}).\]
\end{notation}

Just like Lemma \ref{lem:ipcontrol}, Lemma \ref{lem:ip2control} below follows straight from the work of Green and Tao, by observing that the linear system defining $\IP_2$ has Cauchy-Schwarz complexity 2. This really amounts to saying that the direct argument given for Lemma \ref{lem:ipcontrol} can be extended to $\IP_2$ by carrying out one additional Cauchy-Schwarz step. We leave this as an easy exercise to the reader (but point out that it also follows as a special case from the more sophisticated Proposition \ref{prop:ip2loccontrol}, a proof of which is given in Appendix \ref{app:counting}).

\begin{lemma}[$\IP_2$ is controlled by $U^3$]\label{lem:ip2control}
Let $m\geq 2$. Suppose that for each $i,j\in [m]$, $S\subseteq [m]^2$, $f_{i,j,S}:\F_p^n\ra \C$ is such that $\|f_{i,j,S}\|_{\infty}\leq 1$. Then 
\[|T_{m-\IP_2}((f_{i,j,S})_{i,j\in [m],S\subseteq [m]^2})|\leq \min_{i,j\in [m],S\subseteq [m]^2}\|f_{i,j,S}\|_{U^3}.\]
\end{lemma}

In order to be able to deal with \emph{non-}quadratically uniform sets, we need to understand the structure of functions whose $U^3$ norm is bounded away from zero. One class of functions whose $U^3$ norm is non-negligible is that of \emph{quadratic phase functions}, that is, functions of the form $\omega^{q(x)}$ for some quadratic form $q$ on $\F_p^n$. Indeed, it is straightforward to check that the $U^3$ norm of such functions always equals 1 (see the proof of Lemma \ref{lem:quadvc2} and the discussion of cube auto-completion in \cite[Section 3]{Terry.2021a}/\cite[Section 5]{Terry.2021av2} for the crucial property satisfied by quadratic forms). It would not be unreasonable to conjecture that functions with large $U^3$ norm must always correlate with such ``quadratic phases", in the same way that, as we saw earlier, functions with large $U^2$ norm correlate with linear phases.

Building heavily on prior work of Gowers, such an \emph{inverse theorem} for the $U^3$ norm in $\F_p^n$ was first proved by Green and Tao \cite{Green.2008piq} for characteristic $p>2$, and by Samorodnitsky \cite{Samorodnitsky.2006} in characteristic $p=2$.

\begin{theorem}[Inverse theorem for $U^3$ in $\F_p^n$]\label{thm:u3inv}
Let $f:\F_p^n\ra\C$ be a function satisfying $\|f\|_\infty \leq 1$. Then for all $\eta>0$, there is a constant $c(\eta)$ (tending to zero as $\eta$ tends to zero) such that the following holds. If $\|f\|_{U^3}\geq \eta$, then there exists a quadratic form $q$ on $\F_p^n$ such that $|\E_x f(x)\omega^{q(x)}|\geq c(\eta)$.
\end{theorem}

In arithmetic combinatorics and theoretical computer science, the quantitative form of the dependence of $c(\eta)$ on $\eta$ was an important open question for close to two decades. Indeed, whether or not it is polynomial had been shown to be equivalent to the \emph{Polynomial Freiman-Ruzsa Conjecture} \cite{Green.2010d5}. The very recent resolution of this conjecture in general characteristic $p$ due to Gowers, Green, Manners, and Tao \cite{Gowers.2024} implies that we may indeed take $c(\eta)$ to be a polynomial in $\eta$.

For completeness, we mention that corresponding inverse results are also known in other abelian groups, in particular the cyclic groups \cite{Green.2008piq, Green.2011, Green.2011r7c, Green.2012xge} but also in locally compact abelian groups \cite{Candela.2019}. The first quantitative bounds in the higher-order cyclic case were only relatively recently obtained by Manners \cite{Manners.2018} (see also \cite{Bloom.2020}).

\subsection{A quadratic arithmetic regularity lemma}\label{subsec:quadarl}
The quadratic structure-randomness dichotomy laid out above gives rise to a ``quadratic" regularity lemma. In order to state this rigorously, we shall need to formally define quadratic factors and quadratic atoms.

\begin{definition}[Quadratic factor]\label{def:quadfac}
A \emph{quadratic factor $\calB$ in $\F_p^n$ of complexity $(\ell,q)$} is a pair $(\calL,\calQ)$, where $\calL=\{r_1,\ldots, r_{\ell}\}$ is a linear factor on $\F_p^n$ and $\calQ=\{M_1,\ldots, M_q\}$ is a set of symmetric $n\times n$ matrices over $\F_p$.  The \emph{atoms of $\calB$} are sets of the form
$$
\{x\in \F_p^n: x^T r_1=a_1,\ldots, x^T r_{\ell}=a_{\ell}, x^TM_1x=b_1,\ldots, x^TM_qx=b_q\},
$$
for some $a_1,\ldots, a_{\ell}, b_1,\ldots, b_q\in \F_p$.  We let $\At(\calB)$ denote the set of atoms of $\calB$.\end{definition}

When $\calB$ is a quadratic factor of complexity $(\ell,q)$, $\At(\calB)$ forms a partition of $\F_p^n$ with at most $p^{\ell+q}$ parts.  Computing the size of an atom in a quadratic factor is more subtle than in the purely linear case.  The following notion plays a key role.

\begin{definition}[Rank of a factor]\label{def:rankfac}
Suppose $\calB=(\calL,\calQ)$ is a quadratic factor on $\F_p^n$ of complexity $(\ell,q)$, with $\calL=\{r_1,\ldots, r_{\ell}\}$ and $\calQ=\{M_1,\ldots, M_q\}$.  The \emph{rank of $\calB$} is defined to be
$$
\min\{\rk(\lambda_1M_1+\ldots+\lambda_qM_q): \lambda_1,\ldots, \lambda_q\in \F_p\text{ are not all }0\}.
$$
\end{definition}

The following well-known elementary exponential sum estimate (see e.g. \cite[Lemma 3.1]{Green.2007}) relates the rank of a matrix to the magnitude of the corresponding quadratic exponential sum.

\begin{fact}[Quadratic exponential sums of high rank are small]\label{fct:expsum}
Let $M$ be a symmetric $n\times n$ matrix over $\F_p$ of rank $r$, and let $b\in\F_p^n$. Then
\[\big|\E_x \omega^{x^TMx+b^Tx}\big|\leq p^{-r/2}.\]
\end{fact}

The next lemma now provides a key observation, namely that when a factor has high rank, all its atoms have roughly the same size (see \cite[Lemma 4.2]{Green.2007}).

\begin{lemma}\label{lem:sizeofatoms}
Let $\calB=(\calL,\calQ)$ be a quadratic factor on $\mathbb{F}_p^n$ of complexity $(\ell,q)$ and rank $r$.  Then for any $B\in \At(\calB)$,  
\[|B|=(1+O(p^{\ell+q-r/2}))p^{n-\ell-q}.\]
\end{lemma}

 The following immediate corollary of Lemma \ref{lem:sizeofatoms} will be used throughout the paper.

\begin{corollary}[Size of a quadratic atom]\label{cor:sizeofatoms}
For any function $\tau:\mathbb{\R}^+\rightarrow \mathbb{R}$, there is a growth function $\rho=\rho(\tau)$ such that if $\calB=(\calL,\calQ)$ is a quadratic factor on $\mathbb{F}_p^n$ of complexity $(\ell,q)$ and rank at least $\rho(\ell+q)$, then for any $B\in \At(\calB)$, $\big||B|/p^n- p^{-(\ell+q)}\big|\leq p^{-\tau(\ell+q)}p^{-(\ell+q)}$.
\end{corollary}
When $\tau$ is the constant function $\tau(x)=\log_p(\e^{-1})$, we will write $\rho(\e)$ instead of $\rho(\tau)$.  We will also use the following size estimate on bilinear level sets, which is an   immediate corollary of Lemma \ref{lem:bilsums} (iii) in the appendix.

\begin{corollary}[Size of a bilinear level set]\label{cor:sizeofbilsets}
For any function $\tau:\mathbb{\R}^+\rightarrow \mathbb{R}$, there is a growth function $\rho=\rho(\tau)$ such that if $\calB=(\calL,\calQ)$ is a quadratic factor on $\F_p^n$ of complexity $(\ell,q)$ and rank at least $\rho(\ell+q)$ with $\calQ=\{M_1,\dots,M_q\}$, then for any $b=(b_1,\ldots, b_q)\in\F_p^q$, the bilinear level set
$$
\beta_\calQ(b)=\{(x,y)\in \F_p^n\times \F_p^n: \text{ for each }j\in [q], x^TM_jy=b_j\}
$$
 has size $\big||\beta_{\calQ}(b)|p^{-2n}- p^{-q}\big|\leq p^{-\tau(\ell+q)}p^{-q}$.
\end{corollary}

As with Corollary \ref{cor:sizeofatoms}, when $\tau$ is the constant function $\tau(x)=\log_p(\e^{-1})$, we will write $\rho(\e)$ instead of $\rho(\tau)$.

It turns out that any factor can be refined to have high rank by changing a bounded number of its constraints (see \cite[Lemma 3.11]{Green.2007}), though we shall not need to do so here.

We next state the arithmetic regularity lemma for the $U^3$ norm, which is due to Green and Tao \cite[Proposition 3.9]{Green.2007} and follows from an iterated application of Theorem \ref{thm:u3inv} and \cite[Lemma 3.11]{Green.2007}. Before doing so we need some notation. As with a linear factor, given a quadratic factor $\calB=(\calL,\calQ)$, we write $\E(f|\mathcal{\calB})$ for the projection onto $\calB$. In other words, $\E(f|\mathcal{\calB})$ is the function whose value at $x$ is the average of $f$ over the atom in $\calB$ that contains $x$. That is, given $x\in \F_p^n$, $\E(f|\mathcal{\calB})(x)=\frac{1}{|B|}\sum_{y\in B}f(y)$, where $B$ is the element of $\At(\calB)$ containing $x$. Finally, given two quadratic factors $\calB$ and $\calB'$, we say that $\calB'$ \emph{refines} $\calB$, denoted $\calB'\preceq \calB$, if the partition $\At(\calB')$ of $\F_p^n$ refines the partition $\At(\calB)$ of $\F_p^n$.\footnote{This is also known as \emph{semantic refinement}, see \cite[Definition 3.9]{Bhattacharyya.20126s9}, in contrast to \emph{syntactic refinement}.}

\begin{theorem}[$U^3$ arithmetic regularity lemma]\label{thm:gtarl}
For every $\e>0$ and for all growth functions $\rho, \sigma:\R^+\rightarrow\R^+$ (which may depend on $\e$), there exists $n_0=n_0(\e,\rho,\sigma)$ such that the following holds.
Let $n>n_0$ and let $f:\F_p^n \ra [0,1]$ be a function. Let $(\calL',\calQ')$ be a quadratic factor of complexity $(\ell', q')$. Then there exists $C(\e,\rho,\sigma, \ell', q')$ and a quadratic factor $\calB=(\calL,\calQ)$ of complexity $(\ell, q)$ such that
\begin{enumerate}[label=\normalfont(\roman*)]
\item $(\calL,\calQ)$ refines $(\calL',\calQ')$;
\item the complexity of $(\calL,\calQ)$ satisfies $\ell, q\leq C(\e,\rho,\sigma, \ell', q')$;
\item the rank of $(\calL,\calQ)$ is at least $\rho(\ell+q)$;
\item we can write 
\[f=f_1+f_2+f_3,\]
with $f_1=\E(f|\mathcal{\calB})$, $\|f_2\|_{U^3} \leq 1/\sigma(\ell+q) $, and $\|f_3\|_2\leq \e$.
\end{enumerate}
\end{theorem}

This can be done whilst ensuring that $f_1+f_3$ maps into $[0,1]$ and $f_2$ and $f_3$ map into $[-1,1]$ (see, for example, \cite[Theorem 5.1]{Bhattacharyya.20126s9} or \cite[Theorem 4.1]{Gowers.2011}). Applied to the balanced function of a set $A\subseteq \F_p^n$, this will imply that $A$ is locally (on almost all atoms) well distributed in the sense of $U^3$ (we will make this statement precise in Section \ref{sec:mainproof}).

Finally, we set up some more notation to help us talk about specific atoms more easily.  Suppose $\calB=(\calL,\calQ)$ is a quadratic factor, where  $\calL=\{r_1,\ldots, r_{\ell}\}$ and $\calQ=\{M_1,\ldots, M_q\}$.  Given $a=(a_1,\ldots, a_{\ell})\in \F_p^\ell$, $b=(b_{1},\ldots, b_{q})\in \mathbb{F}_p^{q}$, and $d=(a,b)$, we write 
\begin{align*}
L(a)&=\{x\in \mathbb{F}_p^n: x^Tr_i=a_i\text{ for each }i\in [\ell]\},\\
Q(b)&=\{x\in \mathbb{F}_p^n: x^TM_{j}x=b_{j}\text{ for each }j\in [q]\},\text{ and }\\
B(d)&=L(a)\cap Q(b).
\end{align*}

To conclude this brief overview of higher-order Fourier analysis, some further remarks are in order. First, there are numerous generalizations of the above regularity lemma, both to higher-order uniformity norms, for instance \cite{Gowers.2011, Bhattacharyya.20126s9}, as well as other groups, e.g. \cite{Green.2010jwua,Candela.2019}. 

All these results are obtained iteratively from an inverse theorem (such as Theorem \ref{thm:u3inv}), and therefore the bounds in the arithmetic regularity lemma necessarily depend on those in the inverse theorem as discussed at the end of Section \ref{subsec:uknorms2}. Having said this, regardless of whether the dependence in Theorem \ref{thm:u3inv} is exponential or polynomial (as is now known for $\F_p^n$), for a choice of the rank function $\rho$ common in applications (namely of polynomial form), the bound on the complexity of the factor in Theorem \ref{thm:gtarl} is of Tower type (see for instance \cite[Appendix B]{Gladkova.2025}). We refer the reader to our forthcoming paper \cite{Terry.2024f} for further discussion.

\section{The local Gowers semi-norms and local independence property operators}\label{sec:localnorms}

In this section we define local versions of the Gowers norms and independence property operators encountered in the preceding section. We begin, in Section \ref{subsec:localu2}, by localising the $U^2$ norm and the $\IP$ operator. As a first application of the former, we prove that quadratic atoms of high-rank factors have small local $U^2$ semi-norm. We then show that the local $\IP$-operator is controlled by the local $U^2$ semi-norm. 

In Section \ref{subsec:exaQA}, we go on to show that quadratic atoms of high-rank factors have unbounded $\VC$-dimension. We also establish that they have bounded $\VC_2$-dimension.

In Section \ref{subsec:localu3}, we then set up local versions of the $U^3$ norm and the $\IP_2$ operator, and show that the former controls the latter. This will be helpful in Section \ref{subsec:translating}, where we reformulate the $U^3$ arithmetic regularity lemma to say that for every subset $A\subseteq \F_p^n$, there exists a high-rank, bounded-complexity quadratic factor with the property that $A$ is locally quadratically uniform on most atoms.

\subsection{The local $U^2$ semi-norm and the local $\IP$-operator}\label{subsec:localu2}

We begin by defining the local $U^2$ inner product, relative to an atom of a linear factor $\calL$.

\begin{definition}[Local $U^2$ inner product]\label{def:localipu2}
Given a linear factor $\calL$ of complexity $\ell$, a tuple $d=(a_1,a_2)\in \mathbb{F}_p^{\ell}\times\F_p^\ell$, and functions $(f_{\e})_{\e\in \{0,1\}^2}:\F_p^n\ra \C$, we define
\[\langle (f_{\e})_{\e\in \{0,1\}^2}\rangle_{U^2(d)}=\E_{\begin{subarray}{l} x_0,x_1\in L(a_1)\end{subarray}}\E_{\begin{subarray}{l}y_0,y_1\in L(a_2)\end{subarray}}\prod_{\e\in \{0,1\}^2}\calC^{|\e|}f_{\e}(x_{\e(1)}+y_{\e(2)}).\]
\end{definition}

At first sight this inner product appears to involve two atoms, $L(a_1)$ and $L(a_2)$, but we note that each $f_\omega$ will only ever be evaluated on $L(a_1+a_2)$. Indeed, given $d=(a_1,a_2)\in \F_p^{\ell}\times \F_p^{\ell}$, we observe that the same $U^2(d)$ inner product arises from all pairs $(a_1',a_2')$ such that $a_1'+a_2'=a_1+a_2$. For this reason, the following notation will be useful. 

\begin{notation}[$\Sigma_\calL(d)$]
Given a linear factor $\calL$ of complexity $\ell$ and $d=(a_1,a_2)\in \F_p^{\ell}\times \F_p^{\ell}$, let $\Sigma_{\calL}(d)=a_1+a_2$.  We will drop the subscript when $\calL$ is clear from context.
\end{notation}

We see that the $U^2(d)$ norm in Definition \ref{def:localipu2} depends only on $\Sigma(d)$, rather than $d$. However, this distinction will only matter in the higher-order case. We next define the local $U^2$ semi-norm in the obvious way.

\begin{definition}[Local $U^2$ semi-norm]\label{def:localu2}
Given a linear factor $\calL$ of complexity $\ell$, a tuple $d=(a_1,a_2)\in \mathbb{F}_p^{\ell}\times \F_p^\ell$, and $f:\F_p^n\ra\C$, we define
\[\|f\|_{U^2(d)}= \langle (f)_{\e\in \{0,1\}^2}\rangle_{U^2(d)}^{1/4}.\]
\end{definition}

We need to show that this is well-defined, and indeed a semi-norm.

\begin{lemma}\label{lem:checknorm}
Given a linear factor $\calL$ of complexity $\ell$ and $d=(a_1,a_2)\in \mathbb{F}_p^{\ell}\times \F_p^\ell$, $U^2(d)$ defines a semi-norm on the space of functions from $\F_p^n$ to $\C$, and a norm on the space of functions from $\F_p^n$ to $\C$ supported on $L(a_1+a_2)$.
\end{lemma}
\begin{proof} Let $H=L(0)$ be the atom of $\calL$ labelled by $0$. Note that for $f:\F_p^n\ra\C$ and any fixed $z\in L(a_1+a_2)=z+H$, we may write
\[\|f\|_{U^2(d)}^4=\E_{\begin{subarray}{l} x_0,x_1\in H\end{subarray}}\E_{\begin{subarray}{l}y_0,y_1\in H\end{subarray}}\prod_{\e\in \{0,1\}^2}\calC^{|\e|}f(z+x_{\e(1)}+y_{\e(2)}).\]
Denoting by $f^z$ the function $x\mapsto f(z+x)$ and by $f^z|_H$ the restriction of $f^z$ to $H$, the above can be expressed as
\[\E_{\begin{subarray}{l} x_0,x_1\in H\end{subarray}}\E_{\begin{subarray}{l}y_0,y_1\in H\end{subarray}}\prod_{\e\in \{0,1\}^2}\calC^{|\e|}f^{z}|_{H}(x_{\omega(1)}+y_{\omega(2)}),\]
which, by Fact \ref{fact:fourierequiv}, equals 
\[\sum_{t\in \widehat{H}} |\widehat{f^z |_{H}}(t)|^4,\]
where the Fourier transform is taken relative to $H$ as defined in \cite[Section 2]{Green.2005} (see the remark following Theorem \ref{thm:linearreg}). Since 
\[\|f\|_{U^2(d)}=\big(\sum_{t\in \widehat{H}} |\widehat{f^z |_{H}}(t)|^4\big)^{1/4},\]
and $\ell^4$ is a norm on the space of functions from $H$ to $\C$, $U(d)$ is a semi-norm on the space of functions from $\F_p^n$ to $\C$, and a norm on the space of functions from $\F_p^n$ to $\C$ supported on $L(a_1+a_2)$.\footnote{Alternatively, one may observe that $U^2(d)$ is clearly homogeneous, real-valued and non-negative. It is not difficult to prove a Gowers-Cauchy-Schwarz type inequality (see the more complicated Proposition \ref{prop:locgcs}), which implies that $U^2(d)$ satisfies the triangle inequality.}
\end{proof}

Similar localised or `directional' norms have previously appeared in \cite{Tao.2016} and subsequent quantitative work on polynomial patterns in the integers and the primes, and in the context of $\F_p^n$ most notably and most recently in \cite{Milicevic.2021, Milicevic.2024}.

We shall also need the following standard lemma, 
which shows that when a linear map  is restricted to a subspace, its rank decreases by at most twice the codimension of the subspace. It follows, for instance, from Sylvester's inequality.

\begin{fact}[Rank drop on subspaces is bounded]\label{fct:ranklem}
Let $M$ be a symmetric $n\times n$ matrix of rank $\rho$ on $\F_p^n$, and let $W\leqslant \F_p^n$ be a subspace of codimension $\ell$. Then the rank of the restriction of $M$ to $W$ is at least $\rho-\ell$.
\end{fact}

With these basic facts in hand, we prove that quadratic atoms of high-rank factors have vanishingly small local $U^2$ semi-norm.  

\begin{lemma}[Quadratic atoms of high-rank factors have small local $U^2$ semi-norm]\label{lem:quadatuni}
Let $\rho>0$, let $(\calL,\calQ)$ be a quadratic factor on $\F_p^n$ of complexity $(\ell,q)$ and rank at least $\rho$, and let $A$ be an atom of $(\calL,\calQ)$. Then for any fixed $d=(a_1,a_2)\in \F_p^{\ell}\times \F_p^\ell$, 
\[\|1_A-\alpha_{L(\Sigma(d))}\|_{U^2(d)}=O(p^{-(\rho-\ell)/4}),\]
where $\alpha_{L(\Sigma(d))}=|A\cap L(\Sigma(d))|/|L(\Sigma(d))|$ and $\Sigma(d)=a_1+a_2$.
\end{lemma}

\begin{proof} 
An atom $A$ of $(\calL,\calQ)$ is of the form 
$$
A=\{x\in \F_p^n: x^TM_jx=f_j,x^Tr_i=e_i \mbox{ for all }j\in[q], i\in [\ell]\},
$$
 for some symmetric $n\times n$ matrices $M_j$, vectors $r_i \in \F_p^n$, and $e_i,f_j\in \F_p$. We proceed by considering two cases.  First, if $e=(e_1,\dots,e_\ell)\neq a_1+a_2$, then $\alpha_{L(\Sigma(d))}=0$, and further 
\[\|1_A-\alpha_{L(\Sigma(d))}\|_{U^2(d)}=0.\]

Suppose instead that $e=a_1+a_2$. Let $H=L(0)$ be the atom of $\calL$ labelled by $0$. For the remainder of the proof, we shall continue to use the notation introduced in the proof of Lemma \ref{lem:checknorm}, with $f^z$ denoting the function $x\mapsto f(x+z)$, and $f|_H$ denoting the restriction of $f$ to $H$. 

For any $z\in L(a_1+a_2)=z+H$, we may use (\ref{starstar}) to express the indicator function $1_A$ shifted by $z$ and restricted to $H$ as
\begin{align*}
1_A^z|_{H}(x)=\prod_{j\in [q]}\E_{v\in \F_p}\omega^{v( (x+z)^TM_j(x+z)-f_j) }=\E_{v_j\in \F_p:j\in [q]}\omega^{ (x+z)^T(\sum_j v_jM_j)(x+z)-\sum_j v_j f_j }.
\end{align*}
For fixed $z\in L(a_1+a_2)$, the Fourier coefficient of $1_A^z|_{H}$ relative to $H$ at $t\in \widehat{H}$  can be therefore be written as
\begin{align}\label{eq:fc}
\widehat{1_A^z|_{H}}(t)=\E_{v_j\in \F_p:j\in [q]}\omega^{-\sum_jv_jf_j+z^T\sum_jv_j M_j z}\E_{x\in H}\omega^{ x^T(\sum_j v_j M_j)x +x^T(2\sum_j v_j M_jz-t)}.
\end{align}
By the assumption on the rank of the factor along with Facts \ref{fct:ranklem} and \ref{fct:expsum}, the quadratic exponential sum in $x$ in (\ref{eq:fc}) is bounded by $O(p^{-(\rho-\ell)/2})$ unless $v_j=0$ for all $j\in [q]$, implying that 
\[\big|\widehat{1_A^z|_{H}}(t)\big|=\frac{1}{p^q}\big|\E_{x\in H}\omega^{-x^Tt}\big|+O(p^{-(\rho-\ell)/2}).\]
Further, by (\ref{star}), the remaining linear exponential sum in $x$ is zero unless $t=0\in \widehat{H}$, yielding
\[\big|\widehat{1_A^z|_{H}}(t)\big|=\frac{1}{p^q}1_{\{0\}}(t)+O(p^{-(\rho-\ell)/2}).\]
Thus, setting $f=1_A-\alpha_{L(\Sigma(d))}$, and noting that $\alpha_{L(\Sigma(d))}=\widehat{1_A^z|_{H}}(0)$, we have 
\begin{equation}\label{eq:est}
\big|\widehat{f^z|_{H}}(t)\big|=O(p^{-(\rho-\ell)/2})
\end{equation}
for all $t\in \widehat{H}$. Now 
\[\|1_A-\alpha_{L(\Sigma(d))}\|_{U^2(d)}^4=\sum_{t\in \widehat{H}} |\widehat{f^z |_{H}}(t)|^4\leq \sup_{t\in \widehat{H}} |\widehat{f^z |_{H}}(t)|^2\sum_{t\in \widehat{H}} |\widehat{f^z |_{H}}(t)|^2=O(p^{-(\rho-\ell)}),\]
where the first equality follows from Fact \ref{fact:fourierequiv}  and the final equality from Parseval's identity (see the discussion following Definition \ref{def:ft}), the fact that $f^z|_{H}$ is 1-bounded, and the estimate (\ref{eq:est}).
\end{proof}

We shall be interested in instances of the independence property (Definitions \ref{def:vcdim} and \ref{def:ipop}) whose individual summands are constrained to lie in pre-specified cosets of a subgroup.

\begin{definition}[Local $\IP$-operator]\label{def:localipop}
Given $m\geq 1$, a linear factor $\calL$ on $\F_p^n$ of complexity $\ell$, $d=(a_1,a_2)\in \mathbb{F}_p^{\ell}\times\F_p^\ell$, and functions $(f_{i,S})_{i\in [m],S\subseteq [m]}:\F_p^n\ra \C$, define

\[T_{m-\IP(d)}((f_{i,S})_{i\in [m],S\subseteq [m]})=\E_{\begin{subarray}{l} x_i\in L(a_1) \\ i\in[m]\end{subarray}}\E_{\begin{subarray}{l}y_S\in L(a_2)\\S\subseteq[m] \end{subarray}}\prod_{i\in [m],S\subseteq [m]}f_{i,S}(x_i+y_S).\]
\end{definition}

 Analogous conventions to those introduced in Notations \ref{not:ipop} and \ref{not:ip2op} remain in force.  We are now able to state a local version of Lemma \ref{lem:ipcontrol}, which asserted that the $\IP$-operator is controlled by the $U^2$ norm.

\begin{lemma}[Local $\IP$ is controlled by local $U^2$]\label{lem:iploccontrol}
Let $m\geq 2$, let $\calL$ be a linear factor on $\F_p^n$ of complexity $\ell$, and let $d=(a_1,a_2)\in \mathbb{F}_p^{\ell}\times \F_p^\ell$.
Suppose that for each $i\in [m]$ and $S\subseteq [m]$, $f_{i,S}:\F_p^n\ra \C$ is such that $\|f_{i,S}\|_{\infty}\leq 1$. Then 
\[|T_{m-\IP(d)}((f_{i,S})_{i\in [m],S\subseteq [m]})|\leq \min_{i\in [m],S\subseteq [m]}\|f_{i,S}\|_{U^2(d)}.\]
\end{lemma}

\begin{proof} Inspecting the proof of Lemma \ref{lem:ipcontrol}, we observe that the constraints inherent in the definition of $T_{m-\IP(d)}$ which force $x_{I}$ and $y_{\calS}$ to lie in particular atoms are simply carried through.\end{proof}

Full details of a much more sophisticated version of this argument are provided in Appendix \ref{app:counting}) (see the proof of Proposition \ref{prop:ip2loccontrol}).

\subsection{Interlude: Properties of quadratic atoms}\label{subsec:exaQA}

As an application of the concepts laid out in the preceding section, we now prove\footnote{In \cite[Example 5.1]{Terry.2021av2} we gave a different proof, which used graph regularity. The present argument using the local $U^2$ semi-norm is more direct.} that an atom of a high-rank quadratic factor has large $\VC$-dimension. Note that this is in contrast to an atom of a linear factor, which--being a coset of a subgroup--has $\VC$-dimension $1$. The proof will serve as a blueprint for later arguments (see the proof of Proposition \ref{prop:trivdense}).  Below and throughout, we will write $A^C$ to mean $\F_p^n\setminus A$.
  
\begin{lemma}[Quadratic atoms of high-rank factors have large $\VC$-dimension]\label{lem:quadvc}
For all $m\geq 1$, there exists a growth function $\rho_0=\rho_0(m):\R^+\ra\R^+$ such that  for all growth functions $\rho:\R^+\ra\R^+$ satisfying $\rho\geq \rho_0$ the following holds.

Let $n\geq \ell\geq 1$ be integers and let $(\calL,\calQ)$ be a quadratic factor on $\F_p^n$ of complexity $(\ell,q)$ and rank at least $\rho(\ell+q)$. Then any atom of $(\calL,\calQ)$ has $\VC$-dimension at least $m$. 
\end{lemma}

\begin{proof} Fix an integer $m\geq 1$,  let $\rho_0(x)=C m2^m x$ for some constant $C>0$ to be determined, and let $\rho$ be a growth function satisfying $\rho\geq \rho_0$. Let $n\geq \ell\geq 1$, and suppose $(\calL,\calQ)$ is a quadratic factor of complexity $(\ell,q)$ and rank at least $\rho(\ell+q)$.  Let $A$ be an atom of $(\calL,\calQ)$ labelled by $(e,f)$. Pick $d=(a_1,a_2)\in \F_p^{\ell}\times\F_p^\ell$ such that $e=a_1+a_2$. Let $f_A^{d}=1_A-\alpha_{L(\Sigma(d))}$ and $f_{A^C}^d=1_{A^C}-\overline{\alpha}_{L(\Sigma(d))}$, where $\overline{\alpha}_{L(\Sigma(d))}$ is the density of $A^C$ on $L(\Sigma(d))$.  Observe that by Lemma \ref{lem:quadatuni},
\[\|f_A^{d}\|_{U^2(d)}=O(p^{-(\rho-2\ell)/4}).\]
Recalling Notation \ref{not:ipop}, consider 
\[T_{m-\IP(d)}(1_A |1_{A^C})=T_{m-\IP(d)}(f_{A}^{d}+\alpha_{L(\Sigma(d))} |f_{A^C}^{d}+\overline{\alpha}_{L(\Sigma(d))}).\]
By linearity, this equals 
\[\alpha_{L(\Sigma(d))}^{m 2^{m-1}} (1-\alpha_{L(\Sigma(d))})^{m 2^{m-1}}\]
plus $m 2^{m} -1$  terms of the form
\begin{equation}\label{eq:errorlocip}
T_{m-\IP(d)}((f_{i,S})_{i\in [m],S\subseteq [m]}),
\end{equation}
in which at least one of the input functions $f_{i,S}$ equals $f_A^{d}$ or $f_{A^C}^{d}=-f_A^{d}$. Thus, by Lemma \ref{lem:iploccontrol}, each term of the form (\ref{eq:errorlocip}) is $O(p^{-(\rho-2\ell)/4})$. We also showed in the proof of Lemma \ref{lem:quadatuni} that $\alpha_{L(\Sigma(d))}=p^{-q}+O(p^{-(\rho-2\ell)/2})$, and therefore 
\[T_{m-\IP(d)}(1_A |1_{A^C})\geq (p^{-q})^{m2^{m}} \big(1+O(p^{q-(\rho-2\ell)/2})\big)^{m2^{m}}-m2^m O(p^{-(\rho-2\ell)/4}).\]
Thus, since $\rho$ is a sufficiently fast-growing function of $\ell+q$,  the above difference is positive, and thus $A$ contains an $m$-$\IP$. 
\end{proof}

In contrast to the preceding lemma, we now show that quadratic atoms are well-behaved in terms of their $\VC_2$-dimension (see Definition \ref{def:vc2dim}). This follows from a well-known identity relating the eight corners of the Gowers cube defining the $U^3$ norm.  For a more comprehensive discussion, see \cite[Section 3.2]{Terry.2021a}.

\begin{lemma}[Quadratic atoms have bounded $\VC_2$-dimension]\label{lem:quadvc2}
Let $(\calL,\calQ)$ be a quadratic factor. Then any atom of $(\calL,\calQ)$ has $\VC_2$-dimension at most $1$.
\end{lemma}

\begin{proof} For any $x_0, x_1, y_0, y_1, z_0, z_1 \in \F_p^n$, we have 
\[\sum_{\omega\in \{0,1\}^3} (-1)^{|\omega|}(x_{\omega(1)}+y_{\omega(2)}+z_{\omega(3)})^TM(x_{\omega(1)}+y_{\omega(2)}+z_{\omega(3)}) =0\]
for any $n\times n$ symmetric matrix $M$ with entries in $\F_p$, as well as
\[\sum_{\omega\in \{0,1\}^3} (-1)^{|\omega|}(x_{\omega(1)}+y_{\omega(2)}+z_{\omega(3)})^T r =0\]
for any vector $r\in \F_p^n$.
Thus, if $x_0$, $x_1$, $y_0$, $y_1$, $z_0$, $z_1$ are such that $x_{\omega(1)}+y_{\omega(2)}+z_{\omega(3)}$  lie in an atom $Q$ of $(\calL,\calQ)$ for $\omega \in \{0,1\}^3\setminus\{(0,0,0)\}$, then $x_0+y_0+z_0$ must also lie in $Q$.

Now suppose $Q$ had $2$-$\IP_2$, i.e. suppose there were $x_0$, $x_1$, $y_0$, $y_1$, and for each $S\subseteq \{0,1\}^2$, $z_S$ such that $x_i+y_j+z_S\in Q$ if and only if $(i,j)\in S$. In particular, suppose we had such  $x_0$, $x_1$, $y_0$, $y_1$ and $z_S$, $z_{S'}$ for $S=\{(0,0),(1,0),(0,1),(1,1)\}$ and $S'=\{(1,0),(0,1),(1,1)\}$. This would mean $x_0+y_0+z_S\in Q$, $x_1+y_0+z_S\in Q$, $x_0+y_1+z_S\in Q$, $x_1+y_1+z_S\in Q$, $x_1+y_0+z_{S'}\in Q$, $x_0+y_1+z_{S'}\in Q$, $x_1+y_1+z_{S'}\in Q$ but $x_0+y_0+z_{S'}\notin Q$, which is impossible by the preceding observation with $z_0=z_{S'}$ and $z_1=z_S$.
\end{proof}

\subsection{The local $U^3$ semi-norm and the local $\IP_2$-operator}\label{subsec:localu3}

We now generalise our earlier observations on the local independence property and the local $U^2$ semi-norm to the higher-order setting.

We shall begin by defining our local $U^3$ semi-norm. In contrast to Definition \ref{def:localu2}, this includes constraints on pairs of variables. This is necessary in order to ensure that $f$ is only ever being evaluated on the same quadratic atom. An additional complication arises in that even the constraints on single variables are defined by sets that are not exactly the same size (namely certain pre-specified quadratic atoms). This will be taken care of through the expectation notation, but for pairs of variables we will need to account for this more explicitly, via so-called \emph{characteristic measures}. Indeed, for a set of pairs $\beta\subseteq \F_p^n\times\F_p^n$, we write $\mu_{\beta}$ for the function on $\F_p^n\times \F_p^n$ defined by $(x,y)\mapsto 1_\beta(x,y)  |\F_p^n|^2/|\beta|$, where the normalisation has the effect that $\E_{x,y}\mu_\beta(x,y)=1$.  

In the definition below, we will apply this notation to bilinear level sets associated with a given quadratic factor. Specifically, 
given a quadratic factor $\calB=(\calL,\calQ)$ on $\F_p^n$ with $\calQ=\{M_1,\ldots, M_q\}$, and given $b=(b_1,\ldots, b_q)\in \F_p^q$, we define
$$
\beta_{\calQ}(b)=\{(x,y)\in \F_p^n\times \F_p^n: xM_jy = b_j\text{ for each }j\in [q]\}.
$$
We will drop the subscript $\calQ$ when it is clear from context.  We now give the analogue of Definition \ref{def:localipu2}.

\begin{definition}[Local $U^3$ inner product]\label{def:localu3ip}
Given a quadratic factor $\calB=(\calL,\calQ)$ on $\F_p^n$ of complexity $(\ell,q)$, $d=(a_1,a_2,a_3,b_{12},b_{13},b_{23})\in \F_p^{\ell+q}\times \F_p^{\ell+q}\times \F_p^{\ell+q}\times \F_p^q\times \F_p^q\times \F_p^q $, and functions $(f_{\omega})_{\omega\in \{0,1\}^3}:\F_p^n\ra \C$, define
\begin{align*}\langle (f_{\omega})_{\omega\in \{0,1\}^3}\rangle_{U^3(d)}=\E_{ x_0,x_1\in B(a_1)}&\E_{y_0,y_1\in B(a_2)}\E_{z_0,z_1\in B(a_3)}\\
\prod_{(i,j)\in \{0,1\}^2} &\mu_{\beta(b_{12})}(x_i,y_j)\prod_{(i,k)\in \{0,1\}^2} \mu_{\beta(b_{13})}(x_i,z_k)\prod_{(j,k)\in \{0,1\}^2} \mu_{\beta(b_{23})}(y_j,z_k)\\
&\prod_{\omega\in \{0,1\}^3}\calC^{|\omega|}f_{\omega}(x_{\omega(1)}+y_{\omega(2)}+z_{\omega(3)}).
\end{align*}
\end{definition}

Note that the restriction of singletons and pairs does not affect the linearity properties of the operator. Observe further that in Definition \ref{def:localu3ip}, when computing the $U^3(d)$ inner product for $d=(a_1,a_2,a_3,b_{12},b_{13},b_{23})$, we will only ever evaluate $f$ on the quadratic atom with label $a_1+a_2+a_3+ 2(0b_{12})+2(0b_{13})+2(0b_{23})\in \F_p^{\ell+q}$, where each of the vectors $0b_{ij}$ is $b_{ij}\in\F_p^q$ augmented by $\ell$ initial zeros.   For this reason, we introduce the following notation.

\begin{notation}[$\Sigma_{\calB}(d)$]
Given a quadratic factor $\calB=(\calL,\calQ)$ of complexity $(\ell,q)$ and $d=(a_1,a_2,a_3,b_{12},b_{13},b_{23})\in \F_p^{\ell+q}\times \F_p^{\ell+q}\times \F_p^{\ell+q}\times\F_p^q\times \F_p^q\times \F_p^q$, let 
$$
\Sigma_{\calB}(d)=a_1+a_2+a_3+ 2(0b_{12})+2(0b_{13})+2(0b_{23}),
$$
where each of the vectors $0b_{ij}\in \F_p^{\ell+q}$ consists of $b_{ij}\in\F_p^q$ augmented by $\ell$ initial zeros.  When $\calB$ is clear from context we drop the subscript on $\Sigma$.
\end{notation}

Observe, however, that unlike the linear analogue of this definition (Definition \ref{def:localipu2}), the local $U^3$ inner product (and thus the local $U^3$ semi-norm below that it gives rise to) genuinely depend on the tuple $d=(a_1,a_2,a_3,b_{12},b_{13},b_{23})$ rather than just $\Sigma(d)$.  For instance, the normalization occurring in Definition \ref{def:localu3ip} depends on the specific choice of tuple $d$, since quadratic atoms and bilinear sets with different labels may have slightly different sizes.  Moreover, given an element $g$ in the atom labelled by $\Sigma(d)$, the number of ways to write it as a sum of the form $x+y+z$, where $x,y,z$ and $(x,y),(x,z), (y,z)$ lie in the atoms and bilinear sets specified by $d$ depends on the choice of $d$.

\begin{definition}[Local $U^3$ semi-norm]\label{def:localu3}
Given a quadratic factor $\calB=(\calL,\calQ)$ of complexity $(\ell,q)$, $d=(a_1,a_2,a_3,b_{12},b_{13},b_{23})\in \F_p^{\ell+q}\times \F_p^{\ell+q}\times \F_p^{\ell+q}\times\F_p^q\times \F_p^q\times \F_p^q$ and a function $f:\F_p^n\ra \C$, we define
\[\|f\|_{U^3(d)}=\langle (f)_{\e\in \{0,1\}^3}\rangle_{U^3(d)}^{1/8}.\]
\end{definition}

To see that $U^3(d)$ is well defined, observe that 
\begin{align*}\langle (f)_{\omega\in \{0,1\}^3}\rangle_{U^3(d)}=\E_{ x_0,x_1\in B(a_1)}&\E_{y_0,y_1\in B(a_2)}\E_{z_0,z_1\in B(a_3)}\\
\prod_{(i,j)\in \{0,1\}^2} &\mu_{\beta(b_{12})}(x_i,y_j)\prod_{(i,k)\in \{0,1\}^2} \mu_{\beta(b_{13})}(x_i,z_k)\prod_{(j,k)\in \{0,1\}^2} \mu_{\beta(b_{23})}(y_j,z_k)\\
&\prod_{\omega\in \{0,1\}^3}\calC^{|\omega|}f(x_{\omega(1)}+y_{\omega(2)}+z_{\omega(3)}).
\end{align*}
can be written as
\begin{align*}\E_{ x_0,x_1\in B(a_1)}&\E_{y_0,y_1\in B(a_2)}\prod_{(i,j)\in \{0,1\}^2} \mu_{\beta(b_{12})}(x_i,y_j)\\
&\big|\E_{z_0\in B(a_3)}\prod_{i\in \{0,1\}} \mu_{\beta(b_{13})}(x_i,z_0)\prod_{j\in \{0,1\}} \mu_{\beta(b_{23})}(y_j,z_0)
\prod_{\omega\in \{0,1\}^2}\calC^{|\omega|}f(x_{\omega(1)}+y_{\omega(2)}+z_{0})\big|^2,
\end{align*}
and thus $\langle (f)_{\omega\in \{0,1\}^3}\rangle_{U^3(d)}$ is non-negative. To see that $U^3(d)$ does indeed define a semi-norm, it will be useful to know that the corresponding local inner product satisfies a Gowers-Cauchy-Schwarz inequality. 

\begin{proposition}[Local Gowers-Cauchy-Schwarz]\label{prop:locgcs}
Given a quadratic factor $\calB=(\calL,\calQ)$ of complexity $(\ell,q)$, $d=(a_1,a_2,a_3,b_{12},b_{13},b_{23})\in  \F_p^{\ell+q}\times \F_p^{\ell+q}\times \F_p^{\ell+q}\times\F_p^q\times \F_p^q\times \F_p^q$ and functions $(f_{\e})_{\e\in \{0,1\}^3}:\F_p^n\ra \C$, we have
\[\big|\langle (f_{\e})_{\e\in \{0,1\}^3}\rangle_{U^3(d)}\big|\leq \prod_{\e\in \{0,1\}^3} \|f_\e\|_{U^3(d)}.\]
\end{proposition}

For the proof of Proposition \ref{prop:locgcs}, it will be convenient to have some more notation.

\begin{notation}\label{not:locu3ip} 
For $\e\in \{0,1\}^2$, $j\in \{0,1\}$ and $i\in [3]$, define $\e\wedge_i j$ to be the vector in $\{0,1\}^{3}$ satisfying
\[
(\e\wedge_i j)(k)=\begin{cases}\e(k) &\text{ if }k<i\\
\e(k-1) & \text{ if }k>i\\
j&\text{ if }k=i.
\end{cases}
\]
  In other words, $\e\wedge_i j$ has the effect of inserting $j$ in the $i$th coordinate of $\e$ and pushing the later coordinates along by 1. When $i=3$, we omit the subscript and write $\e\wedge i$ for $\e\wedge_3 i$, i.e. the vector in $\{0,1\}^3$ with coordinates $(\e(1),\e(2),i)$.

For $\e\in \{0,1\}$, $k,j\in \{0,1\}$ and $i<\ell\in [3]$, write $\e\wedge_{(i,\ell)} (j,k)$ for the vector in $\{0,1\}^{3}$ whose $i$th coordinate equals $j$, whose $\ell$th coordinate equals $k$, and whose remaining coordinate is $\e$. 

In the reverse direction, for $\e\in \{0,1\}^3$ and $\emptyset \neq S\subsetneq [3]$, let $\e|_S$ be the vector in $\{0,1\}^{|S|}$ from which all coordinates of $\e$ not indexed by an element in $S$ have been deleted.
\end{notation}

\begin{proofof}{Proposition \ref{prop:locgcs}} The proof will follow the standard Cauchy-Schwarz argument, with extra care taken to account for the bilinear constraints. We begin by rearranging the inner product $\langle (f_{\omega})_{\omega\in \{0,1\}^3}\rangle_{U^3(d)}$ as
\begin{align*}
&\E_{x_0,x_1\in B(a_1),y_0,y_1\in B(a_2)}\prod_{(i,j)\in \{0,1\}^2} \mu_{\beta(b_{12})}(x_i,y_j)^{1/2}\mu_{\beta(b_{12})}(x_i,y_j)^{1/2}\\
&\hspace{10pt}\E_{z_0\in B(a_3)}\prod_{i\in \{0,1\}} \mu_{\beta(b_{13})}(x_i,z_0)\prod_{j\in \{0,1\}}\mu_{\beta(b_{23})}(y_j,z_0)\prod_{\omega\in \{0,1\}^2}\calC^{|\omega|}f_{\omega\wedge 0}(x_{\omega(1)}+y_{\omega(2)}+z_0)\\
&\hspace{20pt}\E_{z_1\in B(a_3)}\prod_{i\in \{0,1\}} \mu_{\beta(b_{13})}(x_i,z_1)\prod_{j\in \{0,1\}}\mu_{\beta(b_{23})}(y_j,z_1)\prod_{\omega\in \{0,1\}^2}\calC^{|\omega|+1}f_{\omega\wedge 1}(x_{\omega(1)}+y_{\omega(2)}+z_1).
\end{align*}
Applying Cauchy-Schwarz in $x_0,x_1,y_0,y_1$, and distributing instances of $\mu_{\beta(b_{12})}^{1/2}(x_i,y_j)$ appropriately yields that $\big|\langle (f_{\omega})_{\omega\in \{0,1\}^3}\rangle_{U^3(d)}\big|^2$
is bounded above by
\begin{align*}\E_{ \begin{subarray}{l}x_0,x_1\in B(a_1)\\y_0,y_1\in B(a_2)\end{subarray}}&\prod_{(i,j)\in \{0,1\}^2} \mu_{\beta(b_{12})}(x_i,y_j)\\
&\big|\E_{z_0\in B(a_3)}\prod_{i\in \{0,1\}} \mu_{\beta(b_{13})}(x_i,z_0)\prod_{j\in \{0,1\}}\mu_{\beta(b_{23})}(y_j,z_0)\prod_{\omega\in \{0,1\}^2}\calC^{|\omega|}f_{\omega\wedge 0}(x_{\omega(1)}+y_{\omega(2)}+z_0)\big|^2
\end{align*}
times
\begin{align*}\E_{ \begin{subarray}{l}x_0,x_1\in B(a_1)\\y_0,y_1\in B(a_2)\end{subarray}}&\prod_{(i,j)\in \{0,1\}^2} \mu_{\beta(b_{12})}(x_i,y_j)\\
&\big|\E_{z_1\in B(a_3)}\prod_{i\in \{0,1\}} \mu_{\beta(b_{13})}(x_i,z_1)\prod_{j\in \{0,1\}}\mu_{\beta(b_{23})}(y_j,z_1)\prod_{\omega\in \{0,1\}^2}\calC^{|\omega|+1}f_{\omega\wedge 1}(x_{\omega(1)}+y_{\omega(2)}+z_1)\big|^2.
\end{align*}
Recalling Notation \ref{not:locu3ip}, it follows that 
\begin{equation}\label{eq:firststep}
\big|\langle (f_{\omega})_{\omega\in \{0,1\}^3}\rangle_{U^3(d)}\big|\leq 
\langle (f_{\omega|_{\{1,2\}}\wedge 0})_{\omega\in \{0,1\}^3}\rangle_{U^3(d)}^{1/2}
\langle (f_{\omega|_{\{1,2\}}\wedge 1})_{\omega\in \{0,1\}^3}\rangle_{U^3(d)}^{1/2}.
\end{equation}
By the exact same argument, for every $i\in [3]$,
\begin{equation}\label{eq:genstep}
\big|\langle (f_{\omega})_{\omega\in \{0,1\}^3}\rangle_{U^3(d)}\big|\leq 
\langle (f_{\omega|_{[3]\setminus\{i\}}\wedge_i 0})_{\omega\in \{0,1\}^3}\rangle_{U^3(d)}^{1/2}
\langle (f_{\omega|_{[3]\setminus\{i\}}\wedge_i 1})_{\omega\in \{0,1\}^3}\rangle_{U^3(d)}^{1/2}.
\end{equation}
Applying (\ref{eq:genstep}) with $i=2$ to the first factor on the right-hand side of (\ref{eq:firststep}), we have
\begin{equation}\label{eq:thirdstep}
\big|\langle (f_{\omega|_{\{1,2\}}\wedge 0})_{\omega\in \{0,1\}^3}\rangle_{U^3(d)}\big|^{1/2}\leq 
\langle (f_{\omega|_{\{1\}}\wedge_{(2,3)}(0, 0)})_{\omega\in \{0,1\}^3}\rangle_{U^3(d)}^{1/4}
\langle (f_{\omega|_{\{1\}}\wedge_{(2,3)}(1, 0)})_{\omega\in \{0,1\}^3}\rangle_{U^3(d)}^{1/4}.
\end{equation}
Applying (\ref{eq:genstep}) again, this time with $i=1$, the first factor on the right of ($\ref{eq:thirdstep}$) can be bounded as
\[\big|\langle (f_{\omega|_{\{1\}}\wedge_{(2,3)}(0, 0)})_{\omega\in \{0,1\}^3}\rangle_{U^3(d)}\big|^{1/4}\leq 
\langle (f_{(0,0, 0)})_{\omega\in \{0,1\}^3}\rangle_{U^3(d)}^{1/8}
\langle (f_{(1,0,0)})_{\omega\in \{0,1\}^3}\rangle_{U^3(d)}^{1/8},\]
which equals $\| f_{(0,0, 0)}\|_{U^3(d)}\|f_{(1,0, 0)}\|_{U^3(d)}$. Similarly, the second factor on the right of ($\ref{eq:thirdstep}$) can be bounded above by $\| f_{(0,1, 0)}\|_{U^3(d)}\|f_{(1,1, 0)}\|_{U^3(d)}$, implying that
\[\langle (f_{\omega|_{\{1,2\}}\wedge 0})_{\omega\in \{0,1\}^3}\rangle_{U^3(d)}^{1/2}\leq \| f_{(0,0, 0)}\|_{U^3(d)}\|f_{(1,0, 0)}\|_{U^3(d)}\| f_{(0,1, 0)}\|_{U^3(d)}\|f_{(1,1, 0)}\|_{U^3(d)}.\]
An upper bound for the second factor on the right in (\ref{eq:firststep}) is derived in an identical fashion.
\end{proofof}

It follows from Proposition \ref{prop:locgcs} in the standard way \cite[page 420]{Tao.2010} that $U^3(d)$ satisfies the triangle inequality. It is also non-negative, and clearly homogeneous, and thus, like $U^2(d)$, a semi-norm on the space of functions $\F_p^n\ra\C$.

\begin{lemma}\label{lem:locu3isnorm}
Given a quadratic factor $\calB=(\calL,\calQ)$ of complexity $(\ell,q)$ and a tuple $d=(a_1,a_2,a_3,b_{12},b_{13},b_{23})\in \F_p^{\ell+q}\times \F_p^{\ell+q}\times \F_p^{\ell+q}\times\F_p^q\times \F_p^q\times \F_p^q$,  $U^3(d)$ defines a semi-norm on the space of functions from $\F_p^n$ to $\C$.
\end{lemma}

Another well known consequence of the standard, global Gowers-Cauchy-Schwarz inequality (Lemma \ref{lem:gcs}) is that the $U^3$ norm dominates the $U^2$ norm. Proposition \ref{prop:locgcs} implies that the local $U^3$ semi-norm, when taken relative to a quadratic factor $\calB=(\calL,\calQ)$ whose quadratic part $\calQ$ is trivial, dominates the local $U^2$ semi-norm with respect to an appropriate atom of $\calL$. Alas, we shall not make use of this fact in this paper.

\begin{lemma}[Local $U^3$ dominates local $U^2$ on linear factors]\label{lem:locu3dom}
Given a quadratic factor $\calB=(\calL,\emptyset)$ of complexity $(\ell,0)$ and $(a_1,a_2,a_3)\in \F_p^{\ell}\times \F_p^{\ell}\times \F_p^{\ell}$,
\[\|f\|_{U^3(d)}\geq \|f\|_{U^2((a_1+a_2,a_3))},\]
where $d=(a_10,a_20,a_30,0,0,0)\in \F_p^{\ell+q}\times \F_p^{\ell+q}\times \F_p^{\ell+q}\times\F_p^q\times \F_p^q\times \F_p^q$ with $a_i0\in \F_p^{\ell+q}$ denoting $a_i$ with $q$ zeros appended.\end{lemma}
\begin{proof}
Applying Proposition \ref{prop:locgcs} with $\calB=(\calL,\emptyset)$ and  $d=(a_10,a_20,a_30,0,0,0)\in \F_p^{\ell+q}\times \F_p^{\ell+q}\times \F_p^{\ell+q}\times\F_p^q\times \F_p^q\times \F_p^q$ to functions $(f_{\e})_{\e\in \{0,1\}^3}:\F_p^n\ra \C$, where $f_{\e}=f$ for some $f:\F_p^n\ra \C$ when $\e(3)=0$ and $f_{\e}=1$ otherwise, we find that 
\[|\langle (f_{\e})_{\e\in \{0,1\}^3}\rangle_{U^3(d)}|\leq \|f\|_{U^3(d)}^4.\]
Observe that
\[\langle (f_{\e})_{\e\in \{0,1\}^3}\rangle_{U^3(d)}=\E_{\begin{subarray}{l} x_0,x_1\in L(a_1)\end{subarray}}\E_{\begin{subarray}{l}y_0,y_1\in L(a_2) \end{subarray}}\E_{\begin{subarray}{l}z_0\in L(a_3) \end{subarray}}\prod_{\e\in \{0,1\}^2}\calC^{|\e|}f(x_{\e(1)}+y_{\e(2)}+z_{0})\]
as the quadratic and bilinear constraints and the expectation in $z_1$ become trivial by assumption.  This equals 
\[\E_{\begin{subarray}{l}z_0\in L(a_3) \end{subarray}}\E_{\begin{subarray}{l} x_0,x_1\in L(a_1)\end{subarray}}\E_{\begin{subarray}{l}y_0,y_1\in L(a_2) \end{subarray}}\prod_{\e\in \{0,1\}^2}\calC^{|\e|}f^{z_0}(x_{\e(1)}+y_{\e(2)})=\E_{\begin{subarray}{l}z_0\in L(a_3) \end{subarray}}\|f^{z_0}\|_{U^2((a_1,a_2))}^4.\]
But as previously observed, $\|\cdot \|_{U^2((a_1,a_2))}$ depends only on $a_1+a_2$, so 
\[\|f\|_{U^3(d)}^4\geq \E_{\begin{subarray}{l}z_0\in L(a_3) \end{subarray}}\|f^{z_0}\|_{U^2((a_1,a_2))}^4 =\|f\|_{U^2(\sigma(a_1,a_2,a_3))}^4,\]
where $\sigma(a_1,a_2,a_3)=a_1+a_2+a_3$.
\end{proof}

We now define the higher-order analogue of the local $\IP$-operator from Definition \ref{def:localipop}.

\begin{definition}[Local $\IP_2$-operator]\label{def:localip2op}
Given a quadratic factor $\calB=(\calL,\calQ)$ of complexity $(\ell,q)$, $d=(a_1,a_2,a_3,b_{12},b_{13},b_{23})\in \F_p^{\ell+q}\times \F_p^{\ell+q}\times \F_p^{\ell+q}\times\F_p^q\times \F_p^q\times \F_p^q$ and functions $(f_{i,j,S})_{i,j\in [m],S\subseteq [m]^2}:\F_p^n\ra \C$, define
\begin{align*}
T_{m-\IP_2(d)}((f_{i,j,S})_{i,j\in [m],S\subseteq [m]^2})=\E_{\begin{subarray}{l}  x_i\in B(a_1)\\i\in[m]\end{subarray}}&\E_{\begin{subarray}{l}y_j\in B(a_2)\\j\in[m]\end{subarray}}\E_{\begin{subarray}{l} z_S\in B(a_3)\\S\subseteq [m]^2 \end{subarray}}\\
\prod_{i,j\in [m]}& \mu_{\beta(b_{12})}(x_i,y_j)\prod_{i\in [m], S\subseteq [m]^2}\mu_{\beta(b_{13})}(x_i,z_S)\\
&\prod_{j\in [m], S\subseteq [m]^2}\mu_{\beta(b_{23})}(y_j,z_S)\prod_{i,j\in [m], S\subseteq [m]^2} f_{i,j,S}(x_i+y_j+z_S).
\end{align*}
\end{definition}
Again, this operator is clearly linear in each input. To conclude this section, we now state a simultaneous extension of Lemmas \ref{lem:ip2control} and \ref{lem:iploccontrol}, namely that  $T_{m-\IP_2(d)}$ is controlled, at least approximately, by the local $U^3(d)$ semi-norm.

\begin{proposition}[Local $\IP_2$ is controlled by local $U^3$]\label{prop:ip2loccontrol}
Let $m\geq 2$, let $\calB=(\calL,\calQ)$ be a quadratic factor on $\F_p^n$ of complexity $(\ell,q)$ and rank at least $\tau$, and let $d=(a_1, a_2, a_3, b_{12}, b_{13}, b_{23})\in \F_p^{\ell+q}\times \F_p^{\ell+q}\times \F_p^{\ell+q}\times\F_p^q\times \F_p^q\times \F_p^q$.

Suppose that for each $i,j\in [m]$, $S\subseteq [m]^2$, $f_{i,j,S}:\F_p^n\ra \C$ is such that $\|f_{i,j,S}\|_{\infty}\leq 1$. Then 
\begin{align*}
|T_{m-\IP_2(d)}((f_{i,j,S})_{i,j\in [m],S\subseteq [m]^2})|\leq (1+O_m(p^{4m2^{m^2}(\ell+q)-\tau/2}))& \min_{i,j\in [m],S\subseteq [m]^2}\|f_{i,j,S}\|_{U^3(d)} \\
&\hspace{15pt}+O_m(p^{2m2^{m^2}(\ell+q)-\tau/128}).
\end{align*}
\end{proposition}

In particular, in factors whose rank is significantly larger than the complexity, the conclusion is asymptotically similar to that of the lower-order Lemma \ref{lem:iploccontrol}.

The proof of Proposition \ref{prop:ip2loccontrol} is given in Appendix \ref{app:counting}. It relies on a Cauchy-Schwarz argument similar to that of the proof of Lemma \ref{lem:iploccontrol}, with added complications that are dealt with similarly to Gowers's \cite[Lemma 6.7]{Gowers.2006}.

To conclude this section, we observe that locally sparse sets are locally uniform, as they should be. 

\begin{proposition}[Locally sparse implies locally uniform]\label{prop:locsparseuni}
For all $\e>0$, there exists a growth function $\rho_0=\rho_0(\e)$ such that for any growth function $\rho\geq \rho_0$   the following holds. 

Let $\calB=(\calL,\calQ)$ be a quadratic factor on $\F_p^n$ of complexity $(\ell,q)$ and rank at least $\rho(\ell+q)$, and let $d=(a_1,a_2,a_3,b_{12},b_{13},b_{23})\in  \F_p^{\ell+q}\times \F_p^{\ell+q}\times \F_p^{\ell+q}\times\F_p^q\times \F_p^q\times \F_p^q$. Given any subset $A\subseteq \F_p^n$, denote by $\alpha_{B(\Sigma(d))}$ the density of $A$ on the atom $B(\Sigma(d))$.

If $\alpha_{B(\Sigma(d))}\leq \e $, then $\|1_A-\alpha_{B(\Sigma(d))}\|_{U^{3}(d)}\leq 2\e^{1/8}$.
\end{proposition}

The same conclusion clearly holds when the set $A$ instead has density near $1$, i.e. when $\alpha_{B(\Sigma(d))}=\E_{x\in B(\Sigma(d))}1_A(x)\geq 1-\epsilon$. Indeed, one simply considers the function $1_{A^C}-\overline{\alpha}_{B(\Sigma(d))}=-(1_A-\alpha_{B(\Sigma(d))})$, where $A^C=\F_p^n\setminus A$ and $\overline{\alpha}_{B(\Sigma(d))}$ is the density of $A^C$ on $B(\Sigma(d))$.

In the linear setting, this fact is extremely straightforward (see \cite[Lemma 5]{Terry.2019}). 
Proving the quadratic analogue is more involved, and will be done in Appendix \ref{app:locsparseuni}.
This result does not play a role in the proof of the main result of this paper, but is crucial in \cite[Section 4]{Terry.2021a}. Since the proof aligns more closely with the methods of the present paper, we have chosen to include it here.

\section{Proving the structure theorem for sets of bounded $\VC_2$-dimension}\label{sec:mainproof}

In this section we prove our main theorem, namely Theorem \ref{cor:vc2}. This asserts that a subset of $\F_p^n$ of bounded $\VC_2$-dimension can be approximated by a union of atoms from a high-rank, bounded-complexity quadratic factor.

In Section \ref{subsec:translating}, we prepare the output of the general quadratic regularity lemma (Theorem \ref{thm:gtarl}) to serve our purposes. The proof of Theorem \ref{cor:vc2} will be given in Section \ref{subsec:firstproof}.

\subsection{Preparing the quadratic regularity lemma}\label{subsec:translating}

We will need to bring the quadratic regularity lemma (Theorem \ref{thm:gtarl}) into a form that will be useful to us in proving Theorem \ref{cor:vc2}. Specifically, we shall show that given a set $A\subseteq \F_p^n$, the high-rank, bounded-complexity quadratic factor produced by Theorem \ref{thm:gtarl} has the property that $A$ is locally quadratically uniform on almost all atoms.

We begin with a technical lemma\footnote{This was first shown in  \cite[Lemma 3.18]{Terry.2021av2}.} showing that functions with small $L_2$ norm are locally quadratically uniform with respect to triples of atoms from a high-rank factor.  This will be used later on to deal with the part of the decomposition produced by Theorem \ref{thm:gtarl} which has small $L_2$ norm.

\begin{lemma}\label{lem:smallpart}
For all $\e>0$ there is a polynomial growth function $\rho_0=\rho_0(\e)$ such that the following holds for all growth functions $\rho\geq \rho_0$.  

Let $n\geq 1$, let $f:\mathbb{F}_p^n\rightarrow [-1,1]$ be such that $||f||_2<\e$, and let $\calB=(\calL,\calQ)$ be a quadratic factor of complexity $(\ell,q)$ and rank at least $\rho(\ell+q)$.    

Then there is a set $\Gamma\subseteq  \F_p^{\ell+q}\times \F_p^{\ell+q}\times \F_p^{\ell+q}\times\F_p^q\times \F_p^q\times \F_p^q$ such that $|\Gamma|\geq (1-8\e)p^{3\ell+6q}$ and such that for every $d=(a_1,a_2,a_3,b_{12},b_{13},b_{23})\in \Gamma$,  
$\|f\|_{U^3(d)}<2\e^{1/16}$. 
\end{lemma}

\begin{proof}
Fix $\e>0$, let $\rho_0=\rho_0(\e)$ be a growth function that grows at least as fast as $\rho(\e)$ from Corollaries  \ref{cor:sizeofatoms} and \ref{cor:sizeofbilsets}, and which simultaneously ensures that the error term in Lemma \ref{lem:genbilsums} with $m=1$ is at most $\e$, and that the proportion of exceptional triples in Lemma \ref{lem:extuples} (i) is at most $\e p^{-q}$. Note that this means that $\rho_0$ can be taken to be a linear function in $x$ whose slope depends on $\e$.

Suppose we are given $\rho\geq \rho_0$, $n\geq 1$, $f:\F_p^n\rightarrow [-1,1]$ with $||f||_2<\e$, and a quadratic factor $\calB=(\calL,\calQ)$ on $\F_p^n$ of complexity $(\ell,q)$ and rank at least $\rho(\ell+q)$.   

Given $d=(a_1,a_2,a_3,b_{12},b_{13},b_{23})\in  \F_p^{\ell+q}\times \F_p^{\ell+q}\times \F_p^{\ell+q}\times\F_p^q\times \F_p^q\times \F_p^q=\F_p^{3\ell+6q}$, define
\begin{align*}
E(d)=\{(x_1,x_2,x_3)\in B(a_1)\times B(a_2)\times B(a_3): \mbox{ for each } 1\leq  i<j\leq 3, \;(x_i,x_j)\in \beta(b_{ij})\}.  
\end{align*}
It follows from our choice of $\rho_0$ and Lemma \ref{lem:genbilsums} with $|I|=|J|=|K|=1$ that for any $d\in \F_p^{3\ell+6q}$, 
\[(1-\e)\prod_{i<j\in [3]}|\beta(b_{ij})|\prod_{i\in [n]}|B(a_i)|p^{-6n}\leq |E(d)|\leq (1+\e)\prod_{i<j\in [3]}|\beta(b_{ij})|\prod_{i\in [n]}|B(a_i)|p^{-6n}.\]
Observe further that for every $(x,y,z)\in (\F_p^n)^3$, there is a unique $d\in \F_p^{3\ell+6q}$ such that $(x,y,z)\in E(d)$.  Therefore, $\calP=\{E(d):d\in \F_p^{3\ell+6q}\}$ is a partition of $(\F_p^n)^3$, all of whose parts have roughly equal size.  We will shortly define $\Gamma$ to be the set of $d$ such that almost all triples in $(x,y,z)\in E(d)$ satisfy $|f(x+y+z)|^2<\e^{1/2}$.  Since $||f||_2<\e$, one should expect most $d$ to be in $\Gamma$.

To this end, define for each $d=(a_1,a_2,a_3,b_{12},b_{13},b_{23})\in \F_p^{3\ell+6q}$ the set
\begin{align*}
I(d)=\{(x,y,z)\in E(d): |f(x+y+z)|^2<\e^{1/2}\},
\end{align*}
and let 
\[\Gamma=\{d\in \F_p^{3\ell+6q}: |I(d)|>(1-\e^{1/2})|E(d)|\}.\] 
We shall show that $|\Gamma|>(1-8\e)p^{3\ell+6q}$.  By assumption, $\E_{x}|f(x)|^2\leq \e^2$, and consequently $\E_{x,y,z}|f(x+y+z)|^2\leq \e^2$.    Combining these observations with the definition of $I(d)$ allows us to deduce that
\[\sum_{d\in \F_p^{3\ell+6q}} \e^{1/2}|E(d)\setminus I(d)|\leq \sum_{d\in \F_p^{3\ell+6q}}\sum_{(x,y,z)\in E(d)}|f(x+y+z)|^2= \sum_{(x,y,z)\in (\F_p^n)^3}|f(x+y+z)|^2,\]
which is at most $\e^2 (p^n)^3$.
Moreover, by definition of $\Gamma$, for all $d\in \F_p^{3\ell+6q}\setminus \Gamma$ we have that
$$ \e^{1/2}|E(d)|\leq |E(d)\setminus I(d)|.
$$
Combining the last two displayed equations yields that 
\begin{align}\label{fgamma}
|\F_p^{3\ell+6q}\setminus \Gamma| \e^{1/2}\min_{d \in \F_p^{3\ell+6q}\setminus \Gamma}|E(d)|\leq \sum_{d\in \F_p^{3\ell+6q}\setminus \Gamma}|E(d)\setminus I(d)|\leq \e^{3/2}(p^n)^3.
\end{align}
By our earlier estimate on the size of $E(d)$, Corollaries \ref{cor:sizeofatoms} and \ref{cor:sizeofbilsets} and our choice of $\rho_0$, 
\[\min_{d \in \F_p^{3\ell+6q}\setminus \Gamma}|E(d)|\geq (1-\e)^4p^{-3q}|B_1||B_2||B_3|\geq (1-\e)^7p^{-3\ell-6q}(p^n)^3.\] 
Combining with (\ref{fgamma}) and rearranging, we have that $|\F_p^{3\ell+6q}\setminus \Gamma|\leq  \e (1+\e)^7p^{3\ell+6q}$, whence $|\Gamma|\geq (1-8\e)p^{3\ell+6q}$ as desired.

From now on fix $d=(a_1,a_2,a_3,b_{12},b_{13},b_{23})\in \Gamma$. We shall show that $\|f\|_{U^3(d)}<2\e^{1/16}$. For ease of notation let  $B_i=B(a_i)$ for each $i\in [3]$, let $\beta_{ij}=\beta(b_{ij})$ for each $1\leq i<j \leq 3$ and let 
\[\gamma=p^{-6n} \prod_{1\leq i<j \leq 3} |\beta_{ij}| \] 
be the product of the densities of the $\beta_{ij}$ in $\F_p^n\times \F_p^n$.

Given $(x_0,y_0,z_0)\in E(d)$, let 
\begin{align*}
O(d)[x_0,y_0,z_0]&=\{(x_1,y_1,z_1) \in B_1\times B_2\times B_3: \mbox{for each }\e\in \{0,1\}^3, (x_{\e(1)},y_{\e(2)},z_{\e(3)})\in E(d)\},
\end{align*}
and define
\[J(d)=\{(x,y,z)\in E(d): |O(d)[x,y,z]|\leq (1+\e)\gamma^3|B_1||B_2||B_3|\}.\]
By  our choice of $\rho_0$ and Lemma \ref{lem:extuples} (i), we have that $|E(d)\setminus J(d)|\leq \e p^{-6q} |E(d)|$. We also trivially have that for all $(x_0,y_0,z_0)\in E(d)$, $|O(d)[x_0,y_0,z_0]|\leq |E(d)|$.

Next observe that, by Cauchy-Schwarz,
\begin{align*}
\|f\|_{U^{3}(d)}^{16}=\big(\E_{\begin{subarray}{l} x_0,x_1\in B_1\end{subarray}}\E_{\begin{subarray}{l} y_0,y_1\in B_2\end{subarray}}\E_{\begin{subarray}{l} z_0,z_1\in B_3\end{subarray}}\\
\prod_{(i,j)\in \{0,1\}^2} \mu_{\beta_{12}}(x_i,y_j)&\prod_{(i,k)\in \{0,1\}^2} \mu_{\beta_{13}}(x_i,z_k)\prod_{(j,k)\in \{0,1\}^2} \mu_{\beta_{23}}(y_j,z_k)\\
&\prod_{\e\in \{0,1\}^3}\calC^{|\e|}f(x_{\e(1)}+y_{\e(2)}+z_{\e(3)})\big)^2
\end{align*}
is at most 
\begin{align*}
\big(\E_{\begin{subarray}{l} x_0,x_1\in B_1\end{subarray}}\E_{\begin{subarray}{l} y_0,y_1\in B_2\end{subarray}}\E_{\begin{subarray}{l} z_0,z_1\in B_3\end{subarray}}&\prod_{(i,j)\in \{0,1\}^2} \mu_{\beta_{12}}(x_i,y_j)\prod_{(i,k)\in \{0,1\}^2}\mu_{\beta_{13}}(x_i,z_k)\\
&\prod_{(j,k)\in \{0,1\}^2}\mu_{\beta_{23}}(y_j,z_k)\prod_{\e\in \{0,1\}^3:|\e|\; \textrm{odd}}|f(x_{\e(1)}+y_{\e(2)}+z_{\e(3)})|^2\big)\\
\cdot \big(\E_{\begin{subarray}{l} x_0,x_1\in B_1\end{subarray}}\E_{\begin{subarray}{l} y_0,y_1\in B_2\end{subarray}}\E_{\begin{subarray}{l} z_0,z_1\in B_3\end{subarray}}&\prod_{(i,j)\in \{0,1\}^2} \mu_{\beta_{12}}(x_i,y_j)\prod_{(i,k)\in \{0,1\}^2}\mu_{\beta_{13}}(x_i,z_k)\\
&\prod_{(j,k)\in \{0,1\}^2}\mu_{\beta_{23}}(y_j,z_k)\prod_{\e\in \{0,1\}^3:|\e|\; \textrm{even}}|f(x_{\e(1)}+y_{\e(2)}+z_{\e(3)})|^2\big),
\end{align*}
which by symmetry equals
\begin{align*}
\Big(\E_{\begin{subarray}{l} x_0,x_1\in B_1\end{subarray}}\E_{\begin{subarray}{l} y_0,y_1\in B_2\end{subarray}}\E_{\begin{subarray}{l} z_0,z_1\in B_3\end{subarray}}&\prod_{(i,j)\in \{0,1\}^2}\mu_{\beta_{12}}(x_i,y_j)\prod_{(i,k)\in \{0,1\}^2}\mu_{\beta_{13}}(x_i,z_k)\\
&\prod_{(j,k)\in \{0,1\}^2}\mu_{\beta_{23}}(y_j,z_k)\prod_{\e\in \{0,1\}^3:|\e|\; \textrm{odd}}|f(x_{\e(1)}+y_{\e(2)}+z_{\e(3)})|^2\Big)^2.
\end{align*}
We thus have that
\begin{align*}
\|f\|_{U^{3}(d)}^8\leq \E_{\begin{subarray}{l} x_0,x_1\in B_1\end{subarray}}\E_{\begin{subarray}{l} y_0,y_1\in B_2\end{subarray}}&\E_{\begin{subarray}{l} z_0,z_1\in B_3\end{subarray}}\prod_{(i,j)\in \{0,1\}^2} \mu_{\beta_{12}}(x_i,y_j)\prod_{(i,k)\in \{0,1\}^2}\mu_{\beta_{13}}(x_i,z_k)\\
&\prod_{(j,k)\in \{0,1\}^2}\mu_{\beta_{23}}(y_j,z_k)\prod_{\e\in \{0,1\}^3:|\e|\; \textrm{even}}|f(x_{\e(1)}+y_{\e(2)}+z_{\e(3)})|^2,
\end{align*}
which equals
\[\gamma^{-4} (|B_1||B_2||B_3|)^{-2} \sum_{(x_0,y_0,z_0)\in E(d)}\sum_{(x_1,y_1,z_1)\in O(d)[x_0,y_0,z_0]}\prod_{\e\in\{0,1\}^3:|\e|\; \textrm{even}}|f(x_{\e(1)}+y_{\e(2)}+z_{\e(3)})|^2.\]
We now split the sum over $(x_0,y_0,z_0)\in E(d)$ using the disjoint union
\[E(d)=J(d)\cup (E(d)\setminus J(d))=(J(d)\cap I(d)) \cup (J(d)\setminus I(d))\cup (E(d)\setminus J(d)).\]
Consider first the sum arising from $E(d)\setminus J(d)$, namely
\[\sum_{(x_0,y_0,z_0)\in E(d)\setminus J(d)}\sum_{(x_1,y_1,z_1)\in O(d)[x_0,y_0,z_0]}\prod_{\e\in\{0,1\}^3:|\e|\; \textrm{even}}|f(x_{\e(1)}+y_{\e(2)}+z_{\e(3)})|^2,\]
which, as $f$ is 1-bounded, is bounded above by
\[|E(d)\setminus J(d)|\max_{(x_0,y_0,z_0)\in E(d)}|O(d)[x_0,y_0,z_0]|\leq \e p^{-6q}|E(d)|^2 \leq \e(1+\e)^4\gamma^4|B_1|^2|B_2|^2|B_3|^2,\]
where the first inequality follows from our bounds on $|E(d)\setminus J(d)|$ and $|O(d)[x_0,y_0,z_0]|$ and the second from our upper bound on $|E(d)|$ and an application of Corollary \ref{cor:sizeofbilsets}.

Next consider the sum arising from $J(d)\setminus I(d)$, that is, 
\[\sum_{(x_0,y_0,z_0)\in J(d)\setminus I(d)}\sum_{(x_1,y_1,z_1)\in O(d)[x_0,y_0,z_0]}\prod_{\e\in\{0,1\}^3:|\e|\; \textrm{even}}|f(x_{\e(1)}+y_{\e(2)}+z_{\e(3)})|^2.\]
By definition of $J(d)$ and the fact that $f$ is 1-bounded, this sum is bounded above by
\[\sum_{(x_0,y_0,z_0)\in J(d)\setminus I(d)}|O(d)[x_0,y_0,z_0]|\leq |J(d)\setminus I(d)|(1+\e)\gamma^3|B_1||B_2||B_3|,\]
which in turn, using the fact that $J(d)\setminus I(d) \subseteq E(d)\setminus I(d)$, $d\in \Gamma$ and our upper bound on $|E(d)|$, is at most
\[\e^{1/2}|E(d)|(1+\e)\gamma^3|B_1||B_2||B_3|\leq \e^{1/2}(1+\e)^2\gamma^4|B_1|^2|B_2|^2|B_3|^2.\]
Finally, consider the third summand, i.e. the sum arising from $J(d)\cap I(d)$. Given $(x_0,y_0,z_0)\in J(d)\cap I(d)$, we have $|f(x_0+y_0+z_0)|^2\leq \e^{1/2}$ by definition of $I(d)$, so for every $(x_1,y_1,z_1)\in O(d)[x_0,y_0,z_0]$, 
$$
\prod_{\e\in\{0,1\}^3:|\e|\; \textrm{even}}|f(x_{\e(1)}+y_{\e(2)}+z_{\e(3)})|^2\leq |f(x_0+y_0+z_0)|^2< \e^{1/2}.
$$
Further, for any such $(x_0,y_0,z_0)$,  we know that $|O(d)[x_0,y_0,z_0]|\leq (1+\e)\gamma^3|B_1||B_2||B_3|$, by definition of $J(d)$.  Consequently, the third summand
\[\sum_{(x_0,y_0,z_0)\in J(d)\cap I(d)}\;\;\sum_{(x_1,y_1,z_1)\in O(d)[x_0,y_0,z_0]}\prod_{\e\in\{0,1\}^3:|\e|\; \textrm{even}}|f(x_{\e(1)}+y_{\e(2)}+z_{\e(3)})|^2\]
is bounded above by 
\[|E(d)|(1+\e)\gamma^3|B_1||B_2||B_3|\e^{1/2}\leq \e^{1/2}(1+\e)^2\gamma^4 |B_1|^2|B_2|^2|B_3|^2.\]
Combining the upper bounds for each summand yields 
\[\|f\|_{U^3(d)}^8\leq \gamma^{-4} \big( \e(1+\e)^4\gamma^4 + \e^{1/2}(1+\e)^2 \gamma^4 + \e^{1/2}(1+\e)^2 \gamma^4\big) \leq 2^8\e^{1/2},\]
as claimed.
\end{proof}

\vspace{2mm}

We are now in a position to prove the repackaging of the arithmetic regularity lemma that we shall use.\footnote{This first appeared as \cite[Proposition 3.20]{Terry.2021av2}.}

 \begin{proposition}\label{prop:unifatom2}
 For all $\e>0$, there is a polynomial growth function $\rho_1=\rho_1(\e)$ such that for all growth functions $\rho\geq\rho_1$, there is $n_1=n_1(\e,\rho)$ such that for all $\ell_0,q_0\geq 0$, there is $D=D(\e,\rho,\ell_0,q_0)$ such that the following holds.  Suppose $n\geq n_1$, $A\subseteq \F_p^n$, and $\calB_0=(\calL_0,\calQ_0)$ is a quadratic factor on $\F_p^n$ of complexity $(\ell_0,q_0)$. Then there are $\ell,q\leq D$, a quadratic factor $\calB=(\calL,\calQ)\preceq \calB_0$ of complexity $(\ell, q)$ together with a set $\Gamma\subseteq \mathbb{F}_p^{3\ell+6q}$ with the property that
\begin{enumerate}[label=\normalfont(\roman*)]
\item  $\calB$ has rank at least $\rho(\ell+q)$;
\item $|\Gamma|\geq (1-\e^{2}/2^{29})p^{3\ell+6q}$;
\item for all $d=(a_1,a_2,a_3,b_{12},b_{13},b_{23})\in \Gamma$, $\|1_A-\alpha_{B(\Sigma(d))}\|_{U^3(d)}<\e$.
\end{enumerate}
\end{proposition}
\begin{proof}
Fix $\e>0$, let $\e'=\e^2/2^{32}$, let $\rho_0=\rho_0(\e')$ be the growth function resulting from an application of Lemma \ref{lem:smallpart}, and let $\rho_2(\e)$ and $\rho_3(\e)$ be the growth functions resulting from applications of Corollary \ref{cor:sizeofatoms} and Corollary \ref{cor:sizeofbilsets}, respectively. Define  $\rho_1(x)=\max\{\rho_0(x),\rho_2(x),\rho_3(x)\}$ and $\sigma_1(x)=2^{26}p^{18x}/\e$. Observe that $\rho_1$ can be taken to be linear in $x$, with slope depending on $\e$. Suppose now that $\rho\geq \rho_1$ is a growth function.  Let $n_1=n_1(\e', \rho,\sigma_1)$ be from Theorem \ref{thm:gtarl}. 

Then by Theorem \ref{thm:gtarl}, given $\ell_0,q_0\geq 1$, $n\geq n_1$, $A\subseteq \F_p^n$, and any quadratic $\calB_0=(\calL_0,\calQ_0)$ in $\F_p^n$ of complexity $(\ell_0,q_0)$, there exists $D=D(\e', \rho, \sigma_1,\ell_0,q_0)$ and a quadratic factor $\calB=(\calL,\calQ)\preceq \calB_0$ of complexity $(\ell,q)$ with $\ell+q\leq D$ and rank at least $\rho(\ell+q)$, together with a decomposition $1_A=f_1+f_2+f_3$ with $f_1=\mathbb{E}(1_A|\calB)$, $||f_2||_{U^3}\leq 1/\sigma_1(\ell+q)$, $||f_3||_2<\e'$, and where $\|f_2\|_\infty, \|f_3\|_\infty\leq 1$.   Let $\Gamma\subseteq \F_p^{\ell+q}$ be as in Lemma \ref{lem:smallpart} applied to $\e'$ and $\rho_0$, $f_3$ and $\calB$, implying that $|\Gamma|\geq (1-8\e')p^{3\ell+6q}$. It suffices to show that for all $d=(a_1,a_2,a_3,b_{12},b_{23},b_{13})\in \Gamma$, $\|1_A-\alpha_{B(\Sigma(d))}\|_{U^3(d)}<\e$.

Fix $d=(a_1,a_2,a_3,b_{12},b_{23},b_{13})\in \Gamma$ and let $\alpha_{B(\Sigma(d))}=|A\cap B(\Sigma(d))|/|B(\Sigma(d))|$. For ease of notation, write $B_i=B(a_i)$ and $\beta_{ij}=\beta(b_{ij})$ for each $1\leq i<j\leq 3$, and let $g_A^{d}=1_A-\alpha_{B(\Sigma(d))}$. It is easy to see that on any $B\in \At(\calB)$, $1_A(x)-\alpha_{B}=f_2(x)+f_3(x)$, where $\alpha_B$ is the density of $A$ on $B$. In particular, since $g_A^{d}(x)=1_A-\alpha_{B(\Sigma(d))}$, we have that for all $x\in B_1$, $y\in B_2$ and $z\in B_3$ with $(x,y)\in \beta_{12}$, $(x,z)\in \beta_{13}$ and $(y,z)\in \beta_{23}$,
\[g_A^{d}(x+y+z)=f_2(x+y+z)+f_3(x+y+z).\]
Therefore, in order to estimate $\|1_A-\alpha_{B(\Sigma(d))}\|_{U^3(d)}^8=\|g_A^{d}\|_{U^3(d)}^8$, it suffices to consider\footnote{By slight abuse of notation, we have fixed any linear order $\omega_1,\dots,\omega_8$ on $\{0,1\}^3$ and defined $\langle f_{\omega_1},f_{\omega_2},f_{\omega_3},f_{\omega_4},f_{\omega_5},f_{\omega_6},f_{\omega_6},f_{\omega_8}\rangle_{U^3(d)}$ to be $\langle (f_{\omega})_{\omega\in \{0,1\}^3}\rangle$.}
\begin{align}\label{al:oct}
\|f_2+f_3\|_{U^3(d)}^8=\sum_{(i_1,\ldots, i_8)\in \{2,3\}^8}\langle f_{i_1},f_{i_2},f_{i_3},f_{i_4},f_{i_5},f_{i_6},f_{i_6},f_{i_8}\rangle_{U^3(d)},
\end{align}

We will give an upper bound for each term in the sum appearing in (\ref{al:oct}). We shall deal with the case when $i_j=3$ for all $j\in[8]$ first.

Indeed, by definition of $\Gamma$ as produced by Lemma \ref{lem:smallpart}, we have 
\begin{equation}\label{eq:firstcase}
\langle f_{3},f_{3},f_{3},f_{3},f_{3},f_{3},f_{3},f_{3}\rangle_{U^3(d)}= ||f_3||_{U^3(d)}^8<2^8\e'^{1/2}=\e/2^8.
\end{equation}
Now fix some $(i_1,\ldots, i_s)$ where $i_j=2$ for some $j\in[8]$. We claim that
\begin{equation}\label{eq:secondcase}
\langle f_{i_1},f_{i_2},f_{i_3},f_{i_4},f_{i_5},f_{i_6},f_{i_7},f_{i_8}\rangle_{U^3(d)}\leq \e/2^8.
\end{equation}
By relabeling if necessary, we may assume that $i_1=2$. Letting 
\[ \eta=p^{-3n}\prod_{i\in [3]}|B_i| \;\;\mbox{ and } \;\;\gamma=p^{-6n}\prod_{1\leq i<j \leq 3}|\beta_{ij}|\]
denote the products of the densities of the $B_i$ in $\F_p^n$ and the $\beta_{ij}$ in $\F_p^n\times \F_p^n$, respectively, observe that we may write the left-hand side of (\ref{eq:secondcase}) as 
\begin{equation}\label{eq:almostthere}
\eta^{-2}\gamma^{-4}\;\;\E_{x_1,y_1,z_1}\E_{x_0,y_0,z_0}f_2(x_0+y_0+z_0)h_{12}(x_0,y_0)h_{13}(x_0,z_0)h_{23}(y_0,z_0),
\end{equation}
where each $h_{ij}$ is a 1-bounded function in two variables depending on one or more of $x_1$, $y_1$ and $z_1$ (but at most two of $x_0$, $y_0$ and $z_0$).

It is well known that averages such as the expectation in (\ref{eq:almostthere}) over $x_0, y_0, z_0$ (for fixed $x_1,y_1,z_1$) are controlled by the $U^3$ norm of $f_2$. Indeed, applying for instance \cite[Theorem 2.3]{Gowers.2010tc}, we have
\[\big|\E_{x_0,y_0,z_0}f_2(x_0+y_0+z_0)h_{12}(x_0,y_0)h_{13}(x_0,z_0)h_{23}(y_0,z_0)\big|\leq \|f_2\|_{U^3},\]
and thus (\ref{eq:almostthere}) is bounded above by 
\begin{equation}\label{eq:f2bound}
\eta^{-2} \gamma^{-4}\;\;\|f_2\|_{U^3}\leq \eta^{-2} \gamma^{-4}\;\sigma_1(\ell+q)^{-1}.
\end{equation}
By our choice of $\rho_1$, we may use Corollary \ref{cor:sizeofatoms} to bound $\eta\geq p^{-3(\ell+q)}(1-\e)^3$ and Corollary \ref{cor:sizeofbilsets} to bound $\gamma\geq p^{-3q}(1-\e)^3$. 
It follows from our choice of $\sigma_1$ that (\ref{eq:f2bound}) is at most 
\[p^{6(\ell+q)}(1+\e)^6\cdot p^{12q}(1+\e)^{12} \cdot 2^{-26}\e p^{-18(\ell+q)} <\e/2^8,\]
as claimed.

Combining (\ref{al:oct}), (\ref{eq:firstcase}), and (\ref{eq:secondcase}) we have that 
\[\|1_A-\alpha_{B(\Sigma(d))}\|_{U^3(d)}^8=\|f_2+f_3\|_{U^3(d)}^8 < 2^8 \e/2^8=\e,\]
as desired.
\end{proof}

\subsection{Proof of Theorem \ref{cor:vc2}}\label{subsec:firstproof}
In this section we prove our main result about sets of bounded $\VC_2$-dimension.  First, we show that if a set $A$ has bounded $\VC_2$-dimension, then it has close to trivial density on any sufficiently uniform atom.\footnote{This result first appeared as \cite[Corollary 3.22]{Terry.2021av2}.}

\begin{proposition}\label{prop:trivdense}
For all $\e>0$ and $m\geq 1$, there exist $\eta=\eta(\e,m)>0$ and a polynomial growth function $\rho_0=\rho_0(\e,m)$ such that for all growth functions $\rho\geq \rho_0$, all $\ell,q\geq 0$ and all $n > \ell+q$ the following holds.

Let $A\subseteq \mathbb{F}_p^n$ be such that $\VC_2(A)<m$. Let $\calB=(\calL,\calQ)$ be a quadratic factor of complexity $(\ell,q)$ and rank at least $\rho(\ell+q)$, and let $d=(a_1,a_2,a_3,b_{12},b_{23},b_{13})\in \F_p^{\ell+q}\times \F_p^{\ell+q}\times \F_p^{\ell+q}\times\F_p^q\times \F_p^q\times \F_p^q$. 

Suppose that $\|1_A-\alpha_{B(\Sigma(d))}\|_{U^3(d)}<\eta$. Then the density $\alpha_{B(\Sigma(d))}$ of $A$ on the atom $B(\Sigma(d))$ satisfies $\alpha_{B(\Sigma(d))} \in [0,\e)\cup (1-\e,1]$.
\end{proposition}

\begin{proof} Fix $\e>0$ and $m\geq 1$, let $\eta = (\e/4)^{m^2 2^{m^2}}$, and let the growth function $\rho_0$ be defined by $\rho_0(x)=128(4m2^{m^2}x+\log_p(\eta^{-1}))$. Let $A\subseteq \mathbb{F}_p^n$ be such that $\VC_2(A)<m$, and let $\calB=(\calL,\calQ)$ be a quadratic factor of complexity $(\ell,q)$ and rank at least $\rho(\ell+q)$, where $\rho\geq \rho_0$.

Let $g_A^{d}(x)=1_A(x)-\alpha_{B(\Sigma(d))}$ and $g_{A^C}^d=1_{A^C}(x)-\overline{\alpha}_{B(\Sigma(d))}$, where $\overline{\alpha}_{B(\Sigma(d))}$ is the density of $A^C$ on $B(\Sigma(d))$. By linearity, 
\[T_{m-\IP_2(d)}(1_A |1_{A^C})=T_{m-\IP_2(d)}(g_{A}^{d}+\alpha_{B(\Sigma(d))} |g_{A^C}^{d}+\overline{\alpha}_{B(\Sigma(d))})\]
equals
\[\alpha_{B(\Sigma(d))}^{m^2 2^{m^2-1}} (1-\alpha_{B(\Sigma(d))})^{m^2 2^{m^2-1}}\]
plus $2^{m^2 2^{m^2}} -1$ terms of the form
\begin{equation}\label{eq:error}
T_{m-\IP_2(d)}((f_{i,j,S})_{i,j\in [m],S\subseteq [m]^2}),
\end{equation}
in which at least one of the input functions $f_{i,j,S}$ equals $g_A^{d}$ or $g_{A^C}^{d}=-g_A^{d}$. Thus, by Proposition \ref{prop:ip2loccontrol}, if $\|g_A^{d}\|_{U^3(d)}<\eta$, then each expression of the form (\ref{eq:error}) is bounded above by 
\[(1+O_m(p^{4m2^{m^2}(\ell+q)-\rho(\ell+q)/2})) \eta +O_m(p^{2m2^{m^2}(\ell+q)-\rho(\ell+q)/128})<3\eta,\]
where the inequality is due to our choice of $\rho_0$. It follows that
\begin{equation}\label{eq:lowerbound}
T_{m-\IP_2(d)}(g_{A}^{d}+\alpha_{B(\Sigma(d))} |g_{A^C}^{d}+\overline{\alpha}_{B(\Sigma(d))})\geq \alpha_{B(\Sigma(d))}^{m^2 2^{m^2-1}} (1-\alpha_{B(\Sigma(d))})^{m^2 2^{m^2-1}}-3(2^{m^2 2^{m^2}} -1)\eta.
\end{equation}
Now suppose towards a contradiction that $\alpha_{B(\Sigma(d))}\in [\e,1-\e]$. Then the first term in (\ref{eq:lowerbound}) is bounded below by $\e^{m^2 2^{m^2}}$,
and since $\eta$ was chosen sufficiently small in terms of $\e$ and $m$, we have that $T_{m-\IP_2(d)}(1_{A}|1_{A^C})>0$. This means that $A$ contains an instance of $m$-$\IP_2$, contradicting the assumption that $\VC_2(A)<m$. 
\end{proof}

We now have all the tools to prove Theorem \ref{thm:vc2}, that is, that a set of bounded $\VC_2$-dimension has near-trivial density on almost all atoms of a high-rank, bounded-complexity quadratic factor.

\vspace{2mm}

\begin{proofof}{Theorem \ref{thm:vc2}}
Fix $m\geq 1$ and $\mu>0$.   Apply Proposition \ref{prop:trivdense} to obtain $\eta=\eta(\mu,m)>0$ and a growth function $\rho_0(\mu,m)$. Let $\e=2^{14}\mu^{1/2}\eta$, and let $\rho_1=\rho_1(\e)$ be the growth function given by  Proposition \ref{prop:unifatom2}. Define $\sigma_1=\max\{\rho_0, \rho_1\}$. Let $\sigma\geq \sigma_1$ be a growth function, and let $n_1=n_1(\e,\sigma)$ and $D=D(\e, \sigma,0,0)$ be given by Proposition \ref{prop:unifatom2}. Note that our choices of $n_1$ and $D$ depended only on $m$, $\mu$, and $\sigma$.

Suppose that $n\geq n_1$ and that $A\subseteq \mathbb{F}_p^n$ has $\VC_2(A)<m$.  By Proposition \ref{prop:unifatom2}, there exists a quadratic factor $\calB=(\calL,\calQ)$ of complexity $(\ell,q)$ and rank at least $\sigma(\ell+q)$, together with a set $\Gamma\subseteq \mathbb{F}_p^{3\ell+6q}$ such that
\begin{enumerate}
\item $\ell,q\leq D$;
\item $|\Gamma|\geq (1-\e^2/2^{29})p^{3\ell+6q}$;
\item for all $d=(a_1,a_2,a_3,b_{12},b_{13},b_{23})\in \Gamma$, 
$$
\|1_A-\alpha_{B(\Sigma(d))}\|_{U^3(d)}<  \e\leq\eta,
$$
\end{enumerate}
where the last inequality holds by our choice of $\e$.  Let 
$$
X=\{ B(b)\in \At(\calB):\text{ there exists $d\in \Gamma$ with $\Sigma(d)=b$}\}.
$$
By Proposition \ref{prop:trivdense}, every $B\in X$ satisfies $|A\cap B|/|B|\in [0,\mu)\cup (1-\mu,1]$.  It is straightforward to show that $|\Gamma|\geq (1-\e^2/2^{29})p^{3\ell+6q}$ implies that 
$$
|X|\geq (1-\e^2/2^{29})p^{\ell+q}\geq (1-\mu)|\At(\calB)|,
$$
where the second inequality is by choice of $\e$.  This finishes the proof. 
\end{proofof}

\vspace{2mm}

We end this section with a proof of Corollary \ref{cor:vc2}, which asserts that sets of bounded $\VC_2$-dimension look approximately like unions of quadratic atoms. 

\vspace{2mm}

\begin{proofof}{Corollary \ref{cor:vc2}}Fix $\mu>0$ and $m\geq 1$. Let $\sigma_1=\sigma_1(\mu^2, m)$ be from Theorem \ref{thm:vc2}, and let $\rho=\rho(\mu^2)$ be as in Corollary \ref{cor:sizeofatoms}.  Let $\sigma_2=\max\{\sigma_1,\rho\}$, and for any $\sigma\geq \sigma_2$, let $n_1=n_1(m,\mu^2,\sigma)$ and $D=D(m,\mu^2,\sigma)$ be from Theorem \ref{thm:vc2}. 

Suppose that $n\geq n_1$ and that $A\subseteq \F_p^n$ is such that $\VC_2(A)<m$.  By Theorem \ref{thm:vc2}, there are $\ell,q\leq D$, and a quadratic factor $\calB$ of complexity at $(\ell,q)$ and rank at least $\sigma(\ell+q)$ together with a set $X\subseteq \At(\calB)$ satisfying $|X|\geq (1-\mu^2)p^{\ell+q}$, such that for all $B\in X$, $|A\cap B|/|B|\in [0,\mu^2)\cup (1-\mu^2,1]$. 

Let $I_0=\{B\in X :|A\cap B|/|B|\in [0,\mu^2)\}$, $I_1=\{B\in X :|A\cap B|/|B|\in (1-\mu^2,1]\}$, and set $I_2=\At(\calB)\setminus X$.  Let $Y=\bigcup_{B\in I_1}B$.  We show that $|A\Delta Y|\leq \mu |\mathbb{F}_p^n|$, which will finish the proof.  Observe that 
\[A\Delta Y=\Big( \bigcup_{B\in I_0\cup I_2}(A\cap B)\Big)\cup  \Big( \bigcup_{B\in I_1}(B\setminus A)\Big),\]
where all unions are disjoint.
Using this along with our assumptions on $I_0,I_1,I_2$, we obtain that
\[|A\Delta Y|=\sum_{B\in I_0\cup I_2}|A\cap B|+\sum_{B\in I_1}|B\setminus A|\leq \sum_{B\in I_2}|B|+\sum_{B\in I_0} \mu^2|B|+ \sum_{B\in I_1}\mu^2|B|.\]
By Corollary \ref{cor:sizeofatoms}, we therefore have that
\[|A\Delta Y|\leq (1+\mu^2)p^{n-\ell-q}\big(|I_2|+\mu^2(|I_0|+|I_1|)\big)\leq (1+\mu^2)p^{n-\ell-q}\cdot2\mu^2 p^{\ell+q}\leq \mu p^n,\]
as claimed.
\end{proofof}

\appendix

\section{Bilinear exponential sum estimates}\label{app:expsums}

To begin with, we need to establish a number of standard facts. These will (hopefully) look amply plausible to the expert, having appeared in various guises in the context of hypergraph counting lemmas and arithmetic counting relative to pseudorandom measures. However, we were unable to locate them in the literature in the exact form needed here.\footnote{Versions of these counting lemmas previously appeared in \cite[Appendix C]{Terry.2021av2}.}

Let $\calB=(\calL,\calQ)$ be a quadratic factor on $\F_p^n$ of complexity $(\ell,q)$, with $\calL=\{r_1,\dots, r_\ell\}$ and $\calQ=\{M_1,\dots,M_q\}$. Recall that for $a=(a_1,\dots,a_{\ell+q})\in \F_p^{\ell+q}$ and $b=(b_1,\ldots, b_q)\in \F_p^q$,  
$$
B(a)=\{x\in \F_p^n: \text{ for each }i\in [\ell], x^Tr=a_i\text{ and for each }j\in [q], x^TM_jx=b_j\}
$$
and
$$
\beta(b)=\{(x,y)\in \F_p^n\times \F_p^n: \text{ for each }j\in [q], x^TM_jy=b_j\}.
$$
We begin with a simple lemma closely related to Fact \ref{fct:expsum}.

\begin{lemma}[Bilinear exponential sums of high rank are small]\label{lem:bilsums} 
Let $\calB=(\calL,\calQ)$ be a quadratic factor on $\F_p^n$ of complexity $(\ell,q)$ and rank at least $\tau$ with $\calL=\{r_1,\dots, r_\ell\}$ and $\calQ=\{M_1,\dots,M_q\}$.
\begin{enumerate}[label=\normalfont(\roman*)]
\item Given any non-trivial linear combination $M$ of $M_1,\dots, M_q$, and any $c,d\in\F_p^n$, we have
\[\big|\E_{x,y\in \F_p^n} \omega^{x^TMy+c^Tx+d^Ty}\big|=O( p^{-\tau}).\]
\item Given $a=(a_1,\dots,a_{\ell+q})\in \F_p^{\ell+q}$ and any non-trivial linear combination $M$ of $M_1,\dots, M_q$, the probability that $x\in B(a)$ lies in $\ker(M)$ is $O(p^{\ell+q-\tau/2})$.
\item Given $b=(b_1,\dots,b_q)\in\F_p^q$, 
\[\E_{x,y\in \F_p^n}1_{\beta(b)}(x,y)=p^{-q}(1+O(p^{q-\tau})).\]
\end{enumerate}
\end{lemma}

\vspace{2mm}

\begin{proofof}{Lemma \ref{lem:bilsums}} For part (i), note that our assumption on the rank of the factor implies $M$ has rank at least $\tau$.  Observe further that
\[\big|\E_{x,y\in \F_p^n} \omega^{x^TMy+c^Tx+d^Ty}\big|=\big|\E_{x\in \F_p^n}\omega^{c^Tx}\E_{y\in \F_p^n} \omega^{(x^TM+d^T)y}\big|.\]
For a fixed $x$, (\ref{star}) implies that exponential sum  over $y$ is $0$ unless $x^TM=-d^T$, in which case it equals 1. Since $M$ has rank at least $\tau$, the proportion of $x\in \F_p^n$ satisfying $x^TM=-d^T$ is at most $p^{-\tau}$. It follows that 
\[\big|\E_{x,y\in \F_p^n} \omega^{x^TMy+c^Tx+d^Ty}\big|\leq \E_{x\in \F_p^n}\big|\E_{y\in \F_p^n} \omega^{(x^TM+d^T)y}\big|=O(p^{-\tau}).\]
This completes the proof of the first part.

We begin the proof of (ii) by observing that for $a=(a_1,\dots,a_{\ell+q})=(a_1,\dots,a_\ell,a_1',\dots,a_q')\in \F_p^{\ell+q}$, we may use (\ref{starstar}) to write the indicator function of the quadratic atom $B(a)$ as
\begin{align}\label{ba}
1_{B(a)}(x)=\E_{\begin{subarray}{l}  v_s\in \F_p\\s\in [\ell]\end{subarray}}\E_{\begin{subarray}{l}  u_t\in \F_p\\t\in [q]\end{subarray}}\omega^{x^T\sum_{t\in [q]}u_tM_t x-\sum_{t\in [q]}u_ta_t'+x^T\sum_{s\in [\ell]}v_sr_s-\sum_{s\in [\ell]}v_sa_s}.
\end{align}
Now using (\ref{ba}), the assumption on the rank of the factor, Lemma \ref{lem:sizeofatoms}, and (\ref{star}), the probability that an element of $B(a)$ lies in $\ker(M)$  is
\begin{align*}
\E_{x\in B(a)}\E_{z\in \F_p^n}\omega^{z^TMx}=p^{\ell+q}(1+O(p^{\ell+q-\tau/2}))\E_{\begin{subarray}{l}  v_s\in \F_p\\s\in [\ell]\end{subarray}}&\E_{\begin{subarray}{l}  u_t\in \F_p\\t\in [q]\end{subarray}}\omega^{-\sum_{t\in [q]}u_ta_t'-\sum_{s\in [\ell]}v_sa_s}\\
&\E_{z,x\in \F_p^n} \omega^{x^T\sum_{t\in [q]}u_tM_t x+x^T\sum_{s\in [\ell]}v_sr_s+x^TMz}.
\end{align*}
By the assumption on the rank of the factor, and Fact \ref{fct:expsum} applied to the exponential sum in $x$, this expression reduces to
\[p^{\ell}(1+O(p^{\ell+q-\tau/2}))\E_{\begin{subarray}{l}  v_s\in \F_p\\s\in [\ell]\end{subarray}}\omega^{-\sum_{s\in [\ell]}v_sa_s}
\E_{x\in \F_p^n} \omega^{x^T\sum_{s\in [\ell]}v_sr_s}\E_{z\in \F_p^n}\omega^{x^TMz} +O(p^{\ell+q-\tau/2}).\]
By the assumption on the rank of the factor, the exponential sum in $z$ is 0 for all but an $O(p^{-\tau})$-proportion of $x\in \F_p^n$, whence 
\[\E_{x\in B(a)}\E_{z\in \F_p^n}\omega^{z^TMx}=p^{\ell-\tau}(1+O(p^{\ell+q-\tau/2}))+ O(p^{\ell+q-\tau/2})=O(p^{\ell+q-\tau/2}).\]
This completes the proof of the second part.

To see (iii), similarly to (ii) note that for $b=(b_1,\dots,b_q)\in\F_p^q$, we can use (\ref{starstar}) to write  the indicator function of the bilinear variety $\beta(b)$ as 
\[1_{\beta(b)}(x,y)=\E_{\begin{subarray}{l}  w_t\in \F_p\\t\in [q]\end{subarray}}\omega^{x^T\sum_{t\in [q]}w_tM_t y-\sum_{t\in [q]}w_tb_t},\]
so 
\[\E_{x,y\in \F_p^n}1_{\beta(b)}(x,y)=\E_{\begin{subarray}{l}  w_t\in \F_p\\t\in [q]\end{subarray}}\omega^{-\sum_{t\in [q]}w_tb_t}\E_{x,y\in \F_p^n}\omega^{x^T\sum_{t\in [q]}w_tM_t y}.\]
By the assumption on the rank of the factor and part (i) with $c=d=0$, the bilinear exponential sum in $x$ and $y$ is $O(p^{-\tau})$ unless $w_t=0$ for all $t\in [q]$, in which case it equals 1. The claim follows.
\end{proofof}

\vspace{3mm}

We use Lemma \ref{lem:bilsums} to prove the following more general lemma about bilinear exponential sums. Its proof consists of further standard (albeit laborious) calculations. The underlying intuition is that if a quadratic factor has high-rank, then the associated bilinear level sets are sufficiently quasirandom for all reasonable counting problems to approximately return the expected frequency.

Recall that for a set $\beta\subseteq \F_p^n\times\F_p^n$, we write $\mu_{\beta}$ for the function on $\F_p^n\times \F_p^n$ defined by $(x,y)\mapsto 1_{\beta}(x,y)  |\F_p^n|^2/|\beta|$.

\begin{lemma}\label{lem:genbilsums}
Let $m\geq 1$, let $\calB=(\calL,\calQ)$ be a quadratic factor on $\F_p^n$ of complexity $(\ell,q)$ and rank at least $\tau$, where $\calL=\{r_1,\dots, r_\ell\}$ and $\calQ=\{M_1,\dots,M_q\}$. Let $d=(a_1, a_2, a_3, b_{12}, b_{13}, b_{23})\in \F_p^{\ell+q}\times \F_p^{\ell+q}\times \F_p^{\ell+q}\times\F_p^q\times \F_p^q\times \F_p^q$. Then for all sets $I, J, K \subseteq [m]$,
\begin{align*}
\E_{\begin{subarray}{l}  x_i\in B(a_1)\\i\in I\end{subarray}}\E_{\begin{subarray}{l}y_j\in B(a_2)\\j\in J\end{subarray}}\E_{\begin{subarray}{l} z_k\in B(a_3)\\k\in K \end{subarray}}\prod_{i\in I, j\in J} &\mu_{\beta(b_{12})}(x_i,y_j)\prod_{i\in I, k\in K}\mu_{\beta(b_{13})}(x_i,z_k)\prod_{ j\in J, k\in K}\mu_{\beta(b_{23})}(y_j,z_k)\\
&=1+O_m(p^{(\ell+q)(|I|+|J|+|K|)+q(|I||J|+|I||K|+|J||K|)-\tau/2}).
\end{align*}
The same holds whenever any instance of $\mu_{\beta(b_{ij})}$ is replaced by the constant function 1.
\end{lemma}

\begin{proof}As in the proof of Lemma \ref{lem:bilsums}, note that for $b=(b_1,\dots,b_q)\in\F_p^q$, we may use (\ref{starstar}) to write the indicator function of the bilinear variety $\beta(b)$  as 
\[1_{\beta(b)}(x,y)=\E_{\begin{subarray}{l}  w_t\in \F_p\\t\in [q]\end{subarray}}\omega^{x^T\sum_{t\in [q]}w_tM_t y-\sum_{t\in [q]}w_tb_t}.\]
Similarly, for $a=(a_1,\dots,a_{\ell+q})=(a_1,\ldots, a_{\ell},a_1',\ldots, a_q')\in \F_p^{\ell+q}$, the indicator function of the quadratic atom $B(a)$ can be written as
\[1_{B(a)}(x)=\E_{\begin{subarray}{l}  v_s\in \F_p\\s\in [\ell]\end{subarray}}\E_{\begin{subarray}{l}  u_t\in \F_p\\t\in [q]\end{subarray}}\omega^{x^T\sum_{t\in [q]}u_tM_t x-\sum_{t\in [q]}u_ta_t'+x^T\sum_{s\in [\ell]}v_sr_s-\sum_{s\in [\ell]}v_sa_s}.\]
We are thus able to rewrite 
\begin{align}\label{eq:notnorm}
\E_{\begin{subarray}{l}  x_i\in \F_p^n\\i\in I\end{subarray}}\E_{\begin{subarray}{l}y_j\in \F_p^n\\j\in J\end{subarray}}\E_{\begin{subarray}{l} z_k\in \F_p^n\\k\in K \end{subarray}}\prod_{i\in I} &1_{B(a_1)}(x_i)\prod_{j\in J}1_{B(a_2)}(y_j)\prod_{k\in K}1_{B(a_3)}(z_k)\\\nonumber
&\prod_{i\in I, j\in J}1_{\beta(b_{12})}(x_i,y_j)\prod_{i\in I, k\in K}1_{\beta(b_{13})}(x_i,z_k)\prod_{j\in J, k\in K}1_{\beta(b_{23})}(y_j,z_k)
\end{align}
as 
\begin{align}\label{eq:massive}
\E_{\begin{subarray}{l}  v_s^i\in \F_p\\s\in [\ell], i\in I\end{subarray}}\E_{\begin{subarray}{l}  u_t^i\in \F_p\\t\in [q],i\in I\end{subarray}}&
\E_{\begin{subarray}{l}  v_s^j\in \F_p\\s\in [\ell],j\in J\end{subarray}}\E_{\begin{subarray}{l}  u_t^j\in \F_p\\t\in [q],j\in J\end{subarray}}
\E_{\begin{subarray}{l}  v_s^k\in \F_p\\s\in [\ell],k\in K\end{subarray}}\E_{\begin{subarray}{l}  u_t^k\in \F_p\\t\in [q],k\in K\end{subarray}}
\E_{\begin{subarray}{l}  w_t^{ij}\in \F_p\\t\in [q], (i,j)\in I\times J\end{subarray}}
\E_{\begin{subarray}{l}  w_t^{ik}\in \F_p\\t\in [q], (i,k)\in I\times K\end{subarray}}
\E_{\begin{subarray}{l}  w_t^{jk}\in \F_p\\t\in [q], (j,k)\in J\times K\end{subarray}}\\\nonumber
\E_{\begin{subarray}{l}  x_i\in \F_p^n\\i\in I\end{subarray}}&\E_{\begin{subarray}{l}y_j\in \F_p^n\\j\in J\end{subarray}}\E_{\begin{subarray}{l} z_k\in \F_p^n\\k\in K \end{subarray}}
\;\omega^{\sum_{i\in I}x_i^T\sum_{t\in [q]}u_t^iM_t x_i -\sum_{i\in I}\sum_{t\in [q]}u_t^i a_t'+\sum_{i\in I}x_i^T\sum_{s\in [\ell]}v_s^ir_s-\sum_{i\in I}\sum_{s\in [\ell]}v_s^ia_s}\\\nonumber
&\hspace{10pt}\omega^{\sum_{j\in J}y_j^T\sum_{t\in [q]}u_t^jM_t y_j -\sum_{j\in J}\sum_{t\in [q]}u_t^j a_t'+\sum_{j\in J}y_j^T\sum_{s\in [\ell]}v_s^jr_s-\sum_{j\in J}\sum_{s\in [\ell]}v_s^ja_s}\\\nonumber
&\hspace{20pt}\omega^{\sum_{k\in K}z_k^T\sum_{t\in [q]}u_t^kM_t z_k -\sum_{k\in K}\sum_{t\in [q]}u_t^k a_t'+\sum_{k\in K}z_k^T\sum_{s\in [\ell]}v_s^kr_s-\sum_{k\in K}\sum_{s\in [\ell]}v_s^ka_s}\\\nonumber
&\hspace{30pt}\omega^{\sum_{i\in I, j\in J}x_i^T\sum_{t\in [q]}w_t^{ij}M_t y_j-\sum_{i\in I, j\in J}\sum_{t\in [q]}w_t^{ij}b_t + \sum_{i\in I, k\in K}x_i^T\sum_{t\in [q]}w_t^{ik}M_t z_k}\\\nonumber
&\hspace{40pt}\omega^{-\sum_{i\in I, k\in K}\sum_{t\in [q]}w_t^{ik}b_t\sum_{j\in J, k\in K}y_j^T\sum_{t\in [q]}w_t^{jk}M_t z_k-\sum_{j\in J, k\in K}\sum_{t\in [q]}w_t^{jk}b_t}.
\end{align}
In the above expression we have not subscripted the labels $a$ and $b$ of the quadratic atoms and bilinear varieties, respectively, since these turn out to be irrelevant in the course of the argument, as the reader will readily verify. 

Fortunately from this point onwards things simplify quickly: first, note that by Fact \ref{fct:expsum} and the assumption that $\calB$ has rank at least $\tau$, for any fixed $i\in I$, and any fixed choices for the $y_j$'s and $z_k$'s,
\begin{equation}\label{eq:singsum} 
\E_{x_i\in \F_p^n}\omega^{x_i^T\sum_{t\in [q]}u_t^iM_t x_i+ x_i^T\sum_{s\in [\ell]}v_s^ir_s+ x_i^T\sum_{t\in [q],j\in J}w_t^{ij}M_t y_j+x_i^T\sum_{t\in [q],k\in K}w_t^{ik}M_t z_k}=O(p^{-\tau/2})
\end{equation}
(here no summation convention has been employed) unless $u_t^i=0$ for all $t\in [q]$. Since the analogous facts about $y_j$, $z_k$ hold for $j\in J$, $k\in K$, respectively, we have that up to an additive error of $O_m(p^{-\tau/2})$,(\ref{eq:massive}) equals 
\begin{align}\label{eq:lessmassive}
p^{-q(|I|+|J|+|K|)}&\E_{\begin{subarray}{l}  v_s^i\in \F_p\\s\in [\ell], i\in I\end{subarray}}
\E_{\begin{subarray}{l}  v_s^j\in \F_p\\s\in [\ell],j\in J\end{subarray}}
\E_{\begin{subarray}{l}  v_s^k\in \F_p\\s\in [\ell],k\in K\end{subarray}}
\omega^{-\sum_{i\in I}\sum_{s\in [\ell]}v_s^ia_s-\sum_{j\in J}\sum_{s\in [\ell]}v_s^ja_s-\sum_{k\in K}\sum_{s\in [\ell]}v_s^ka_s} \\\nonumber
\E_{\begin{subarray}{l}  w_t^{ij}\in \F_p\\t\in [q], (i,j)\in I\times J\end{subarray}}&
\E_{\begin{subarray}{l}  w_t^{ik}\in \F_p\\t\in [q], (i,k)\in I\times K\end{subarray}}
\E_{\begin{subarray}{l}  w_t^{jk}\in \F_p\\t\in [q], (j,k)\in J\times K\end{subarray}}\\\nonumber
\hspace{10pt}&\omega^{-\sum_{i\in I, j\in J}\sum_{t\in [q]}w_t^{ij}b_t-\sum_{i\in I, k\in K}\sum_{t\in [q]}w_t^{ik}b_t-\sum_{j\in J, k\in K}\sum_{t\in [q]}w_t^{jk}b_t}\\\nonumber
&\hspace{10pt}\E_{\begin{subarray}{l}  x_i\in \F_p^n\\i\in I\end{subarray}}\E_{\begin{subarray}{l}y_j\in \F_p^n\\j\in J\end{subarray}}\E_{\begin{subarray}{l} z_k\in \F_p^n\\k\in K \end{subarray}}
\;\omega^{\sum_{i\in I}x_i^T\sum_{s\in [\ell]}v_s^ir_s+\sum_{j\in J} y_j^T\sum_{s\in [\ell]}v_s^jr_s+\sum_{k\in K}z_k^T\sum_{s\in [\ell]}v_s^kr_s}\\\nonumber
&\hspace{30pt}\omega^{\sum_{i\in I, j\in J}x_i^T\sum_{t\in [q]}w_t^{ij}M_t y_j + \sum_{i\in I, k\in K}x_i^T\sum_{t\in [q]}w_t^{ik}M_t z_k+\sum_{j\in J, k\in K}y_j^T\sum_{t\in [q]}w_t^{jk}M_t z_k}.
\end{align}
Now fix $(x_i)_{i\in I}$ and $(y_j)_{j\in J}$, and consider for any $k\in K$ the sum 
\begin{equation}\label{eq:linsum}
\E_{z_k \in \F_p^n}\omega^{z_k^T(\sum_{s\in [\ell]}v_s^k r_s+ \sum_{t\in [q],i\in I}w_t^{ik}M_t x_i+\sum_{t\in [q],j\in J}w_t^{jk}M_t y_j)}.\end{equation}
By (\ref{star}), this exponential sum will be zero unless 
\begin{align}\label{sr}
\sum_{s\in [\ell]}v_s^k r_s+ \sum_{t\in [q],i\in I}w_t^{ik}M_t x_i+\sum_{t\in [q],j\in J}w_t^{jk}M_t y_j= 0.
\end{align}
 But for fixed $j_0\in J$, $(x_i)_{i\in I}$ and $(y_j)_{j\in J\setminus\{j_0\}}$, if $w_t^{j_0k}\neq 0$ for some $t\in [q]$, then the proportion of $y_{j_0}\in \F_p^n$ for which (\ref{sr}) holds is $O(p^{-\tau})$. Since this holds for any $j_0\in J$, and symmetrically for any $i_0\in I$, it follows that up to an additive error of $O_m(p^{-\tau/2})$,  (\ref{eq:massive}) equals
\begin{align}\label{eq:evenlessmassive}
&p^{-q(|I|+|J|+|K|)-q(|I||J|+|I||K|+|J||K|)}\\\nonumber
&\hspace{40pt}\E_{\begin{subarray}{l}  v_s^i\in \F_p\\s\in [\ell], i\in I\end{subarray}}
\E_{\begin{subarray}{l}  v_s^j\in \F_p\\s\in [\ell],j\in J\end{subarray}}
\E_{\begin{subarray}{l}  v_s^k\in \F_p\\s\in [\ell],k\in K\end{subarray}}
\omega^{-\sum_{i\in I}\sum_{s\in [\ell]}v_s^ia_s-\sum_{j\in J}\sum_{s\in [\ell]}v_s^ja_s-\sum_{k\in K}\sum_{s\in [\ell]}v_s^ka_s}\\\nonumber
&\hspace{60pt}\E_{\begin{subarray}{l}  x_i\in \F_p^n\\i\in I\end{subarray}}\E_{\begin{subarray}{l}y_j\in \F_p^n\\j\in J\end{subarray}}\E_{\begin{subarray}{l} z_k\in \F_p^n\\k\in K \end{subarray}}
\omega^{\sum_{i\in I}x_i^T\sum_{s\in [\ell]}v_s^ir_s+\sum_{j\in J} y_j^T\sum_{s\in [\ell]}v_s^jr_s+\sum_{k\in K}z_k^T\sum_{s\in [\ell]}v_s^kr_s}.
\end{align}
The latter is now a straightforward linear exponential sum, which equals $0$ unless $v_s^i$, $v_s^j$, $v_s^k$ are $0$ for all $s\in [\ell]$ and $i\in I$, $j\in J$, $k\in K$, respectively. It follows that (\ref{eq:notnorm}) equals $p^{-(\ell+q)(|I|+|J|+|K|)-q(|I||J|+|I||K|+|J||K|)}+O_m(p^{-\tau/2})$.

It remains to observe that by Lemma \ref{lem:sizeofatoms}, for every $i\in [3]$, 
$$
|B(a_i)|=p^{-(\ell+q)}(1+O(p^{\ell+q-\tau/2})p^n,
$$
 and by Lemma \ref{lem:bilsums} (iii), for all $i<j\in [3]$, $|\beta(b_{ij})|=p^{-q}(1+O(p^{q-\tau}))p^{2n}$. Therefore the sought-after average
\begin{equation}\label{eq:aim}
\E_{\begin{subarray}{l}  x_i\in B(a_1)\\i\in I\end{subarray}}\E_{\begin{subarray}{l}y_j\in B(a_2)\\j\in J\end{subarray}}\E_{\begin{subarray}{l} z_k\in B(a_3)\\k\in K \end{subarray}}\prod_{i\in I, j\in J} \mu_{\beta(b_{12})}(x_i,y_j)\prod_{i\in I,k\in K}\mu_{\beta(b_{13})}(x_i,z_k)\prod_{ j\in J, k\in K}\mu_{\beta(b_{23})}(y_j,z_k)
\end{equation}
equals 
\begin{align*} 
\big(p^{-(\ell+q)(|I|+|J|+|K|)-q(|I||J|+|I||K|+|J||K|)}&+O_m(p^{-\tau/2})\big)\\
\big(p^{\ell+q}(1+O(p^{\ell+q-\tau/2})& \big)^{|I|+|J|+K|}\big(p^{q}(1+O(p^{q-\tau}) \big)^{|I||J|+|I||K|+|J||K|}+O_m(p^{-\tau/2}),
\end{align*}
where the second factor in the product addresses the normalisation inherent in the expectation notation and the third factor that in the definition of the $\mu_{\beta(b_{ij})}$. A simple computation yields the bound on (\ref{eq:aim}) asserted in the statement of the lemma.

It is not difficult to verify that when any instance of $\mu_{\beta(b_{ij})}$ is replaced by the constant function 1, the same argument goes through with one fewer constraint, yielding the same bound with a slightly improved error term (which will irrelevant for us).
\end{proof}

It turns out that we will also need to be able to count configurations of the type appearing in Lemma \ref{lem:genbilsums} whilst holding some variables constant.

\begin{lemma}\label{lem:extuples}
Let $m\geq 1$ and let $\calB=(\calL,\calQ)$ be a quadratic factor on $\F_p^n$ of complexity $(\ell,q)$ and rank at least $\tau$ with $\calL=\{r_1,\dots, r_\ell\}$ and $\calQ=\{M_1,\dots,M_q\}$. Let $d=(a_1, a_2, a_3, b_{12}, b_{13}, b_{23})\in \F_p^{\ell+q}\times \F_p^{\ell+q}\times \F_p^{\ell+q}\times\F_p^q\times \F_p^q\times \F_p^q$. Then 
\begin{enumerate}[label=\normalfont(\roman*)]
\item for all but an $O(p^{5\ell+18q-\tau/4})$-proportion of $(x',y',z')$ with $x' \in B(a_1)$, $y' \in B(a_2)$, $z' \in B(a_3)$ and $(x',y')\in \beta(b_{12})$, $(x',z')\in \beta(b_{13})$, $(y',z')\in \beta(b_{23})$, the function $g$ defined by
\begin{align*}
g(x',y',z')= \E_{x\in B(a_1)}\E_{y\in B(a_2)}\E_{z\in B(a_3)}&\mu_{\beta(b_{12})}(x,y)\mu_{\beta(b_{13})}(x,z)\mu_{\beta(b_{23})}(y,z)\\
\mu_{\beta(b_{13})}(x',z)&\mu_{\beta(b_{23})}(y',z)\mu_{\beta(b_{12})}(x',y)\\\nonumber
&\mu_{\beta(b_{23})}(y,z')\mu_{\beta(b_{12})}(x,y')\mu_{\beta(b_{13})}(x,z')
\end{align*}
satisfies 
\[g(x',y',z')=1+O(p^{3\ell+9q-\tau/8});\]
\item for all sets $I, J, K \subseteq [m]$ and for all but an $O_m(p^{10(\ell+q)(|I|+|J|+|K|)+2q(|I||J|+|I||K|+|J||K|)-\tau/4})$-proportion of tuples $(x,x',y,y',z,z')$ with $x,x' \in B(a_1)$, $y,y' \in B(a_2)$ and $z,z' \in B(a_3)$, the function $h$ defined by
\begin{align*}
h(x,x',y, y',z,z')=  &\E_{\begin{subarray}{l}  x_i\in B(a_1)\\i\in I\end{subarray}}\E_{\begin{subarray}{l}y_j\in B(a_2)\\j\in J\end{subarray}}\E_{\begin{subarray}{l} z_k\in B(a_3)\\k\in K \end{subarray}}\\
&\hspace{10pt}\prod_{i\in I, j\in J} \mu_{\beta(b_{12})}(x_i,y_j)\prod_{i\in I, k\in K}\mu_{\beta(b_{13})}(x_i,z_k)\prod_{j\in J, k\in K}\mu_{\beta(b_{23})}(y_j,z_k)\\
&\hspace{25pt}\prod_{k\in K}\mu_{\beta(b_{13})}(x,z_k)\mu_{\beta(b_{13})}(x',z_k)\mu_{\beta(b_{23})}(y,z_k)\mu_{\beta(b_{23})}(y',z_k)\\
&\hspace{40pt}\prod_{j\in J}\mu_{\beta(b_{12})}(x,y_j)\mu_{\beta(b_{12})}(x',y_j)\mu_{\beta(b_{23})}(y_j,z)\mu_{\beta(b_{23})}(y_j,z')\\
&\hspace{55pt}\prod_{i\in I}\mu_{\beta(b_{12})}(x_i,y)\mu_{\beta(b_{12})}(x_i,y')\mu_{\beta(b_{13})}(x_i,z)\mu_{\beta(b_{13})}(x_i,z')
\end{align*}
satisfies 
\[h(x,x',y, y',z,z')=1+O_m(p^{2(\ell+q)(|I|+|J|+|K|)+q(|I||J|+|I||K|+|J||K|)-\tau/8}).\]
\end{enumerate}
\end{lemma}

In the proof of Lemma \ref{lem:extuples} we will use the following elementary fact about averaging, which is a weighted variant of \cite[Lemma 6.5]{Gowers.2006}.

\begin{lemma}\label{fct:avg}
Let $\eta \in (0,1]$, $\delta\in (0,1/4)$, $\theta>0$, and let $f,g$ be real-valued functions on some finite domain $D$ such that $|\mathrm{supp}(f)|\geq\theta|D|$, $f(z)\geq (2\theta)^{-1}$ for all $z\in \mathrm{supp}(f)$, $\E_{z\in D}f(z)\leq 1+\eta$, $\E_{z\in D}f(z)g(z)\geq \delta-\eta$ and $\E_{z\in D}f(z)g(z)^2 \leq \delta^2+\eta$. Then for all but an $8\eta^{1/2}$-proportion of $z\in \mathrm{supp}(f)$, $|g(z)-\delta|<\eta^{1/4}$.
\end{lemma}

\begin{proof}We expand
\[\E_{z\in D}f(z)(g(z)-\delta)^2=\E_{z\in D}f(z)g(z)^2-2\delta \E_{z\in D}f(z)g(z)+\delta^2\E_{z\in D}f(z),\]
which, using the assumptions, is at most 
\[(\delta^2+\eta)-2\delta(\delta-\eta)+\delta^2(1+\eta)=(1+2\delta+\delta^2)\eta\leq 4\eta,\]
where the final inequality uses the fact that $\delta<1/4$. Now for each $z\in \mathrm{supp}(f)$, $f(z)\geq (2\theta)^{-1}$, so if more than an $8\eta^{1/2}$-proportion of $z\in \mathrm{supp}(f)$ satisfied $|g(z)-\delta|\geq \eta^{1/4}$, we would have
\[4\eta\geq \E_{z\in D}f(z)(g(z)-\delta)^2 = \frac{1}{|D|}\sum_{z\in \mathrm{supp}(f)} f(z) (g(z)-\delta)^2> \frac{8\eta^{1/2}|\mathrm{supp}(f)|}{|D|}(2\theta)^{-1}\eta^{1/2}, \]
a contradiction since $\theta\leq|\mathrm{supp}(f)|/|D|$.
\end{proof}
\vspace{2mm}
We are now able to give a proof of Lemma \ref{lem:extuples}.
\vspace{3mm}

\begin{proofof}{Lemma \ref{lem:extuples}}
To see (i), first note that by Lemma \ref{lem:genbilsums} with $|I|=|J|=|K|=1$
\[\E_{x'\in B(a_1), y'\in B(a_2), z'\in B(a_3)}\mu_{\beta(b_{12})}(x',y')\mu_{\beta(b_{13})}(x',z')\mu_{\beta(b_{23})}(y',z')=1+O(p^{3\ell+6q-\tau/2}),\]
and by the same lemma with $|I|=|J|=|K|=2$,
\[\E_{x'\in B(a_1), y'\in B(a_2), z'\in B(a_3)}\mu_{\beta(b_{12})}(x',y')\mu_{\beta(b_{13})}(x',z')\mu_{\beta(b_{23})}(y',z')g(x',y',z')=1+O(p^{6\ell+18q-\tau/2}).\]
Finally, by the same reasoning with $|I|=|J|=|K|=3$,
\[\E_{x'\in B(a_1), y'\in B(a_2), z'\in B(a_3)}\mu_{\beta(b_{12})}(x',y')\mu_{\beta(b_{13})}(x',z')\mu_{\beta(b_{23})}(y',z')g(x',y',z')^2=1+O(p^{9\ell+36q-\tau/2}).\]
The conditions of Lemma \ref{fct:avg} are therefore satisfied with 
\[f(x',y',z')=\mu_{\beta(b_{12})}(x',y')\mu_{\beta(b_{13})}(x',z')\mu_{\beta(b_{23})}(y',z')\]
defined on $D=B(a_1)\times B(a_2)\times B(a_3)$, $\delta=1$, $\eta= O(p^{9\ell+36q-\tau/2})$ and $\theta=p^{-3q}(1+O(p^{q-\tau}))$.
It follows that for all but an $ O(p^{5\ell+18q-\tau/4})$-proportion of $(x',y',z')\in \mathrm{supp}(f)$, $|g(x',y',z')-1|=O(p^{3\ell+9q-\tau/8})$.

The proof of (ii) is similar, except this time we will take the function $f$ defined on $D=B(a_1)^2\times B(a_2)^2 \times B(a_3)^2$ to be the constant function 1. Note that by Lemma \ref{lem:genbilsums},
\[\E_{x,x'\in B(a_1), y,y'\in B(a_2), z,z'\in B(a_3)}h(x,x',y,y',z,z')\]
and
\[\E_{x,x'\in B(a_1), y,y'\in B(a_2), z,z'\in B(a_3)}h(x,x',y,y',z,z')^2\]
both equal
\[1+O_m(p^{(\ell+q)(2(|I|+1+|J|+1+|K|+1))+4q((|I|+1)(|J|+1)+(|I|+1)(|K|+1)+(|J|+1)(|K|+1))-\tau/2}).\]
Thus Lemma \ref{fct:avg} with $f=1$, $\theta=\delta=1$, $g=h$, and 
\[\eta=O_m(p^{(\ell+q)(2(|I|+1+|J|+1+|K|+1))+4q((|I|+1)(|J|+1)+(|I|+1)(|K|+1)+(|J|+1)(|K|+1))-\tau/2})\] 
yields that for all but an $O_m(p^{10(\ell+q)(|I|+|J|+|K|)+2q(|I||J|+|I||K|+|J||K|)-\tau/4})$-proportion of tuples $(x,x',y,y',z,z')$ with $x,x' \in B(a_1)$, $y,y' \in B(a_2)$ and $z,z' \in B(a_3)$, the function $h$ satisfies 
\[h(x,x',y, y',z,z')=1+O_m(p^{2(\ell+q)(|I|+|J|+|K|)+q(|I||J|+|I||K|+|J||K|)-\tau/8}),\]
as claimed.
\end{proofof}

Finally, we prove another lemma in the same vein which will be used in the proof of Proposition \ref{prop:locsparseuni} in Appendix \ref{app:locsparseuni}.

\begin{lemma}\label{lem:bilsumssparse}Let $\calB=(\calL,\calQ)$ be a quadratic factor on $\F_p^n$ of complexity $(\ell,q)$ and rank at least $\tau$ consisting of  vectors $r_1,\dots, r_\ell$ and matrices $M_1,\dots,M_q$. Let $d=(a_1, a_2, a_3, b_{12}, b_{13}, b_{23})\in \F_p^{\ell+q}\times \F_p^{\ell+q}\times \F_p^{\ell+q}\times\F_p^q\times \F_p^q\times \F_p^q$. Then
\begin{enumerate}
\item[(i)] for all but at most an $O(p^{O(\ell+q)-\tau/4})$-proportion of of $y_0,y_1\in B(a_2)$ and $z_0,z_1\in B(a_3)$ with $(y_0,z_0),(y_0,z_1),(y_1,z_0),(y_1,z_1)\in \beta(b_{23})$, the function $g$ defined by 
\[g(y_0,y_1,z_0,z_1)=\E_{ x\in B(a_1)}\prod_{j\in \{0,1\}} \mu_{\beta(b_{12})}(x,y_j)\prod_{k\in \{0,1\}} \mu_{\beta(b_{13})}(x,z_k)\]
satisfies
\[g(y_0,y_1,z_0,z_1)=1+O(p^{O(\ell+q)-\tau/8});\]
\item[(ii)] for all but at most an $O(p^{O(\ell+q)-\tau/4})$-proportion of $x_0\in B(a_1)$, $y_0\in B(a_2)$ and $z_0\in B(a_3)$ with $(x_0,y_0)\in \beta(b_{12})$, $(x_0,z_0)\in \beta(b_{13})$, and $(y_0,z_0)\in \beta(b_{23})$, the function $h$ defined by 
\[h(x_0,y_0,z_0)=\E_{y_1\in B(a_2)}\E_{z_1\in B(a_3)}\mu_{\beta(b_{12})}(x_0,y_1)\mu_{\beta(b_{13})}(x_0,z_1)\prod_{(j,k)\in \{0,1\}^2\setminus\{(0,0)\}} \mu_{\beta(b_{23})}(y_j,z_k)\]
satisfies 
\[h(x_0,y_0,z_0)=1+O(p^{O(\ell+q)-\tau/8});\]
\item[(iii)] for all but at most an $O(p^{O(\ell+q)-\tau/4})$-proportion of $w\in B(\Sigma(d))$, we have 
\[\E_{y_0\in B(a_2)}\E_{z_0\in B(a_3)}\mu_{\beta(b_{12}+b_{23}+a_2')}(w,y_0)\mu_{\beta(b_{13}+b_{23}+a_3')}(w,z_0)\mu_{\beta(b_{23})}(y_0,z_0)=1+O(p^{O(\ell+q)-\tau/8}),\]
where $a_2'=(a_{2,\ell+1},a_{2,\ell+2},\dots,a_{2,\ell+q} )$ and  $a_3'=(a_{3,\ell+1},a_{3,\ell+2},\dots,a_{3,\ell+q} )$.
\end{enumerate}
\end{lemma}

\begin{proof}
To see part (i), apply Lemma \ref{lem:genbilsums} with $|I|=1$ and $|J|=|K|=2$ to get
\[\E_{\begin{subarray}{l}y_j\in B(a_2)\\j\in [2]\end{subarray}}\E_{\begin{subarray}{l} z_k\in B(a_3)\\k\in [2] \end{subarray}}\prod_{(j,k)\in \{0,1\}^2}\mu(y_j,z_k) \;\;\E_{x\in B(a_1)}\prod_{j\in \{0,1\}} \mu(x,y_j)\prod_{k\in \{0,1\}}\mu(x,z_k)=1+O(p^{5(\ell+q)+8q-\tau/2}).\]
But also note that the same lemma applied with $|I|=|J|=|K|=2$ yields
\[\E_{\begin{subarray}{l}y_j\in B(a_2)\\j\in [2]\end{subarray}}\E_{\begin{subarray}{l} z_k\in B(a_3)\\k\in [2] \end{subarray}}\prod_{(j,k)\in \{0,1\}^2}\mu(y_j,z_k) \left(\E_{x\in B(a_1)}\prod_{j\in \{0,1\}} \mu(x,y_j)\prod_{k\in \{0,1\}}\mu(x,z_k)\right)^2\]
equals $1+O(p^{6\ell+18q-\tau/2})$. 

The conditions of Lemma \ref{fct:avg} are therefore satisfied with 
\[f(y_0,y_1,z_0,z_1)=\mu(y_0,z_0)\mu(y_0,z_1)\mu(y_1,z_0)\mu(y_1,z_1)\]
defined on $D=B(a_2)^2\times B(a_3)^2$, $\delta=1$, $\eta= O(p^{6\ell+18q-\tau/2})$ and $\theta=p^{-4q}(1+O(p^{q-\tau}))$. It follows that for all but an $O(p^{3\ell+9q-\tau/4})$-proportion of $y_0,y_1\in B(a_2)$ and $z_0,z_1\in B(a_3)$ with $(y_0,z_0),(y_0,z_1),(y_1,z_0),(y_1,z_1)\in \beta(b_{23})$, $g(y_0,y_1,z_0,z_1)=1+O(p^{2\ell+5q-\tau/8})$.

The proof part (ii) is identical (and essentially Lemma \ref{lem:extuples} (i)). 

For (iii), we again apply Lemma \ref{lem:genbilsums} with the tuple of labels $(\Sigma(d), a_2,a_3,b_{12}+b_{23}+a_2',b_{13}+b_{23}+a_3',b_{23})$,  $f$ equal to the constant function 1 supported on $B(\Sigma(d))$, and $|I|=|J|=|K|=1$ and $|I|=1, |J|=|K|=2$, respectively, to arrive at the conclusion.
\end{proof}

\section{Proof of Proposition \ref{prop:ip2loccontrol}}\label{app:counting}

We are now in a position to prove Proposition \ref{prop:ip2loccontrol}, i.e. the approximate control of the local $\IP_2$-operator by the local $U^3$ semi-norm. The estimates from the preceding appendix will play a crucial role.

Since the sheer number of letters that appear in the course of the proof easily becomes overwhelming, we shall suppress subscripts on $\mu$ and $f$ as no ambiguity arises from doing so (e.g. variables labelled $x$ will always lie in $B(a_1)$ and be indexed by the letter $i$, and pairs $(x,y)$ will always be constrained to lie in $\beta(b_{12})$, etc.).  Thus, for example, when we write $\mu(x_i,y_j)$, we will mean $\mu_{\beta(b_{12})}(x_i,y_j)$, and when we write $f(x_i+y_j+z_S)$, we will mean $f_{i,j,S}( x_i+y_j+z_S)$.

\vspace{3mm}

\begin{proofof}{Proposition \ref{prop:ip2loccontrol}}
Let $I, J, \calS$ denote a choice of $i,j \in [m], S\subseteq [m]^2$ for which $\|f_{I,J,\calS}\|_{U^3(d)}=\min\{\|f_{i,j,S}\|_{U^3(d)}: i,j\in [m], \calS\subseteq [m]^2\}$ is minimal.  Then with the above notational conventions, $T_{m-\IP_2(d)}((f_{i,j,S})_{i,j\in [m],S\subseteq [m]^2})$ from Definition \ref{def:localip2op} can be expressed as
\begin{align}\label{eq:target}
\E_{\begin{subarray}{l}  x_i\in B(a_1)\\i\in [m]\end{subarray}}\E_{\begin{subarray}{l}y_j\in B(a_2)\\j\in [m]\end{subarray}}\E_{\begin{subarray}{l} z_S\in B(a_3)\\S\neq \calS \end{subarray}}&
\prod_{i, j\in [m]}\mu(x_i,y_j) \prod_{i\in [m], S\neq \calS} \mu(x_ i,z_S)\\\nonumber
\prod_{j\in [m], S\neq \calS}&\mu(y_j,z_S)\prod_{i,j\in [m], S\neq \calS}f(x_i+y_j+z_S)\\\nonumber
&\E_{ z_\calS\in B(a_3)}\;\prod_{i\in [m]}\mu(x_i,z_\calS)\prod_{j\in [m]}\mu(y_j, z_\calS)\prod_{i,j\in [m]} f(x_i+y_j+z_\calS).
\end{align}
Taking (\ref{eq:target}) to the 8th power, writing  
\[\prod_{i, j\in [m]}\mu(x_i,y_j) \prod_{i\in [m], S\neq \calS} \mu(x_ i,z_S)\prod_{j\in [m], S\neq \calS}\mu(y_j,z_S)\prod_{i,j\in [m], S\neq \calS}f(x_i+y_j+z_S)\]
as
\begin{align*}
\prod_{i, j\in [m]}\mu(x_i,y_j)^{1/2} \prod_{i\in [m], S\neq \calS} &\mu(x_ i,z_S)^{1/2}\prod_{j\in [m], S\neq \calS}\mu(y_j,z_S)^{1/2}\prod_{i,j\in [m], S\neq \calS}f(x_i+y_j+z_S)\\
&\cdot \prod_{i, j\in [m]}\mu(x_i,y_j)^{1/2} \prod_{i\in [m], S\neq \calS} \mu(x_ i,z_S)^{1/2}\prod_{j\in [m], S\neq \calS}\mu(y_j,z_S)^{1/2}
\end{align*}
and applying Cauchy-Schwarz, we obtain an upper bound on $|T_{m-\IP_2(d)}((f_{i,j,S})_{i,j\in [m],S\subseteq [m]^2})|^8$ of
\begin{align}\label{al:7}
\Big(\E_{\begin{subarray}{l}  x_i\in B(a_1)\\i\in [m]\end{subarray}}\E_{\begin{subarray}{l}y_j\in B(a_2)\\j\in [m]\end{subarray}}\E_{\begin{subarray}{l} z_S\in B(a_3)\\S\neq \calS \end{subarray}}
\prod_{i, j\in [m]}\mu(x_i,y_j)& \prod_{i\in [m], S\neq \calS} \mu(x_ i,z_S)\\\nonumber
&\prod_{j\in [m], S\neq \calS}\mu(y_j,z_S)\prod_{i,j\in [m], S\neq \calS}|f(x_i+y_j+z_S)|^2 \Big)^4
\end{align}
times
\begin{align}\label{al:7prime}
\Big(\E_{\begin{subarray}{l}  x_i\in B(a_1)\\i\in [m]\end{subarray}}\E_{\begin{subarray}{l}y_j\in B(a_2)\\j\in [m]\end{subarray}}&\E_{\begin{subarray}{l} z_S\in B(a_3)\\S\neq \calS \end{subarray}}
\prod_{i, j\in [m]}\mu(x_i,y_j) \prod_{i\in [m], S\neq \calS} \mu(x_ i,z_S)\prod_{j\in [m], S\neq \calS}\mu(y_j,z_S)\\\nonumber
&\big|\E_{z_\calS\in B(a_3)}\prod_{i\in [m]}\mu(x_i,z_\calS)\prod_{j\in [m]}\mu(y_j,z_\calS)\prod_{i,j\in [m]}f(x_i+y_j+z_\calS)\big|^2 \Big)^4.
\end{align}
Since the family of functions $f_{i,j,S}$ is 1-bounded and $\mu$ is non-negative, the first term in the product, (\ref{al:7}), is at most 
\[\Big(\E_{\begin{subarray}{l}  x_i\in B(a_1)\\i\in [m]\end{subarray}}\E_{\begin{subarray}{l}y_j\in B(a_2)\\j\in [m]\end{subarray}}\E_{\begin{subarray}{l} z_S\in B(a_3)\\S\neq \calS \end{subarray}}
\prod_{i, j\in [m]}\mu(x_i,y_j) \prod_{i\in [m], S\neq \calS} \mu(x_ i,z_S)\prod_{j\in [m], S\neq \calS}\mu(y_j,z_S) \Big)^4,\]
which by Lemma \ref{lem:genbilsums} with $|I|=|J|=m$ and $|K|=2^m-1$ equals
\[1+O_m(p^{(\ell+q)(2m+2^{m^2}-1)+q(m^2+2m(2^{m^2}-1))-\tau/2})=1+O_m(p^{4m2^{m^2}(\ell+q)-\tau/2}).\]
We expand the second term in the product, (\ref{al:7prime}), as
\begin{align*}
\Big(&\E_{\begin{subarray}{l}  x_i\in B(a_1)\\i\in [m]\end{subarray}}\E_{\begin{subarray}{l}y_j\in B(a_2)\\j\in [m]\end{subarray}}\E_{\begin{subarray}{l} z_S\in B(a_3)\\S\neq \calS \end{subarray}}
\prod_{i, j\in [m]}\mu(x_i,y_j) \prod_{i\in [m], S\neq \calS} \mu(x_ i,z_S)\prod_{j\in [m], S\neq \calS}\mu(y_j,z_S)\\\nonumber
&\E_{z_\calS,z_\calS'\in B(a_3)}\prod_{i\in [m]}\mu(x_i,z_\calS)\mu(x_i,z_\calS')\prod_{j\in [m]}\mu(y_j,z_\calS)\mu(y_j,z_\calS')\prod_{i,j\in [m]}f(x_i+y_j+z_\calS)\overline{f(x_i+y_j+z_\calS')} \Big)^4,
\end{align*}
which can be rearranged as
\begin{align*}
\Big(\E_{\begin{subarray}{l}  x_i\in B(a_1)\\i\in [m]\end{subarray}}&\E_{\begin{subarray}{l}y_j\in B(a_2)\\j\neq J\end{subarray}}\E_{\begin{subarray}{l} z_S\in B(a_3)\\S\neq \calS \end{subarray}}\E_{z_\calS,z_\calS'\in B(a_3)}
\prod_{i\in [m], j\neq J}\mu(x_i,y_j) \prod_{i\in [m], S\neq \calS} \mu(x_ i,z_S)\prod_{j\neq J, S\neq \calS}\mu(y_j,z_S)\\
\prod_{i\in [m]}&\mu(x_i,z_\calS)\mu(x_i,z_\calS')\prod_{j\neq J}\mu(y_j,z_\calS)\mu(y_j,z_\calS')\prod_{i\in [m],j\neq J}f(x_i+y_j+z_\calS)\overline{f(x_i+y_j+z_\calS')}\\
&\E_{y_J\in B(a_2)}\mu(y_J,z_\calS)\mu(y_J,z_\calS')\prod_{S\neq \calS} \mu(y_J,z_S)\prod_{i\in [m]} \mu(x_i,y_J)f(x_i+y_J+z_\calS)\overline{f(x_i+y_J+z_\calS')}\Big)^4.
\end{align*}
We again regard 
\begin{align*}\prod_{i\in [m], j\neq J}\mu(x_i,y_j)& \prod_{i\in [m], S\neq \calS} \mu(x_ i,z_S)\prod_{j\neq J, S\neq \calS}\mu(y_j,z_S)\\
&\prod_{i\in [m]}\mu(x_i,z_\calS)\mu(x_i,z_\calS')\prod_{j\neq J}\mu(y_j,z_\calS)\mu(y_j,z_\calS')\prod_{i\in [m],j\neq J}f(x_i+y_j+z_\calS)\overline{f(x_i+y_j+z_\calS')}
\end{align*} 
as a product of two factors, each containing the square root of the product of instances of $\mu$, and one containing the product over instances of $f$. This allows us to apply Cauchy-Schwarz again to obtain an upper bound on (\ref{al:7prime}) of
\begin{align}\label{al:8}
\Big(\E_{\begin{subarray}{l}  x_i\in B(a_1)\\i\in [m]\end{subarray}}&\E_{\begin{subarray}{l}y_j\in B(a_2)\\j\neq J\end{subarray}}\E_{\begin{subarray}{l} z_S\in B(a_3)\\S\neq \calS \end{subarray}}
\E_{z_\calS,z_\calS'\in B(a_3)}
\prod_{i\in [m], j\neq J}\mu(x_i,y_j) \prod_{i\in [m], S\neq \calS} \mu(x_ i,z_S)\prod_{j\neq J, S\neq \calS}\mu(y_j,z_S)\\\nonumber
&\prod_{i\in [m]}\mu(x_i,z_\calS)\mu(x_i,z_\calS')\prod_{j\neq J}\mu(y_j,z_\calS)\mu(y_j,z_\calS')\prod_{i\in [m],j\neq J}|f(x_i+y_j+z_\calS)\overline{f(x_i+y_j+z_\calS')}|^2\Big)^2
\end{align}
times
\begin{align}\label{al:8prime}
\Big(&\E_{\begin{subarray}{l}  x_i\in B(a_1)\\i\in [m]\end{subarray}}\E_{\begin{subarray}{l}y_j\in B(a_2)\\j\neq J\end{subarray}}\E_{\begin{subarray}{l} z_S\in B(a_3)\\S\neq \calS \end{subarray}}
\E_{z_\calS,z_\calS'\in B(a_3)}\\\nonumber
&\prod_{i\in [m], j\neq J}\mu(x_i,y_j) \prod_{i\in [m], S\neq \calS} \mu(x_ i,z_S)\prod_{j\neq J, S\neq \calS}\mu(y_j,z_S)\prod_{i\in [m]}\mu(x_i,z_\calS)\mu(x_i,z_\calS')\prod_{j\neq J}\mu(y_j,z_\calS)\mu(y_j,z_\calS')\\\nonumber
&\big|\E_{y_J\in B(a_2)}\mu(y_J,z_\calS)\mu(y_J,z_\calS')\prod_{S\neq \calS} \mu(y_J,z_S)\prod_{i\in [m]} \mu(x_i,y_J)f(x_i+y_J+z_\calS)\overline{f(x_i+y_J+z_\calS')}\big|^2\Big)^2.
\end{align}
As before, the first term in this product, (\ref{al:8}), is at most 
\begin{align*}
\Big(\E_{\begin{subarray}{l}  x_i\in B(a_1)\\i\in [m]\end{subarray}}\E_{\begin{subarray}{l}y_j\in B(a_2)\\j\neq J\end{subarray}}\E_{\begin{subarray}{l} z_S\in B(a_3)\\S\neq \calS \end{subarray}}
\E_{z_\calS,z_\calS'\in B(a_3)}
\prod_{i\in [m], j\neq J}&\mu(x_i,y_j) \prod_{i\in [m], S\neq \calS} \mu(x_ i,z_S)\prod_{j\neq J, S\neq \calS}\mu(y_j,z_S)\\\nonumber
&\prod_{i\in [m]}\mu(x_i,z_\calS)\mu(x_i,z_\calS')\prod_{j\neq J}\mu(y_j,z_\calS)\mu(y_j,z_\calS')\Big)^2,
\end{align*}
which equals $1+O_m(p^{(\ell+q)(2m+2^{m^2})+q(m(m-1)+(2m-1)(2^{m^2}+1))-\tau/2})=1+O_m(p^{4m2^{m^2}(\ell+q)-\tau/2})$, by an application of Lemma \ref{lem:genbilsums} with $|I|=m$, $|J|=m-1$ and $|K|=2^{m^2}+1$.

The second factor, (\ref{al:8prime}), can be expanded and rearranged to give
\begin{align*}
\Big(\E_{\begin{subarray}{l}  x_i\in B(a_1)\\i\in [m]\end{subarray}}\E_{\begin{subarray}{l}y_j\in B(a_2)\\j\neq J\end{subarray}}\E_{\begin{subarray}{l} z_S\in B(a_3)\\S\neq \calS \end{subarray}}
\E_{z_\calS,z_\calS'\in B(a_3)} &\prod_{i\in [m], j\neq J}\mu(x_i,y_j) \prod_{i\in [m], S\neq \calS} \mu(x_ i,z_S)\prod_{j\neq J, S\neq \calS}\mu(y_j,z_S) \\
&\hspace{20pt}\prod_{i\in [m]}\mu(x_i,z_\calS)\mu(x_i,z_\calS')\prod_{j\neq J}\mu(y_j,z_\calS)\mu(y_j,z_\calS')\\
\E_{y_J,y_J'\in B(a_2)}\mu(y_J,z_\calS)\mu(y_J,z_\calS')&\mu(y_J',z_\calS)\mu(y_J',z_\calS')\prod_{S\neq \calS} \mu(y_J,z_S)\mu(y_J',z_S)\\
\prod_{i\in [m]} \mu(x_i,y_J) \mu(x_i,y_J')f(x_i&+y_J+z_\calS)\overline{f(x_i+y_J+z_\calS')f(x_i+y_J'+z_\calS)}f(x_i+y_J'+z_\calS')\Big)^2,
\end{align*}
or 
\begin{align*}
\Big(\E_{\begin{subarray}{l}  x_i\in B(a_1)\\i\neq I\end{subarray}}&\E_{\begin{subarray}{l}y_j\in B(a_2)\\j\neq J\end{subarray}}\E_{\begin{subarray}{l} z_S\in B(a_3)\\S\neq \calS \end{subarray}}
\E_{z_\calS,z_\calS'\in B(a_3)}\E_{y_J,y_J'\in B(a_2)}\mu(y_J,z_\calS)\mu(y_J,z_\calS')\mu(y_J',z_\calS)\mu(y_J',z_\calS')  \\
\prod_{i\neq I, j\neq J} &\mu(x_i,y_j) \prod_{i\neq I, S\neq \calS} \mu(x_i,z_S)\prod_{j\neq J, S\neq \calS} \mu(y_j,z_S) \\
\prod_{j\neq J}&\mu(y_j,z_\calS)\mu(y_j,z_\calS')\prod_{S\neq \calS}\mu(y_J,z_S)\mu(y_J',z_S)\prod_{i\neq I}\mu(x_i,y_J)\mu(x_i,y_J')\mu(x_i,z_\calS)\mu(x_i,z_\calS')\\
&\prod_{i\neq I}f(x_i+y_J+z_\calS)\overline{f(x_i+y_J+z_\calS')f(x_i+y_J'+z_\calS)}f(x_i+y_J'+z_\calS')\\
&\hspace{20pt}\E_{x_I\in B(a_1)}\mu(x_I,y_J) \mu(x_I,y_J')\mu(x_I,z_\calS) \mu(x_I,z_\calS')\prod_{j\neq J}\mu(x_I,y_j)\prod_{S\neq \calS}\mu(x_I,z_S)\\
&\hspace{100pt}f(x_I+y_J+z_\calS)\overline{f(x_I+y_J+z_\calS')f(x_I+y_J'+z_\calS)}f(x_I+y_J'+z_\calS')\Big)^2.
\end{align*}
For the final time, we regard the combined product of instances of $\mu$ and $f$ in the first, second, third and fourth line as an appropriate product involving square-roots of $\mu$, and apply Cauchy-Schwarz to obtain, via Lemma \ref{lem:genbilsums}, an upper bound of $1+O_m(p^{4m2^{m^2}(\ell+q)-\tau/2})$ times
\begin{align*}
\E_{\begin{subarray}{l}  x_i\in B(a_1)\\i\neq I\end{subarray}}&\E_{\begin{subarray}{l}y_j\in B(a_2)\\j\neq J\end{subarray}}\E_{\begin{subarray}{l} z_S\in B(a_3)\\S\neq \calS \end{subarray}}
\E_{z_\calS,z_\calS'\in B(a_3)}\E_{y_J,y_J'\in B(a_2)}\mu(y_J,z_\calS)\mu(y_J,z_\calS')\mu(y_J',z_\calS)\mu(y_J',z_\calS')  \\
\prod_{i\neq I, j\neq J} &\mu(x_i,y_j) \prod_{i\neq I, S\neq \calS} \mu(x_i,z_S)\prod_{j\neq J, S\neq \calS} \mu(y_j,z_S) \\
\prod_{j\neq J}&\mu(y_j,z_\calS)\mu(y_j,z_\calS')\prod_{S\neq \calS}\mu(y_J,z_S)\mu(y_J',z_S)\prod_{i\neq I}\mu(x_i,y_J)\mu(x_i,y_J')\mu(x_i,z_\calS)\mu(x_i,z_\calS')\\
&|\E_{x_I\in B(a_1)}\mu(x_I,y_J) \mu(x_I,y_J') \mu(x_I,z_\calS) \mu(x_I,z_\calS')\prod_{j\neq J}\mu(x_I,y_j)\prod_{S\neq \calS}\mu(x_I,z_S)\\
&\hspace{40pt}f(x_I+y_J+z_\calS)\overline{f(x_I+y_J+z_\calS')f(x_I+y_J'+z_\calS)}f(x_I+y_J'+z_\calS')|^2.
\end{align*}
Upon expanding, this expression equals the inner product
\begin{equation}\label{eq:innerprod}
\langle h,F\rangle:=\E_{x_I,x_I'\in B(a_1)}\E_{y_J,y_J'\in B(a_2)}\E_{z_\calS,z_\calS'\in B(a_3)} h(x_I,x_I',y_J, y_J',z_\calS,z_\calS')\cdot F(x_I,x_I',y_J, y_J',z_\calS,z_\calS'),
\end{equation}
where
\begin{align*}
h(x_I,x_I',y_J, y_J',z_\calS,z_\calS')&=  \E_{\begin{subarray}{l}  x_i\in B(a_1)\\i\neq I\end{subarray}}\E_{\begin{subarray}{l}y_j\in B(a_2)\\j\neq J\end{subarray}}\E_{\begin{subarray}{l} z_S\in B(a_3)\\S\neq \calS \end{subarray}}\prod_{i\neq I, j\neq J} \mu(x_i,y_j)\prod_{i\neq I, S\neq \calS}\mu(x_i,z_S)\prod_{ j\neq J, S\neq \calS}\mu(y_j,z_S)\\
\prod_{S\neq \calS}\mu(x_I,z_S)&\mu(x_I',z_S)\mu(y_J,z_S)\mu(y_J',z_S)
\prod_{j\neq J}\mu(x_I,y_j)\mu(x_I',y_j)\mu(y_j,z_\calS)\mu(y_j,z_\calS')\\
&\prod_{i\neq I}\mu(x_i,y_J)\mu(x_i,y_J')\mu(x_i,z_\calS)\mu(x_i,z_\calS')
\end{align*}
and 
\begin{align*}
F(x_I,x_I',y_J, y_J',& z_\calS ,z_\calS')= \mu(x_I,y_J)\mu(x_I',y_J)\mu(x_I,y_J')\mu(x_I',y_J')  \\
\mu(x_I,z_\calS)&\mu(x_I',z_\calS)\mu(x_I,z_\calS')\mu(x_I',z_\calS')  \mu(y_J,z_\calS)\mu(y_J',z_\calS)\mu(y_J,z_\calS')\mu(y_J',z_\calS')  \\
f_{I,J,\calS}&(x_I+y_J+z_\calS)\overline{f_{I,J,\calS}(x_I+y_J+z_\calS')}\overline{f_{I,J,\calS}(x_I+y_J'+z_\calS)}f_{I,J,\calS}(x_I+y_J'+z_\calS')\\
&\overline{f_{I,J,\calS}(x_I'+y_J+z_\calS)}f_{I,J,\calS}(x_I'+y_J+z_\calS')f_{I,J,\calS}(x_I'+y_J'+z_\calS)\overline{f_{I,J,\calS}(x_I'+y_J'+z_\calS')}.
\end{align*}
We conclude the proof of Proposition \ref{prop:ip2loccontrol} by showing that the inner product in (\ref{eq:innerprod}) is approximately equal to 
\[\E_{x_I,x_I'\in B(a_1)}\E_{y_J,y_J'\in B(a_2)}\E_{z_\calS,z_\calS'\in B(a_3)} F(x_I,x_I',y_J, y_J',z_\calS,z_\calS'),\]
which equals $\min_{i,j\in [m],S\subseteq [m]^2}\|f_{i,j,S}\|_{U^3(d)}^8$ by our choice of $I,J,\calS$. In doing so, we shall follow the argument in \cite[Lemma 6.7]{Gowers.2006}. 

First, note that Lemma \ref{lem:extuples} (ii) with $|I|=|J|=m-1$ and $|K|=2^{m^2}-1$ states that for all but an $O_m(p^{24m2^{m^2}(\ell+q)-\tau/4})$-proportion of choices of $x_I,x_I'\in B(a_1)$, $y_J,y_J'\in B(a_2)$, $z_\calS,z_\calS'\in B(a_3)$, we have that
\begin{equation}\label{eq:goodsix}
h(x_I,x_I',y_J, y_J',z_\calS,z_\calS')=1+O_m(p^{6m2^{m^2}(\ell+q)-\tau/8}).
\end{equation}
Call those six-tuples $x_I,x_I'\in B(a_1)$, $y_J,y_J'\in B(a_2)$, $z_\calS,z_\calS'\in B(a_3)$ for which (\ref{eq:goodsix}) holds `good', and define a function $h'$ by setting $h'(x_I,x_I',y_J, y_J',z_\calS,z_\calS')=h(x_I,x_I',y_J, y_J',z_\calS,z_\calS')$ for all good six-tuples, and $h'(x_I,x_I',y_J, y_J',z_\calS,z_\calS')=1$ otherwise. By definition of $h'$ and (\ref{eq:goodsix}),
\[\|h'-1\|_\infty = O_m(p^{6m2^{m^2}(\ell+q)-\tau/8}).\]
Moreover, since all but an $O_m(p^{24m2^{m^2}(\ell+q)-\tau/4})$-proportion of $x_I,x_I'\in B(a_1)$, $y_J,y_J'\in B(a_2)$, $z_\calS,z_\calS'\in B(a_3)$ are good,
\[\|h-h'\|_1 =\E_{x_I,x_I'\in B(a_1)}\E_{y_J,y_J'\in B(a_2)}\E_{z_\calS,z_\calS'\in B(a_3)} |h(x_I,x_I',y_J, y_J',z_\calS,z_\calS')-h'(x_I,x_I',y_J, y_J',z_\calS,z_\calS')|\]
is at most $O_m(p^{24m2^{m^2}(\ell+q)-\tau/4})$.
Moreover,  by Lemma \ref{lem:bilsums} (iii), $\|\mu\|_{\infty}= p^q(1+O(p^{q-\tau}))$, so since $f_{I,J,\calS}$ is 1-bounded we have that $\|F\|_\infty\leq p^{12q}(1+O(p^{q-\tau}))$. By Lemma \ref{lem:genbilsums} with $I=J=K=[2]$, we also have the bound
\[\|F\|_1\leq 1+O(p^{6\ell+18q-\tau/2}).\] 
By linearity of the inner product and the triangle inequality then, 
\[|\langle h,F\rangle-\langle 1,F\rangle|\leq |\langle h-h',F\rangle| + |\langle h'- 1,F\rangle|\leq \|h-h'\|_1\|F\|_\infty + \|h'-1\|_\infty \|F\|_1,\]
which by the above observations is at most $O_m(p^{32m2^{m^2}(\ell+q)-\tau/8})$. 

We have thus shown that the 8th power of (\ref{eq:target}) is bounded above in absolute value by 
\[(1+O_m(p^{4m2^{m^2}(\ell+q)-\tau/2}))|\langle h,F\rangle|=(1+O_m(p^{4m2^{m^2}(\ell+q)-\tau/2}))(|\langle 1,F\rangle| +O_m(p^{32m2^{m^2}(\ell+q)-\tau/8})).\]
The observation that $\langle 1,F\rangle= \|f_{I,J,\calS}\|_{U^3(d)}^8$
now finishes the proof of Proposition \ref{prop:ip2loccontrol}.
\end{proofof}

\section{Proof of Proposition \ref{prop:locsparseuni}}\label{app:locsparseuni}

The proof of Proposition \ref{prop:locsparseuni} below makes crucial use of Lemma \ref{lem:bilsumssparse} from the preceding appendix.

\begin{proofof}{Proposition \ref{prop:locsparseuni}} Fix $\e>0$. We shall choose the growth function $\rho_0$ at the end of the proof. Let $\rho\geq \rho_0$ be a growth function, and let $\calB=(\calL,\calQ)$ be a quadratic factor on $\F_p^n$ of complexity $(\ell,q)$ and rank at least $\rho(\ell+q)$, and let $d=(a_1,a_2,a_3,b_{12},b_{13},b_{23})\in \mathbb{F}_p^{3(\ell+q)+3q}$. 

Let $A\subseteq \F_p^n$, and without loss of generality assume that $\alpha_{B(\Sigma(d))}\in [0,\e)$ (else apply the argument below to the complement of $A$). Expanding out 
\[\|1_A-\alpha_{B(\Sigma(d))}\|_{U^{3}(d)}^8\]
as a local Gowers inner product, we find that it equals $\|1_A\|_{U^{3}(d)}^8$ plus $2^8-1$ terms in which at least one of the eight inputs equals the constant function $\alpha_{B(\Sigma(d))}$. By the triangle inequality, our assumption on $\alpha_{B(\Sigma(d))}$, and  Lemma \ref{lem:genbilsums} with $|I|=|J|=|K|=2$, each of these terms can therefore be bounded above by
\begin{equation}\label{eq:initialobs}
\e(1+O(p^{6\ell+18q-\rho/2})).
\end{equation}
Now by definition, 
\begin{align*}\|1_A\|_{U^3(d)}^8=\E_{ x_0,x_1\in B(a_1)}\E_{y_0,y_1\in B(a_2)}&\E_{z_0,z_1\in B(a_3)}\\
\prod_{(i,j)\in \{0,1\}^2} &\mu_{\beta(b_{12})}(x_i,y_j)\prod_{(i,k)\in \{0,1\}^2} \mu_{\beta(b_{13})}(x_i,z_k)\prod_{(j,k)\in \{0,1\}^2} \mu_{\beta(b_{23})}(y_j,z_k)\\
&\prod_{\omega\in \{0,1\}^3}1_A(x_{\omega(1)}+y_{\omega(2)}+z_{\omega(3)}),
\end{align*}
which is trivially at most 
\begin{align*}
\E_{ x_0,x_1\in B(a_1)}&\E_{y_0,y_1\in B(a_2)}\E_{z_0,z_1\in B(a_3)}\\
&\prod_{(i,j)\in \{0,1\}^2}  \mu_{\beta(b_{12})}(x_i,y_j)\prod_{(i,k)\in \{0,1\}^2} \mu_{\beta(b_{13})}(x_i,z_k)\prod_{(j,k)\in \{0,1\}^2} \mu_{\beta(b_{23})}(y_j,z_k)1_A(x_{0}+y_{0}+z_{0}).
\end{align*}
Rearranging, we find that $\|1_A\|_{U^3(d)}^8$ is at most 
\begin{align*}
\E_{x_0\in B(a_1)}\E_{y_0\in B(a_2)}&\E_{z_0\in B(a_3)}1_A(x_0+y_0+z_0)\mu_{\beta(b_{12})}(x_0,y_0)\mu_{\beta(b_{13})}(x_0,z_0)\mu_{\beta(b_{23})}(y_0,z_0)\\
\E_{y_1\in B(a_2)}&\E_{z_1\in B(a_3)}\mu_{\beta(b_{12})}(x_0,y_1)\mu_{\beta(b_{13})}(x_0,z_1)\prod_{(j,k)\in \{0,1\}^2\setminus\{(0,0)\}} \mu_{\beta(b_{23})}(y_j,z_k)\\
&\E_{ x_1\in B(a_1)}\prod_{j\in \{0,1\}} \mu_{\beta(b_{12})}(x_1,y_j)\prod_{k\in \{0,1\}} \mu_{\beta(b_{13})}(x_1,z_k)
\end{align*}
By Lemma \ref{lem:bilsumssparse} (i), for all but an $O(p^{O(\ell+q)-\rho(\ell+q)/4})$-proportion of $y_0,y_1\in B(a_2)$ and $z_0,z_1\in B(a_3)$ with $(y_0,z_0),(y_0,z_1),(y_1,z_0),(y_1,z_1)\in \beta(b_{23})$, the expectation in $x_1$ in the final line is $1+O(p^{O(\ell+q)-\rho(\ell+q)/8})$.
Similarly,  by Lemma \ref{lem:bilsumssparse} (ii), for all but an $O(p^{O(\ell+q)-\rho(\ell+q)/4})$-proportion of $x_0\in B(a_1)$, $y_0\in B(a_2)$ and $z_0\in B(a_2)$ with $x_0\in B(a_1)$, $y_0\in B(a_2)$ and $z_0\in B(a_3)$ with $(x_0,y_0)\in \beta(b_{12})$, $(x_0,y_0)\in \beta(b_{13})$, and $(y_0,z_0)\in \beta(b_{23})$, we have that
\[\E_{y_1\in B(a_2)}\E_{z_1\in B(a_3)}\mu_{\beta(b_{12})}(x_0,y_1)\mu_{\beta(b_{13})}(x_0,z_1)\prod_{(j,k)\in \{0,1\}^2\setminus\{(0,0)\}} \mu_{\beta(b_{23})}(y_j,z_k)\]
equals $1+O(p^{O(\ell+q)-\rho(\ell+q)/8})$. Combining with the preceding estimate, we find that  $\|1_A\|_{U^3(d)}^8$ is at most 
\begin{align*}(1+O(p^{O(\ell+q)-\rho(\ell+q)/8}))\E_{x_0\in B(a_1)}&\E_{y_0\in B(a_2)}\E_{z_0\in B(a_3)}1_A(x_0+y_0+z_0)\\
&\mu_{\beta(b_{12})}(x_0,y_0)\mu_{\beta(b_{13})}(x_0,z_0)\mu_{\beta(b_{23})}(y_0,z_0)+O(p^{O(\ell+q)-\rho(\ell+q)/4}).
\end{align*}
Note that if $x_0\in B(a_1)$, $y_0\in B(a_2)$, $z_0\in B(a_3)$, $(x_0,y_0)\in \beta(b_{12})$, $(x_0,z_0)\in \beta(b_{13})$, and $(y_0,z_0)\in \beta(b_{23})$, then $x_{0}+y_{0}+z_{0}\in B(\Sigma(d))$. Conversely, if $y_0\in B(a_2)$, $z_0\in B(a_3)$, $(x_0,y_0)\in \beta(b_{12})$, $(x_0,z_0)\in \beta(b_{13})$, $(y_0,z_0)\in \beta(b_{23})$, and $x_{0}+y_{0}+z_{0}\in B(\Sigma(d))$, then $x_0\in B(a_1)$. We may therefore make a change of variable $w=x_0+y_0+z_0$ and conclude, noting that all atoms of $\calB$ are approximately the same size by Lemma \ref{lem:sizeofatoms}, that
\begin{align*}
&\E_{x_0\in B(a_1)}\E_{y_0\in B(a_2)}\E_{z_0\in B(a_3)}1_A(x_0+y_0+z_0)\mu_{\beta(b_{12})}(x_0,y_0)\mu_{\beta(b_{13})}(x_0,z_0)\mu_{\beta(b_{23})}(y_0,z_0)\\
&=(1+O(p^{\ell+q-\rho(\ell+q)/2}))\E_{w\in B(\Sigma(d))}\E_{y_0\in B(a_2)}\E_{z_0\in B(a_3)}1_A(w)\mu_{\beta(b_{12})}(w-y_0-z_0,y_0)\\
&\hspace{280pt}\mu_{\beta(b_{13})}(w-y_0-z_0,z_0)\mu_{\beta(b_{23})}(y_0,z_0)\\
&=(1+O(p^{\ell+q-\rho(\ell+q)/2}))\E_{w\in B(\Sigma(d))}1_A(w)\E_{y_0\in B(a_2)}\E_{z_0\in B(a_3)}\mu_{\beta(b_{12}+b_{23}+a_2')}(w,y_0)\\
&\hspace{280pt}\mu_{\beta(b_{13}+b_{23}+a_3')}(w,z_0)\mu_{\beta(b_{23})}(y_0,z_0).
\end{align*}
Finally, by Lemma \ref{lem:bilsumssparse} (iii), for all but an $O(p^{O(\ell+q)-\rho(\ell+q)/4})$-proportion of $w\in B(\Sigma(d))$, we have 
\[\E_{y_0\in B(a_2)}\E_{z_0\in B(a_3)}\mu_{\beta(b_{12}+b_{23}+a_2')}(w,y_0)\mu_{\beta(b_{13}+b_{23}+a_3')}(w,z_0)\mu_{\beta(b_{23})}(y_0,z_0)=1+O(p^{O(\ell+q)-\rho(\ell+q)/8}),\]
and thus
\[\|1_A\|_{U^3(d)}^8\leq (1+O(p^{O(\ell+q)-\rho(\ell+q)/8}))\E_{w\in B(\Sigma(d))}1_A(w)+O(p^{O(\ell+q)-\rho(\ell+q)/8}).\]
Now provided that $\rho_0$ is chosen such that both $O(p^{O(\ell+q)-\rho(\ell+q)/8}))$ in the preceding line and $O(p^{6\ell+18q-\rho/2})$ from (\ref{eq:initialobs}) are at most $\e$ for any $\rho\geq \rho_0$, we have that $\|1_A\|_{U^3(d)}^8\leq (1+\e)\alpha_{B(\Sigma(d))}+\e \leq 3\e$, and thus
\[\|1_A-\alpha_{B(\Sigma(d))}\|_{U^{3}(d)}^8\leq 3\e + (2^8-1)\e(1+\e)\leq 2^{10}\e.\]
This completes the proof on taking 8th roots.
\end{proofof}

For reference in \cite{Terry.2021a} and \cite{Terry.2024f}, we remark that in the above proof we have implicitly shown the following.

\begin{proposition}
Let $A\subseteq \F_p^n$ and let $\calB=(\calL,\calQ)$ be a quadratic factor on $\F_p^n$ of complexity $(\ell,q)$ and rank at least $\rho$. Then for any $d=(a_1,a_2,a_3,b_{12},b_{13},b_{23})\in  \F_p^{\ell+q}\times \F_p^{\ell+q}\times \F_p^{\ell+q}\times\F_p^q\times \F_p^q\times \F_p^q$, we have
\begin{align*}
\E_{x\in B(a_1)}\E_{y\in B(a_2)}\E_{z\in B(a_3)}1_A(x+y+z)\mu_{\beta(b_{12})}(x,y)&\mu_{\beta(b_{13})}(x,z)\mu_{\beta(b_{23})}(y,z)
\\&=\alpha_{B(\Sigma(d))}+O(p^{O(\ell+q)-\rho/8}),
\end{align*}
where $\alpha_{B(\Sigma(d))}$ denotes the density of $A$ on the atom $B(\Sigma(d))$.
\end{proposition}

\section{Counting across distinct atoms}\label{app:countinggen}

In this appendix we prove a counting lemma for ternary sums arising from distinct quadratic atoms, Proposition \ref{prop:countinggp}. This will play a crucial role in the revised version of our companion paper \cite{Terry.2021a}. Whilst the proof in earlier versions of \cite{Terry.2021a} (see \cite{Terry.2021av2}) passed through hypergraph regularity, the new proof we give here has been substantially simplified by the use of the local Gowers $U^3$ semi-norms developed in this paper. Because it is very similar in spirit to others in earlier sections, we have chosen to include it here rather than in \cite{Terry.2021a}. We could, of course, have set up the machinery at this level of generality earlier in the paper, but it seemed prudent not to overload the reader from the outset.

We will begin by generalising the counting lemma relative to linear atoms associated to the local Gowers $U^2$ semi-norm.  This serves as a warm up for the higher order analogues needed for Proposition \ref{prop:countinggp}, the main result of this section.

The ``multi-local $F$ operator" defined below is a generalization of the local IP-operator in Definition \ref{def:localipop} to arbitrary bipartite graphs $F$.

\begin{definition}[Multi-local $F$-operator]\label{def:multlocalipop}
Let $F=(U\cup V, E)$ be a bipartite graph.  Given a linear factor $\calL$ on $\F_p^n$ of complexity $\ell$, a tuple of labels $d=(d_u)_{u\in U}(d_v)_{v\in V}\in \F_p^{|U|\ell+|V|\ell}$,  and functions $(f_{u,v})_{u\in U, v\in V}:\F_p^n\ra \C$, define

\[T_{F(d)}((f_{u,v})_{u\in U, v\in V})=\E_{\begin{subarray}{l} x_u:u\in U\\ x_u\in L(d_u)\end{subarray}}\E_{\begin{subarray}{l}y_v:v\in V \\ y_v\in L(d_v) \end{subarray}}\prod_{u\in U, v\in V}f_{u,v}(x_u+y_v).\]
\end{definition}

To ease notation, we set the following convention for special cases of Definition \ref{def:multlocalipop}.  

\begin{notation}\label{not:ff}
In the notation of Definition \ref{def:multlocalipop}, given functions $f,g$ with the property that for all $uv\in E$, $f_{u,v}=f$ and for all $uv\notin E$, $f_{u,v}=g$, we write $T_{F(d)}(f|g)$ to mean $T_{F(d)}((f_{u,v})_{u\in U, v\in V})$.
\end{notation}

The operator in Definition \ref{def:multlocalipop} allows us to count  instances of a general bipartite graph $F$ in sum-graphs associated to a set $A$.  Indeed, given a bipartite graph  $F=(U\cup V, E)$, a set $A\subseteq \F_p^n$, a linear factor $\calL$ on $\F_p^n$ of complexity $\ell$, and a tuple of labels $d=(d_u)_{u\in U}(d_v)_{v\in V}\in \F_p^{|U|\ell+|V|\ell}$, observe that
$$
|\{(a_u)_{u\in U}(b_v)_{v\in V}\in \prod_{u\in U}L(d_u)\times \prod_{v\in V}L(d_v): a_u+b_v\in A\text{ if and only if }uv\in E\}|
$$
is equal to  
$$
T_{F(d)}(1_A| 1_{A^C})\prod_{u\in U}|L(d_u)|\prod_{v\in V}|L(d_v)|=T_{F(\dbar)}(1_A|1_{A^C})|L(0)|^{|U|+|V|},
$$
where $A^C=\F_p^n\setminus A$ and we have used the fact that all atoms of $\calL$ have the same size.

Our first lemma states that the multilocal $F$-operator is controlled by the local $U^2$ semi-norm. Its proof is identical to that of Lemma \ref{lem:iploccontrol}, so we leave it as an exercise to the interested reader. 

\begin{lemma}[Multi-local $F$-operator is controlled by local $U^2$]\label{lem:ipmultloccontrol}
Let $F=(U\cup V, E)$ be a bipartite graph.  Let $\ell\geq 2$, let $\calL$ be a linear factor on $\F_p^n$ of complexity $\ell$, and let $d=(d_u)_{u\in U}(d_v)_{v\in V}\in \F_p^{|U|\ell+|V|\ell}$ be a tuple of labels. Suppose that for each $u\in U$ and $v\in V$, $f_{u,v}:\F_p^n\ra \C$ is such that $\|f_{u,v}\|_{\infty}\leq 1$. Then 
\[|T_{F(d)}((f_{u,v})_{u\in U, v\in V})|\leq \min_{u\in U, v\in V}\|f_{u,v}\|_{U^2((d_u,d_v))}.\]
\end{lemma}

We immediately deduce the following counting lemma for general binary sums across possibly distinct atoms.

\begin{proposition}[Counting lemma for induced binary sums across linear atoms]\label{prop:countinggplinear}
For all $\e>0$ and bipartite graphs $F=(U\cup V,E)$ the following holds.  Let $A\subseteq \mathbb{F}_p^n$,  let $\calL$ be a linear factor of complexity $\ell$ on $\F_p^n$, and let $d=(d_u)_{u\in U}(d_v)_{v\in V}\in \F_p^{|U|\ell+|V|\ell}$ be a tuple of labels. For each $(u,v)\in U\times V$, define 
$$
\alpha_{uv}=\frac{|A\cap L(d_u+d_v)|}{|L(d_u+d_v))|}.
$$
Suppose that for each $(u,v)\in U\times V$, $\|1_A-\alpha_{uv}\|_{U^2((d_u,d_v))}\leq \e$.  Then the number of $(x_u)_{u\in U}(y_v)_{v\in V}\in \prod_{u\in U}L(d_u)\times \prod_{v\in V}L(d_v)$ such that $x_u+y_v\in A$ if and only if $uv\in E$ is within $\e |U||V|\prod_{u\in U}|L(d_u)| \prod_{v\in V}|L(d_v)|$ of 
\[\prod_{u\in U, v\in V: uv\in E}\alpha_{uv}\prod_{u\in U, v\in V: uv\notin E} (1-\alpha_{uv})\prod_{u\in U}|L(d_u)| \prod_{v\in V}|L(d_v)|.\]
\end{proposition}
\begin{proof}
For each $u\in U$ and $v\in V$  let $g_{A}^{uv}=1_A-\alpha_{uv}$, and $g_{A^C}^{uv}=1_{A^C}-\overline{\alpha}_{uv}$, where $A^C=G\setminus A$ and $\overline{\alpha}_{uv}$ is the density of $A^C$ on $L(d_u+d_v)$.  Note that $g_{A^C}^{uv}=-g_{A}^{uv}$.
Recalling Notation \ref{not:ff}, we observe that 
\begin{align}\label{fe1}
T_{F(d)}(1_A |1_{A^C})=T_{F(d)}((g_{A}^{uv}+\alpha_{uv})_{u\in U, v\in V: uv\in E} |(g_{A^C}^{uv}+\overline{\alpha}_{uv})_{u\in U, v\in V: uv\notin E}).
\end{align}
By linearity, (\ref{fe1}) is the sum of the term  
\begin{align}\label{al:11}\prod_{u\in U, v\in V: uv\in E}\alpha_{uv}\prod_{u\in U, v\in V: uv\notin E} (1-\alpha_{uv})
\end{align}
plus $2^{|U||V|}-1$ terms of the form
\begin{align}\label{al:22} T_{F(d)}((f_{u,v})_{u\in U, v\in V}),
\end{align}
in which at least one of the input functions $f_{u,v}$ equals $g_A^{uv}$ or $-g_A^{uv}=g_{A^C}^{uv}$. By Lemma \ref{lem:ipmultloccontrol}, each term of the form (\ref{al:22}) is at most $\|g_A^{uv}\|_{U^2((d_u,d_v))}$ for some choice of $u\in U$ and $v\in V$.   Now by assumption, each of these terms satisfies $\|g_A^{uv}\|_{U^2((d_u,d_v))}<\e$. Consequently, $T_{F(d)}(1_A |1_{A^C})$ is within $(2^{|U||V|}-1)\e$ of (\ref{al:11}). This implies the desired conclusion by definition of $T_{F(d)}(1_A |1_{A^C})$ (see the remark following Notation \ref{not:ff}).
\end{proof}

We next define an operator that allows us to count copies of a given $3$-partite $3$-uniform hypergraph in a ternary sum-graph.  This definition is based on the $\IP_2$-operator introduced in Definition \ref{def:localip2op}, but differs in that it is designed to help us count copies of an arbitrary 3-partite 3-uniform hypergraph $F$ where the vertices lie in pre-specified but possibly distinct atoms.  

\begin{definition}[Ternary multi-local $F$-operator]\label{def:multiF}
Let $F=(U\cup V\cup W,E)$ be a $3$-partite $3$-uniform hypergraph, and let $\calB=(\calL,\calQ)$ be a quadratic factor on $\F_p^n$ of complexity $(\ell,q)$.  Assume that for each $u\in U, v\in V, w\in W$, we have a tuple 
$$
e_{uvw}=(a_u,b_v,c_w,d_{uv},d_{uw},d_{vw})\in \F_p^{\ell+q}\times \F_p^{\ell+q}\times \F_p^{\ell+q}\times \F_p^q\times \F_p^q\times \F_p^q.
$$
 Given $e=(e_{uvw})_{u\in U, v\in V, w\in W}$ and functions $(f_{u,v,w})_{u\in U, v\in V, w\in W}: \F_p^n \ra\C$, define
\begin{align*}
T_{F(e)}((f_{u,v,w})_{u\in U, v\in V, w\in W})=\E_{x_u\in B(a_u):u\in U}\E_{y_v\in B(b_v): v\in V}&\E_{z_w\in B(c_w): w\in W}\prod_{u\in U, v\in V} \mu_{\beta(d_{uv})}(x_u,y_v)\\
\prod_{u\in U, w\in W}&\mu_{\beta(d_{uw})}(x_u,z_w)\prod_{ v\in V, w\in W}\mu_{\beta(d_{vw})}(y_v,z_w)\\
&\prod_{u\in U, v\in V, w\in W}f_{u,v,w}(x_u+y_v+z_w).
\end{align*}
\end{definition}

We note that $T_{F(e)}$ is linear in its input functions.  As before, we will use the following simplified notation for special cases of Definition \ref{def:multiF}.  

\begin{notation}
In the notation of Definition \ref{def:multiF}, given functions $f,g$ with the property that for all $uvw\in E$, $f_{u,v,w}$ is equal to $f$ and for all $uvw\notin E$, $f_{u,v,w}=g$, we will write $T_{F(e)}(f|g)$ to denote $T_{F(e)}((f_{uvw})_{u\in U,v\in V,w\in W})$.
\end{notation}

In analogy to the linear case, the operator in Definition \ref{def:multiF} allows us to count copies of $F$ with vertices behaving in pre-specified ways with respect to the factor $\calB$.  We now describe the relevant type of counting problem in more detail.  

\begin{definition}\label{def:if}
Let $F=(U\cup V\cup W,E)$ be a $3$-partite $3$-uniform hypergraph, let  $\calB=(\calL,\calQ)$ be a quadratic factor on $\F_p^n$ of complexity $(\ell,q)$, and let $A\subseteq \F_p^n$.  Assume  that for each $(u,v,w)\in U\times V
\times W$,  we have a tuple 
$$
e_{uvw}=(a_u,b_v,c_w,d_{uv},d_{uw},d_{vw})\in \F_p^{\ell+q}\times \F_p^{\ell+q}\times \F_p^{\ell+q}\times \F_p^q\times \F_p^q\times \F_p^q.
$$
 Given $e=(e_{uvw})_{u\in U, v\in V, w\in W }$, define
\begin{align*}
\calI_F(e)=\{(x_u)_{u\in U}(y_v)_{v\in V}(z_w)_{w\in W}\in  \prod_{u\in U}B(a_u)&\times\prod_{v\in V}B(b_v)\times\prod_{w\in W}B(c_w):\\
 \text{ for each }&(u,v,w)\in U\times V\times W,
(x_u,y_v)\in \beta_{\calQ}(d_{uv}), \\
&(x_u,z_w)\in \beta_{\calQ}(d_{uw}), \text{ and }(y_v,z_w)\in \beta_{\calQ}(d_{vw})\}.
\end{align*}
\end{definition}

When the factor has sufficiently high rank, the sets of the form $\calI_F(e)$ are always non-empty. This can be deduced by combining Lemma \ref{lem:genbilsums} with Corollaries \ref{cor:sizeofatoms} and \ref{cor:sizeofbilsets}.

\begin{lemma}\label{lem:Isize}
Let $F=(U\cup V\cup W,E)$ be a $3$-partite $3$-uniform hypergraph. There is a growth function $\tau$ such that the following holds. 

Let  $\calB=(\calL,\calQ)$ be a quadratic factor on $\F_p^n$ of complexity $(\ell,q)$ and rank at least $\tau(\ell+q)$. Assume that for each $(u,v,w)\in U\times V
\times W$,  we have a tuple of labels
$$
e_{uvw}=(a_u,b_v,c_w,d_{uv},d_{uw},d_{vw})\in \F_p^{\ell+q}\times \F_p^{\ell+q}\times \F_p^{\ell+q}\times \F_p^q\times \F_p^q\times \F_p^q ,
$$
and let $e=(e_{uvw})_{u\in U, v\in V, w\in W }$. Then $|\calI_F(e)|>0$.
\end{lemma}

In the notation of Definition \ref{def:if}, we are interested in understanding what proportion of tuples $(x_u)_{u\in U}(y_v)_{v\in V}(z_w)_{w\in W}\in \calI_F(e)$ satisfy $ x_u+y_v+z_w\in A$ if and only if $uvw\in E$, for a fixed subset $A\subseteq \F_p^n$.  

\begin{proposition}\label{prop:opmeaning}
Let $F=(U\cup V\cup W,E)$ be a $3$-partite $3$-uniform hypergraph, let $\calB=(\calL,\calQ)$ be a quadratic factor on $\F_p^n$ of complexity $(\ell,q)$ and rank $\tau$, and let $A\subseteq \F_p^n$.  Assume  that for each $(u,v,w)\in U\times V
\times W$,  we have a tuple 
$$
e_{uvw}=(a_u,b_v,c_w,d_{uv},d_{uw},d_{vw})\in \F_p^{\ell+q}\times \F_p^{\ell+q}\times \F_p^{\ell+q}\times \F_p^q\times \F_p^q\times \F_p^q .
$$ 
Then given $e=(e_{uvw})_{u\in U, v\in V, w\in W }$, the number of $(x_u)_{u\in U}(y_v)_{v\in V}(z_w)_{w\in W}\in \calI_F(e)$ such that $ x_u+y_v+z_w\in A$ if and only if $uvw\in E$ is equal to
$$
(1+O(p^{3(|U|+|V|+|W|)\ell+3q(|U||V|+|U||W|+|V||W|)-\tau}))|T_{F(e)}(1_A|1_{A^C})||\calI_F(e)|.
$$
\end{proposition}
\begin{proof}
This follows from the definition of $T_{F(e)}(1_A|1_{A^C})$ and Lemma \ref{lem:extuples}  
\end{proof}

We will need the following quadratic analogue of Lemma \ref{lem:ipmultloccontrol}. The astute reader will notice that the proof of Lemma \ref{lem:ip2multcontrol} is essentially identical to the proof Proposition \ref{prop:ip2loccontrol}. It is thus omitted. 

\begin{lemma}[Multi-local $F$-operator is controlled by local $U^3$]\label{lem:ip2multcontrol}
Let $m\geq 2$, let $U,V,W\subseteq [m]$ and let $F=(U\cup V\cup W,E)$ be a $3$-partite $3$-uniform hypergraph.
Let $\calB=(\calL,\calQ)$ be a quadratic factor on $\F_p^n$  of complexity $(\ell,q)$ and rank at least $\tau$. 
For each $(u,v,w)\in U\times V\times W$, let $e_{uvw}=(a_u,b_v,c_w,d_{uv},d_{uw},d_{vw})\in\F_p^{\ell+q}\times \F_p^{\ell+q}\times \F_p^{\ell+q}\times \F_p^q\times \F_p^q\times \F_p^q$ be a tuple of labels, and set $e=(e_{uvw})_{u\in U, v\in V, w\in W }$. 

Suppose that for each $u\in U$, $v\in V$, and $w\in W$, $f_{u,v,w}:\F_p^n\ra \C$ is such that $\|f_{u,v,w}\|_{\infty}\leq 1$. Then 
\[|T_{F(e)}((f_{u,v,w})_{u\in U, v\in V, w\in W})|\leq (1+O(p^{O_m(\ell+q)-\tau/2})) \min_{u\in U, v\in V, w\in W}\|f_{u,v,w}\|_{U^3(e_{uvw})}+O(p^{O_m(\ell+q)-\tau/128}).\]
\end{lemma}

We now state and prove the ternary counting lemma we shall need.   

\begin{proposition}[Counting lemma for induced ternary sums across quadratic atoms]\label{prop:countinggp}
For every positive integer $m$ and all $\e>0$, there exists a growth function $\tau_0=\tau_0(m,\e)$ such that the following holds for any growth function $\tau\geq \tau_0$. 

Let $U,V,W\subseteq [m]$, and let $F=(U\cup V\cup W,E_F)$ be a $3$-partite $3$-uniform hypergraph. Let $\calB=(\calL,\calQ)$ be a quadratic factor on $\F_p^n$ of complexity $(\ell,q)$ and rank at least $\tau$, and let $A\subseteq \mathbb{F}_p^n$. For each $(u,v,w)\in U\times V
\times W$,  let 
$$
e_{uvw}=(a_u,b_v,c_w,d_{uv},d_{uw},d_{vw})\in \F_p^{\ell+q}\times \F_p^{\ell+q}\times \F_p^{\ell+q}\times \F_p^q\times \F_p^q\times \F_p^q 
$$
 be a tuple of labels, and let $e=(e_{uvw})_{u\in U, v\in V, w\in W }$.

Suppose that for each $(u,v,w)\in U\times V\times W$, $\|1_A-\alpha_{uvw}\|_{U^3(e_{uvw})}\leq \e$, where
$$
\alpha_{uvw}=\frac{|A\cap B(\Sigma(e_{uvw}))|}{|B(\Sigma(e_{uvw}))|}.
$$  
Then we have
$$
\Big|T_{F(e)}(1_A| 1_{A^C})-\prod_{u\in U, v\in V, w\in W: uvw\in E_F}\alpha_{uvw}\prod_{u\in U, v\in V, w\in W: uvw\notin E_F} (1-\alpha_{uvw}) \Big|<3\e m^3.
$$
Moreover, the number of tuples $(x_u)_{u\in U}(y_v)_{v\in V}(z_w)_{w\in W}\in \calI_F(e)$ such that $ x_u+y_v+z_w\in A$ if and only if $uvw\in E_F$ differs from 
\[\prod_{u\in U, v\in V, w\in W: uvw\in E_F}\alpha_{uvw}\prod_{u\in U, v\in V, w\in W: uvw\notin E_F} (1-\alpha_{uvw}) |\calI_F(e)|\]
by at most $O_m(\e)|\calI_F(e)|$.
\end{proposition}

\begin{proof} Given $m\geq 1$, let the growth function $\tau_0$ be such that for all $\tau\geq \tau_0$, the error term $O(p^{O_m(\ell+q)-\tau(\ell+q)/2})$ arising from Lemma \ref{lem:ip2multcontrol} is at most $\e$. For each $u\in U$, $v\in V$ and $w\in W$, let $g_{A}^{e_{uvw}}=1_A-\alpha_{uvw}$, and $g_{A^C}^{e_{uvw}}=1_{A^C}-\overline{\alpha}_{uvw}$, where $A^C=G\setminus A$ and $\overline{\alpha}_{uvw}=1-\alpha_{uvw}$ is the density of $A^C$ on $B(\Sigma(e_{uvw}))$.  Note $g_{A^C}^{e_{uvw}}=-g_{A}^{e_{uvw}}$.
Recalling the notation following Definition \ref{def:multiF}, consider 
\[T_{F(e)}(1_A |1_{A^C})=T_{F(e)}((g_{A}^{e_{uvw}}+\alpha_{uvw})_{u\in U, v\in V, w\in W: uvw\in E_F} |(g_{A^C}^{e_{uvw}}+\overline{\alpha}_{uvw})_{u\in U, v\in V, w\in W: uvw\notin E_F}).\]
By linearity, this equals the term 
\begin{align}\label{al:1}\prod_{u\in U, v\in V, w\in W: uvw\in E_F}\alpha_{ijk}\prod_{u\in U, v\in V, w\in W: uvw\notin E_F} (1-\alpha_{uvw})
\end{align}
plus $|U||V||W|-1$ terms of the form
\begin{align}\label{al:2}T_{F(e)}((f_{u,v,w})_{u\in U, v\in V, w\in W}),
\end{align}
in which at least one of the input functions $f_{u,v,w}$ equals $g_A^{e_{uvw}}$ or $-g_A^{e_{uvw}}=g_{A^C}^{e_{uvw}}$. By definition of $\tau_0$, by Lemma \ref{lem:ip2multcontrol} that each term of the form (\ref{al:2}) is bounded above by $2\|g_A^{e_{uvw}}\|_{U^3(e_{uvw})}+\e$, for some choice of $u\in U$, $v\in V$, and $w\in W$. Now by assumption, for each $u\in U$, $v\in V$, and $w\in W$, $\|g_A^{e_{uvw}}\|_{U^3(e_{uvw})}<\e$. This shows that in absolute value,
$$
T_{F(e)}(1_A| 1_{A^C})-\prod_{u\in U, v\in V, w\in W: uvw\in E_F}\alpha_{uvw}\prod_{u\in U, v\in V, w\in W: uvw\notin E_F} (1-\alpha_{uvw})
$$
is bounded above by $3\e m^3$. 

The `moreover' part follows immediately from Proposition \ref{prop:opmeaning} on collecting error terms.\end{proof}

In particular, we shall use the following corollary in \cite{Terry.2021b}.

\begin{corollary}\label{cor:countinggp}
For every positive integer $m$, there exists a constant $C=C(m)$ and for all $\e>0$, there exists a growth function $\tau_0=\tau_0(m,\e)$ such that the following holds for any growth function $\tau\geq \tau_0$. 

Let $U,V,W\subseteq [m]$, and let $F=(U\cup V\cup W,E_F)$ be a $3$-partite $3$-uniform hypergraph. Let $\calB=(\calL,\calQ)$ be a quadratic factor on $\F_p^n$ of complexity $(\ell,q)$ and rank at least $\tau$, and let $A\subseteq \mathbb{F}_p^n$. For each $(u,v,w)\in U\times V
\times W$,  let 
$$
e_{uvw}=(a_u,b_v,c_w,d_{uv},d_{uw},d_{vw})\in \F_p^{\ell+q}\times \F_p^{\ell+q}\times \F_p^{\ell+q}\times \F_p^q\times \F_p^q\times \F_p^q 
$$
 be a tuple of labels, and let $e=(e_{uvw})_{u\in U, v\in V, w\in W }$.

Suppose that for each $(u,v,w)\in U\times V\times W$, $\|1_A-\alpha_{uvw}\|_{U^3(e_{uvw})}\leq C \e^{m^2 2^{m^2}}$, where
$$
\alpha_{uvw}=\frac{|A\cap B(\Sigma(e_{uvw}))|}{|B(\Sigma(e_{uvw}))|}\in (\e,1-\e).
$$ 
Then there exists a tuple $(x_u)_{u\in U}(y_v)_{v\in V}(z_w)_{w\in W}\in \calI_F(e)$ such that $ x_u+y_v+z_w\in A$ if and only if $uvw\in E_F$.
\end{corollary}

\begin{proof}
This follows immediately from Proposition \ref{prop:countinggp} combined with Lemma \ref{lem:Isize}.
\end{proof}

%\bibliography{VC2bibi.bib}

\begin{thebibliography}{10}

\bibitem{Alon.2007}
Noga Alon, Eldar Fischer, and Ilan Newman, \emph{{Efficient testing of
  bipartite graphs for forbidden induced subgraphs}}, SIAM Journal on Computing
  \textbf{37} (2007), no.~3, 959--976.

\bibitem{Alon.2018is}
Noga Alon, Jacob Fox, and Yufei Zhao, \emph{{Efficient arithmetic regularity
  and removal lemmas for induced bipartite patterns}}, Discrete Analysis
(2019), no.~3, 14pp.

\bibitem{Baldwin.1976}
John Baldwin and Jan Saxl, \emph{{Logical stability in group theory}}, Journal
  of the Australian Mathematical Society \textbf{21} (1976), no.~3, 267--276.

\bibitem{Bhattacharyya.20126s9}
Arnab Bhattacharyya, Eldar Fischer, Hamed Hatami, Pooya Hatami, and Shachar
  Lovett, \emph{{Every locally characterized affine-invariant property is
  testable}}, Proceedings of the 45th ACM Symposium on Theory of Computing
  (STOC 2013) (2013), 429--436.

\bibitem{Bloom.2020}
Thomas Bloom, \emph{{Quantitative inverse theory of Gowers uniformity norms
  [after F. Manners]}}, S\'eminaire Bourbaki \textbf{73} (2020), no.~1174, 1--31.

\bibitem{Candela.2019}
Pablo Candela and Bal\'azs Szegedy, \emph{{Regularity and inverse theorems for
  uniformity norms on compact abelian groups and nilmanifolds}},
Journal f\"ur die Reine und Angewandte Mathematik (Crelle's Journal) \textbf{789} (2022), 1--42.

\bibitem{Chernikov.2019b}
Artem Chernikov and Nadja Hempel, \emph{{On n-dependent groups and fields II}},
  Forum of Mathematics, Sigma \textbf{E38} (2021), no.~9, 1--51.

\bibitem{Chernikov.2019}
Artem Chernikov, Daniel Palacin, and Kota Takeuchi, \emph{{On n-dependence}},
  Notre Dame Journal of Formal Logic \textbf{60} (2019), no.~2, 195--214.

\bibitem{Chernikov.2020}
Artem Chernikov and Henry Towsner, \emph{{Hypergraph regularity and higher
  arity VC-dimension}}, arXiv:2010.00726 (2020).

\bibitem{Conant.2021}
Gabriel Conant, \emph{{Quantitative structure of stable sets in arbitrary
  finite groups}}, Proceedings of the American Mathematical Society
  \textbf{149} (2021), no.~9, 4015--4028.

\bibitem{Conant.2018n5}
Gabriel Conant and Anand Pillay, \emph{{Pseudofinite groups and VC-dimension}},
  Journal of Mathematical Logic \textbf{21} (2021), no.~2, 2150009, 23pp.

\bibitem{Conant.2017p4}
Gabriel Conant, Anand Pillay, and Caroline Terry, \emph{{A group version of
  stable regularity}}, Mathematical Proceedings of the Cambridge Philosophical
  Society \textbf{168} (2020), 405--413.

\bibitem{Conant.2018zd}
\bysame, \emph{{Structure and regularity for subsets of groups with finite
  $\VC$-dimension}}, Journal of the European Mathematical Society \textbf{24} (2022), no. 2, pp. 583--621.
  
  \bibitem{Gladkova.2025}
  Val Gladkova, \emph{{A note on lower bounds for arithmetic regularity partitions}}, arXiv: (2025).

\bibitem{Gowers.1998}
Timothy Gowers, \emph{{A new proof of Szemer\'edi's theorem for arithmetic
  progressions of length four}}, Geometric and Functional Analysis \textbf{8}
  (1998), no.~3, 529--551.

\bibitem{Gowers.2001}
\bysame, \emph{{A new proof of Szemer\'edi's theorem}}, Geometric and Functional
  Analysis \textbf{11} (2001), no.~3, 465--588.
  
  \bibitem{Gowers.2006}
\bysame, \emph{{Quasirandomness, counting and regularity for 3-uniform hypergraphs}}, Combinatorics, Probability and Computing \textbf{15} (2006), no.~3, 143--184.

\bibitem{Gowers.2016s0e}
\bysame, \emph{{Generalizations of Fourier analysis, and how to apply them}},
  Bulletin of the American Mathematical Society \textbf{54} (2016), no.~1,
  1--44.
  
\bibitem{Gowers.2024}
Timothy Gowers, Ben Green, Freddie Manners, and Terence Tao, \emph{{Marton's
  conjecture in abelian groups with bounded torsion}}, arXiv:2404.02244 (2024).

\bibitem{Gowers.2010tc}
Timothy Gowers and Julia Wolf, \emph{{The true complexity of a system of linear
  equations}}, Proceedings of the London Mathematical Society
  \textbf{100} (2010), no.~1, 155--176.

\bibitem{Gowers.2011}
\bysame, \emph{{Linear forms and higher-degree uniformity for functions on
  $\F_p^n$}}, Geometric and Functional Analysis \textbf{21} (2011), no.~1, 36--69.

\bibitem{Green.2005s}
Ben Green, \emph{{Finite field models in additive combinatorics}}, Surveys in Combinatorics 2005, London Math. Soc. Lecture Notes \textbf{327}, 1--27.

\bibitem{Green.2005}
Ben Green, \emph{{A Szemer\'edi-type regularity lemma in abelian groups, with
  applications}}, Geometric and Functional Analysis \textbf{15} (2005), no.~2,
  340--376.

\bibitem{Green.2007}
\bysame, \emph{{Montr\'eal notes on quadratic Fourier analysis}}, Additive
  combinatorics (Montr\'eal 2006, ed. Granville et al.), CRM Proceedings \textbf{43}, 69--102, AMS, 2007.
  
  \bibitem{Green.2015}
  Ben Green and Tom Sanders, \emph{Fourier uniformity on subspaces}, arXiv:1607.07701 (2015).

\bibitem{Green.2008piq}
Ben Green and Terence Tao, \emph{{An inverse theorem for the Gowers
  $U^3(G)$ norm}}, Proceedings of the Edinburgh Mathematical
  Society \textbf{51} (2008), no.~1, 73--153.

\bibitem{Green.2008}
\bysame, \emph{{The primes contain arbitrarily long arithmetic progressions}},
  Annals of Mathematics \textbf{167} (2008), no.~2, 481--547.

\bibitem{Green.2010jwua}
\bysame, \emph{{An arithmetic regularity lemma, an associated counting lemma,
  and applications}}, An Irregular Mind, Bolyai Soc. Math. Stud., no.~21, Janos
  Bolyai Math. Soc., Budapest, 2010, pp.~261--334.

\bibitem{Green.2010d5}
\bysame, \emph{{An equivalence between inverse sumset theorems and inverse
  conjectures for the $U^3$ norm}}, Mathematical Proceedings of the Cambridge Philosophical Society \textbf{149} (2010), no.~1, 1--19.

\bibitem{Green.2010ta4}
\bysame, \emph{{Linear equations in primes}}, Annals of Mathematics \textbf{171} (2010), no.~3, 1753--1850.

\bibitem{Green.2011r7c}
Ben Green, Terence Tao, and Tamar Ziegler, \emph{{An inverse theorem for the
  Gowers $U^{s+1}[N]$ norm}}, Electronic Research
  Announcements in Mathematical Sciences \textbf{18} (2011), no.~0, 69--90.

\bibitem{Green.2012xge}
\bysame, \emph{{An inverse theorem for the Gowers
  $U^{s+1}[N]$ norm}}, Annals of Mathematics
  \textbf{176} (2012), no.~2, 1231--1372.

\bibitem{Green.2011}
\bysame, \emph{{An inverse theorem for the
  Gowers $U^4$ norm}}, Glasgow Mathematical Journal
  \textbf{53} (2011), no.~1, 1--50.

\bibitem{Hempel.2016}
Nadja Hempel, \emph{{On n-dependent groups and fields}}, Mathematical Logic
  Quarterly \textbf{62} (2016), no.~3, 215--224.

\bibitem{Host.2005}
Bernard Host and Bryna Kra, \emph{{Nonconventional ergodic averages and
  nilmanifolds}}, Annals of Mathematics \textbf{161} (2005),
  no.~1, 397--488.

\bibitem{Kra.2006}
Bryna Kra, \emph{{From combinatorics to ergodic theory and back again}},
  Proceedings of International Congress of Mathematicians, Madrid 2006
  \textbf{III} (2006), 57--76.

\bibitem{Lovasz.2010}
Laszlo Lovasz and Balazs Szegedy, \emph{{Regularity partitions and the topology
  of graphons}}, An Irregular Mind, Bolyai Soc. Math. Stud., no.~21, Janos
  Bolyai Math. Soc., Budapest, 2010, pp.~415--446.

\bibitem{Malliaris.2014}
Maryanthe Malliaris and Saharon Shelah, \emph{{Regularity lemmas for stable
  graphs}}, Transactions of the American Mathematical Society \textbf{366}
  (2014), no.~3, 1551--1585.

\bibitem{Manners.2018}
Frederick Manners, \emph{{Quantitative bounds in the inverse theorem for the
  Gowers $U^{s+1}$ norms over cyclic groups}},
  arXiv:1811.00718 (2018).

\bibitem{Milicevic.2021}
Luka Mili\'cevi\'c, \emph{{An Inverse Theorem for Certain Directional Gowers
  Uniformity Norms}}, Publications de l'Institut Math\'ematique \textbf{113} (2023),
  1--56.

\bibitem{Milicevic.2024}
\bysame, \emph{{Good bounds for sets lacking skew corners}}, arXiv:2404.07180 (2024).

\bibitem{Peluse.2024}
Sarah Peluse, \emph{{Finite field models in additive combinatorics--twenty years on}}, Surveys in Combinatorics 2024, London Math. Soc. Lecture Notes \textbf{493}, 159--199.

\bibitem{Roth.1953}
K.~F. Roth, \emph{{On certain sets of integers}}, Journal of the London
  Mathematical Society. Second Series \textbf{28} (1953), no.~1, 104--109.

\bibitem{Samorodnitsky.2006}
Alex Samorodnitsky, \emph{{Low-degree tests at large distances}}, The 39th
  Annual ACM Symposium on Theory of Computing (STOC 2007) (2007), 506--515.

\bibitem{Sanders.2018}
Tom Sanders, \emph{{The coset and stability rings}}, Online Journal of Analytic
  Combinatorics \textbf{15} (2020), no.~2020, 10pp.

\bibitem{Shelah.2014}
Saharon Shelah, \emph{{Strongly dependent theories}}, Israel Journal of
  Mathematics \textbf{204} (2014), no.~1, 1--83.

\bibitem{Shelah.2007}
\bysame, \emph{{Definable groups for dependent and 2-dependent theories}},
  Sarajevo Journal of Mathematics \textbf{25} (2017), no.~13, 3--25.

\bibitem{Sisask.2018}
Olof Sisask, \emph{{Convolutions of sets with bounded VC-dimension are
  uniformly continuous}}, Discrete Analysis (2021), no.~1, 25pp.

\bibitem{Szemeredi.1975}
Endre Szemer\'edi, \emph{{On sets of integers containing no $k$ elements in
  arithmetic progression}}, Acta Arithmetica, no.~27 (1975), 199--245.

\bibitem{Tao.2012cv}
Terence Tao, \emph{{Higher order Fourier analysis}}, Graduate Studies in
  Mathematics, vol. 142, American Mathematical Society, Providence, RI, 2012.

\bibitem{Tao.2010}
Terence Tao and Van Vu, \emph{{Additive combinatorics}}, Cambridge University
  Press, vol. 105, Cambridge University Press, 2010.

\bibitem{Tao.2016}
Terence Tao and Tamar Ziegler, \emph{{Polynomial patterns in the primes}},
 Forum of Mathematics, Pi (2018), no.~6:e1, 60pp. 

\bibitem{Terry.2018}
Caroline Terry, \emph{{$\VC_l$-dimension and the jump to the fastest speed of a
  hereditary L-property}}, Proc. Amer. Math. Soc. \textbf{146} (2018),
  3111--3126.

\bibitem{Terry.2019}
Caroline Terry and Julia Wolf, \emph{{Stable arithmetic regularity in the
  finite field model}}, Bulletin of the London Mathematical Society \textbf{51}
  (2019), 70--88.

\bibitem{Terry.2021a}
\bysame, \emph{{Higher-order generalizations of stability and arithmetic
  regularity}}, arXiv:2111.01739 (for the earlier version referenced in the text, see entry below)
  (2025).
  
  \bibitem{Terry.2021av2}
\bysame, \emph{{Higher-order generalizations of stability and arithmetic
  regularity}}, arXiv:2111.01739v2 (for the updated version of this manuscript, see entry above)
  (2021).

\bibitem{Terry.2021b}
\bysame, \emph{{Irregular triads in 3-uniform hypergraphs}}, To appear, Memoirs of the American Mathematical Society, arXiv:2111.01737
  (2021).

\bibitem{Terry.2024f}
\bysame, \emph{{On the quadratic complexity of subsets of $\F_p^n$ of bounded $\mathrm{VC_{2}}$-dimension}}, arXiv:2510.12767 (2025).

\bibitem{Wolf.2015}
Julia Wolf, \emph{{Finite field models in arithmetic combinatorics -- ten years
  on}}, Finite Fields and their Applications \textbf{32} (2015), 233--274.

\end{thebibliography}
%\bibliographystyle{amsplain}

\providecommand{\bysame}{\leavevmode\hbox to3em{\hrulefill}\thinspace}
\providecommand{\MR}{\relax\ifhmode\unskip\space\fi MR }
\providecommand{\MRhref}[2]{%
  \href{http://www.ams.org/mathscinet-getitem?mr=#1}{#2}
}
\providecommand{\href}[2]{#2}

\end{document}